\documentclass[11pt,letterpaper]{amsart}

\usepackage[T1]{fontenc}
\usepackage[utf8]{inputenc}   % keep if you use pdfLaTeX
\usepackage{textcomp}         % extra glyphs in T1
\usepackage{palatino}         % URW Palladio (ppl)
\usepackage[euler-hat-accent]{eulervm}  % Euler math
\usepackage{tikz-cd}
\usepackage{hyperref}
\usepackage{relsize}
\usepackage{silence}                       % load very early, before the culprit
\usepackage{amsmath,amssymb,amsthm,amscd}
\usepackage{stmaryrd}
\usepackage{xspace,verbatim}
\usepackage{latexsym}
\usepackage{microtype}
\usepackage{enumitem}
\usepackage{epsfig}    % legacy; consider \usepackage{graphicx}
\usepackage{xcolor}
\usepackage{url}
\usepackage{mathtools}

\setlist[enumerate]{itemsep=4pt, topsep=4pt} % more vertical space between items & before/after list
\setlist[itemize]{itemsep=4pt, topsep=4pt}
\usepackage{xypic,epic,eepic,floatflt}
\input xy
\xyoption{all}

\newlength{\baseunit}
\newcount{\numlines}
\numberwithin{equation}{section}
\newtheorem{proposition}{Proposition}[section]
\newtheorem{theorem}[proposition]{Theorem}
\newtheorem{corollary}[proposition]{Corollary}
\newtheorem{lemma}[proposition]{Lemma}

\newtheorem*{claim}{Claim}

\theoremstyle{definition}
\newtheorem{definition}[proposition]{Definition}
\newtheorem{notation}[proposition]{Notation}
\newtheorem{example}[proposition]{Example}

\theoremstyle{remark}
\newtheorem{remark}[proposition]{Remark}

\makeatletter
\newcommand{\rightcomment}[1]{%
  \begingroup
  \tagsleft@false
  \tag*{\makebox[0pt][r]{\text{#1}}}%
  \endgroup
}
\makeatother

\newcommand\codim{\operatorname{codim}}
\newcommand\Bl{\operatorname{Bl}}

\title[The Chow ring of a wonderful compactification]{The Chow ring of a wonderful compactification}

\author[P. Gallardo]{Patricio Gallardo}
\address{Department of Mathematics, University of California, Riverside, Riverside, CA 92521, United States.}
\email{pgallard@ucr.edu}

\author[J. Gonz\'alez-Anaya]{Javier Gonz\'alez-Anaya}
\address{Department of Mathematics and Computer Science, Santa Clara University, Santa Clara, CA 95053, United States.}
\email{jgonzalezanaya@scu.edu}

\author[J. L. Gonz\'alez]{Jos\'e Luis Gonz\'alez}
\address{Department of Mathematics, University of California, Riverside, Riverside, CA 92521, United States.}
\email{joselg@ucr.edu}

\date{\today}

\begin{document}
\begin{abstract}
The theory of wonderful compactifications provides a common framework encompassing the constructions of important varieties such as the Fulton--MacPherson configuration spaces, the moduli space of $n$-pointed stable rational curves, and the De Concini--Procesi models of subspace arrangements, and, more generally, of smooth compactifications in which the subvarieties in a building set are replaced by a simple normal crossings boundary divisor. 
We give a presentation of the Chow ring of a wonderful compactification, assuming only that the Chow ring of the ambient variety surjects onto that of each subvariety in the induced arrangement. 
Our results recover and extend many previously known presentations and apply to a broad range of new examples.

\end{abstract}

% \dedicatory{Optional dedication.}

\maketitle
% \setcounter{tocdepth}{1}
% \tableofcontents

%%%%%%%%%%%%%%%%%%%%%%%%%%%%%%%%%%%%%%%%%%%%%%%%%%%%%%%%%%%%%%%%
%%%%%%%%%%%%%%%%%%%%%%%%%%%%%%%%%%%%%%%%%%%%%%%%%%%%%%%%%%%%%%%%
%%%%                                                        %%%%
%%%%               1. Introduction    
%%%%                                                        %%%%
%%%%%%%%%%%%%%%%%%%%%%%%%%%%%%%%%%%%%%%%%%%%%%%%%%%%%%%%%%%%%%%%
%%%%%%%%%%%%%%%%%%%%%%%%%%%%%%%%%%%%%%%%%%%%%%%%%%%%%%%%%%%%%%%%

\section{Introduction}
\label{section.introduction}

The theory of wonderful compactifications generalizes several influential constructions in algebraic geometry, using iterated blowups to produce a smooth compactification with a simple normal crossings boundary. Classical examples include De Concini and Procesi's models of subspace arrangements \cite{de1995wonderful}; the Fulton--MacPherson configuration spaces \cite{Fulton-MacPherson}; and the moduli space of stable rational curves $\overline{M}_{0,n}$ \cite{Knudsen, Kapranov-veronese, Keel-thesis}, alongside numerous others (see e.g., 
\cite{ulyanov,hu2003compactification,Chen-Gibney-Krashen,routis2014weighted,clader2023wonderful,polymatroids-moduli}). 

Li~\cite{LiLi} developed a general framework for wonderful compactifications using building sets of smooth subvarieties that encompasses all of these examples. Although Li later obtained a general additive description of their Chow groups~\cite{li2009chow}, presentations of their Chow rings have previously been available only in special cases; see, for example, \cite{petersen2017chow,hu2003compactification,routis2014weighted}.

In this article we give a presentation of the Chow ring of the wonderful compactification associated to a building set under the assumption that the restriction maps on Chow rings from the ambient variety to the subvarieties in the associated arrangement are surjective. 
A similar surjectivity assumption appears in Keel's presentation of the Chow ring of a blowup, which is a key tool in our approach.  
Throughout the article we work over an algebraically closed field of arbitrary characteristic. We recall the basic notions needed to state the main result.

%%%%%%%%%%%%%%%%%%%%%%%%%%%%%%%%%%%%%%%%%%%

Let $X$ be a smooth variety. 
In what follows, all intersections of varieties are scheme-theoretic.  
An \emph{arrangement} $\mathcal{A}$ on $X$ is a finite collection of smooth subvarieties properly contained in $X$ that is closed under nonempty intersections. A subcollection $\mathcal{G}\subseteq\mathcal{A}$ is a \emph{building set} if, for every $Z\in\mathcal{A}$, the minimal elements of $\mathcal{G}$ containing $Z$ intersect transversally with intersection $Z$. These minimal elements are called the \emph{$\mathcal{G}$-factors} of $Z$. 

In the above setting, there exists a smooth variety $X_{\mathcal{G}}$, called the \emph{wonderful compactification} of $\mathcal{G}$ on $X$, together with a morphism $X_{\mathcal{G}}\rightarrow X$ that is an isomorphism over $X^{\circ} := X \smallsetminus \bigcup_{Y \in \mathcal{G}} Y$. The \emph{boundary} $X_{\mathcal{G}}\smallsetminus X^\circ$ is the union of smooth divisors $D_Y$, one for each $Y\in\mathcal{G}$, with $D_Y$ mapping onto $Y$, and every nonempty intersection of these boundary divisors is transversal. 
The incidence relations among the divisors $D_Y$ are determined by those in the original building set via the concept of \emph{$\mathcal{G}$-nested sets}, as in Definition~\ref{definition.G.nested}. 
A collection of boundary divisors has nonempty intersection if and only if the corresponding elements of $\mathcal{G}$ form a $\mathcal{G}$-nested set.

Throughout, we denote the Chow ring of any smooth variety $W$ by $A^\bullet(W)$. 
For $Z\in\mathcal{A}$ let $J_Z$ denote the kernel of the pullback map $A^\bullet(X) \rightarrow A^\bullet(Z)$. For each $Y \in \mathcal{G}$ fix a Chern polynomial $P_{X,Y}(t) \in A^{\bullet}(X)[t]$ in the sense of Definition~\ref{definition.Chern.polynomial}, 
and for each inclusion $Y \subsetneq Z$ with $Y \in \mathcal{G}$ and $Z \in \mathcal{A}$ fix a Chern polynomial $P_{Z,Y}(t)\in A^\bullet(Z)[t]$ and a lift $P_{X,Z,Y}(t)$ of it to $A^{\bullet}(X)[t]$. 
In particular, these Chern polynomials exist if the pullback $A^{\bullet}(X) \rightarrow A^{\bullet}(Z)$ is surjective for each $Z \in \mathcal{A}$.

Let $A^\bullet(X)[E_Y \, | \, Y \in \mathcal{G}]$ denote the polynomial ring  over $A^\bullet(X)$ in the variables $E_Y$ for $Y \in \mathcal{G}$. This polynomial ring maps to $A^\bullet(X_{\mathcal{G}})$ via the pullback $A^\bullet(X)\to A^\bullet(X_{\mathcal{G}})$ on coefficients and by sending each variable $E_Y$ to the class of the boundary divisor $D_Y$ of $X_{\mathcal{G}}$ corresponding to $Y$. For each $Z \in \mathcal{A}$ denote by $E_Z \in A^\bullet(X)[E_Y \, | \, Y \in \mathcal{G}]$ the product of the variables $E_Y$ as $Y$ ranges over the $\mathcal{G}$-factors of $Z$.  
The following is the main theorem in this article. 

\begin{theorem} \label{theorem.main.intro}
    Let $\mathcal{G}$ be a building set on a smooth variety $X$ and let $\mathcal{A}$ be its associated arrangement.
    Assume that the pullback $A^{\bullet}(X) \rightarrow A^{\bullet}(Z)$ is surjective for all $Z \in \mathcal{A}$. 
    Then, the Chow ring of the wonderful compactification $X_\mathcal{G}$ is  $A^\bullet(X_{\mathcal{G}}) \cong A^\bullet(X)[E_Y \, | \, Y \in \mathcal{G}]/I$, where $I$ is the ideal generated by: % (1)-(4) below: 
    \begin{enumerate}
    \setlength\itemsep{0.7em}
        \item \ \ $\prod_{Y\in\mathcal{T}} E_Y$ \quad for all $\mathcal{T}\subseteq\mathcal{G}$ such that $\mathcal{T}$ is \emph{not} $\mathcal{G}$-nested;
        \item \ \ $J_Z \cdot E_Z$  \quad for all $Z$ in $\mathcal{A}$;
        \item \ \ $P_{X,Y}\, \raisebox{-1.5mm}{$\Biggl($}
        -\mathlarger{\sum}\nolimits_{\substack{Y'\in\mathcal{G}\\Y'\subseteq Y}} \, E_{Y'}
        \raisebox{-1.5mm}{$\Biggr)$}$
        \quad for all $Y\in\mathcal{G}$; 
        \item \ \ $E_Z \cdot P_{X,Z,Y}\, \raisebox{-1.5mm}{$\Biggl($}
        -\mathlarger{\sum}\nolimits_{\substack{Y'\in\mathcal{G}\\Y'\subseteq Y}} \, E_{Y'}
        \raisebox{-1.5mm}{$\Biggr)$}$
        \quad
        for all $Z\in\mathcal{A}$ and $Y\in\mathcal{G}$
        such that $Y\subsetneq Z$.
    \end{enumerate}
\end{theorem}
Let us discuss the intuition behind relations (1)--(4). Given $Z\in\mathcal{A}$, let $D_Z$ denote the intersection of the boundary divisors corresponding to the $\mathcal{G}$-factors of $Z$, and let $i_Z:D_Z\rightarrow X_{\mathcal{G}}$ be the inclusion. Since these boundary divisors intersect transversally, the image of $E_Z$ in $A^\bullet(X_{\mathcal{G}})$ is the class of $D_Z$, and the projection formula gives $(i_Z)_*i_Z^*(\alpha)=\alpha\cdot E_Z$ for every $\alpha\in A^\bullet(X_{\mathcal{G}})$.

Relations (1) capture the fact that a collection of boundary divisors has empty intersection if and only if the corresponding collection of elements of the building set is not $\mathcal{G}$-nested, while relations (2) reflect the fact that the elements of $J_Z\subseteq A^\bullet(X)\subseteq A^\bullet(X_{\mathcal{G}})$ restrict to zero in $A^\bullet(D_Z)$. Under the hypotheses of Keel's formula for a blowup $\operatorname{Bl}_W V$ with exceptional divisor $E$, the Chern polynomial satisfies 
$P_{V,W}(-E)=0\in A^\bullet(\operatorname{Bl}_W V)$.
Relations (3) keep track of the vanishing of these expressions $P_{V,W}(-E)$ through the iterated blowup construction of $X_{\mathcal{G}}$. Lastly, relations (4) are the analogues of the relations (3) on the boundary strata $D_Z$, where the Chern polynomial $P_{X,Z,Y}$ records the blowup relation arising from $Y\subsetneq Z$.

The form of our presentation in Theorem~\ref{theorem.main.intro} is influenced by earlier presentations of the Chow rings of particular  wonderful compactifications.
Fulton and MacPherson proposed a presentation of the Chow rings of their compactifications in \cite[Theorem 6]{Fulton-MacPherson}. Later, Routis introduced a weighted variant of Fulton--MacPherson compactifications and wrote a presentation of their Chow rings in \cite[Theorem 8]{routis2014weighted}. The first complete proof of these two presentations was given by Petersen in \cite[Theorems 3.2 and 3.4]{petersen2017chow}; our proof is inspired by Petersen's inductive argument. Weighted Fulton--MacPherson spaces are wonderful compactifications, and our relations (1)--(4) recover the presentations found by Fulton--MacPherson and Routis (Example~\ref{example.WFM}).

Despite the similarities, our proof encounters numerous technical obstacles that can occur for general building sets, but are not present for the building sets arising from collections of diagonals in the $n$-fold self-product $X^n$ of a smooth algebraic variety $X$, as in the weighted Fulton--MacPherson case. For example, the incidence relations among elements of a general building set play an important role in the induction, and they can be significantly more complex than those among the diagonals in $X^n$. Additionally, Chern polynomials lifted from elements of the arrangement to the ambient variety provide the new relations (4), which have a simpler form using only the building set in the weighted Fulton--MacPherson case. To illustrate such differences, in Section~\ref{section.examples} we give examples to show that using direct analogues of the relations in those previous works does not suffice to describe the Chow rings of wonderful compactifications in general.  

The proof of Theorem~\ref{theorem.main.intro} relies on a delicate induction on the cardinality of the building set, whose inductive step requires identifying certain strict transforms of elements of the building set as wonderful compactifications to which the induction hypothesis can be applied. This motivates us to introduce and study the notion of \emph{restriction of building sets} (Definition~\ref{definition.restriction.building.set}) as a tool to carry out general inductive arguments for wonderful compactifications. We show that restriction gives rise to new building sets (Theorem~\ref{theorem.restriction.of.building.set.is.building.set}) and we describe their poset structures under containment. We also show that the wonderful compactification of a restricted building set can be obtained as a strict transform (Theorem~\ref{strict.transforms.are.wonderful}). Moreover, we show that the iterated blowup constructions of the wonderful compactifications of a building set and of its restriction to a maximal member can be carried out compatibly (Theorem~\ref{theorem.compatibility.restriction.iterated.blowups}; cf.~Example~\ref{example.order.respecting.containments.needed}). This can be summarized as follows.

\begin{theorem}[Theorems \ref{theorem.restriction.of.building.set.is.building.set}, \ref{strict.transforms.are.wonderful} and \ref{theorem.compatibility.restriction.iterated.blowups}
] \label{theorem.restriction.building.sets.intro}
    Let $X$ be a smooth variety, $\mathcal{G}$ be a building set on $X$, and $Y$ be an element of $\mathcal{G}$. 
    \begin{itemize}[leftmargin=2em]
        \item[(i)] The restriction $\mathcal{G}|_{Y} := \left\{ Y \cap Z \ | \ Z \in \mathcal{G} \textnormal{ such that } \emptyset \neq Y \cap Z \subsetneq Y \right\}$ of $\mathcal{G}$ to $Y$ is  a building set on $Y$. 
    \end{itemize}

    \noindent Moreover, if $Y$ is maximal in $\mathcal{G}$,  then   $\mathcal{G}':=\mathcal{G} \smallsetminus \{Y \}$ is a building set on $X$, and we have:  
    
    \begin{itemize}[leftmargin=2em] 
        \item[(ii)] There is a canonical isomorphism $Y_{\mathcal{G}|_Y} \cong
    \widetilde{Y}$ between the wonderful compactification $Y_{\mathcal{G}|_Y}$ of $\mathcal{G}|_Y$ on $Y$ and the strict transform $\widetilde{Y}$ of $Y$ in the wonderful compactification $X_{\mathcal{G}'}$.
        \item[(iii)] The iterated blowup constructions of wonderful compactifications for $X_{\mathcal{G}'}$ and $Y_{\mathcal{G}|_Y}$ are compatible under restriction for any order of $\mathcal{G}'$ that respects containments. 
    \end{itemize}
\end{theorem}
While Theorem~\ref{theorem.main.intro} provides a general presentation, many wonderful compactifications in the literature are endowed with additional geometric properties that simplify relations (1)--(4). In Section~\ref{section.examples} we identify natural properties of building sets (see Definition~\ref{definition.complete.regular.retractable}) that simplify these relations. For example, we define a \emph{complete} building set as one for which $\mathcal{G}=\mathcal{A}$. Proposition~\ref{prop.complete.relations.1} proves that for complete building sets relations (1) are equivalent to  
\begin{itemize}
    \item[(1')] \ \ $E_Y\cdot E_{Y'}$ \quad for all $Y,Y'\in \mathcal{G}$ such that $\{Y,Y'\}$ is \emph{not} $\mathcal{G}$-nested.
\end{itemize}
Similarly, by Proposition~\ref{prop.reduction.relations.2}, for building sets with \emph{retractable} arrangements relations (2) are equivalent to 
\begin{itemize}
    \item[(2')] \ \ $J_Y \cdot E_Y$  \quad for all $Y$ in $\mathcal{G}$.
\end{itemize}
We also discuss conditions under which relations (4) simplify. For example, by Proposition~\ref{prop.transversally.regular.relations.4}, in a \emph{transversally regular} building set relations (4) are implied by relations (1)--(3). Moreover, by Proposition~\ref{prop.weak.transversally.regular.relations.4}, for a building set that admits \emph{weak transversal witnesses}, in the presence of relations (1)--(3), relations (4) are equivalent to
\begin{itemize}
    \item[(4')] \ \ $E_Z \cdot P_{X,Z,Y} \raisebox{-1.5mm}{$\Biggl($}
-\mathlarger{\sum}\nolimits_{\substack{Y'\in\mathcal{G}\\ Y'\subseteq Y}} E_{Y'}\raisebox{-1.5mm}{$\Biggr)$}$ for all $Z,Y\in\mathcal{G}$ such that $Y\subsetneq Z$.
\end{itemize}

Lastly, while reviewing the literature, we learned from \cite{petersen2017chow} that the original statement of \cite[Lemma~5.4]{Fulton-MacPherson} omits one relation; a corrected version is given in \cite[Lemma~2.2]{petersen2017chow}. We also observed that the presentation in \cite[Theorem~1.8]{hu2003compactification} of the Chow rings of what the author calls regular simple arrangements of subvarieties omits relations (4). Neither result is used in our arguments. Instead, we give a more general version of the former in Lemma~\ref{lemma.Jideal.blowup} and a corrected version of the latter in Section~\ref{sec:based.hu.example}.

This article is organized as follows. 
In Section~\ref{section.background} we review background material on transversality, Chow theory, and wonderful compactifications. 
In Section~\ref{section.restriction} we study restrictions of building sets and their main properties, proving Theorem~\ref{theorem.restriction.building.sets.intro}. 
In Section~\ref{section.wonderful} we present the inductive proof of the main theorem.
Finally, in Section~\ref{section.examples} we present examples that illustrate the main theorem and compare it with analogues for particular wonderful compactifications in the literature.

\begingroup
\renewcommand{\addcontentsline}[3]{}
\subsection*{Acknowledgments}
The authors are grateful for the supportive environment provided by
their home institutions, the University of California, Riverside (PG and JLG), and Santa Clara University (JGA). 
Jos\'e Gonz\'alez was supported by a grant from the Simons Foundation
(Award Number 710443). Patricio Gallardo is partially supported by the National Science
Foundation under Grant No. DMS-2316749. 
Patricio Gallardo thanks Long Horizon Research, Inc. for access to the workspace Sundial. The authors are grateful to Chris Eur for asking us about the Chow rings of the spaces we introduced in \cite{polymatroids-moduli}, which are examples of wonderful compactifications.
\endgroup

\begingroup
\renewcommand{\addcontentsline}[3]{}
\subsection*{Statement on the use of AI} 
Publicly available AI Models from OpenAI and Anthropic, such as ChatGPT 5.5 and Claude 4.8, were used to answer general questions and search for references.  All logical arguments, proofs, and writing are the authors' own.
\endgroup

\tableofcontents

%%%%%%%%%%%%%%%%%%%%%%%%%%%%%%%%%%%%%%%%%%%%%%%%%%%%%%%%%%%%%%%%
%%%%%%%%%%%%%%%%%%%%%%%%%%%%%%%%%%%%%%%%%%%%%%%%%%%%%%%%%%%%%%%%
%%%%                                                        %%%%
%%%%     2. Background    
%%%%                                                        %%%%
%%%%%%%%%%%%%%%%%%%%%%%%%%%%%%%%%%%%%%%%%%%%%%%%%%%%%%%%%%%%%%%%
%%%%%%%%%%%%%%%%%%%%%%%%%%%%%%%%%%%%%%%%%%%%%%%%%%%%%%%%%%%%%%%%

\section{Background: Chow theory and wonderful compactifications} \label{section.background}

This section provides a brief review of the core concepts used throughout this article. In the first subsection, we outline Li's theory of wonderful compactifications. Transversality is the cornerstone of this theory and thus plays a central role in proving our main result. The next two subsections examine transversality and present several results used later. Finally, we establish the necessary notation, terminology, and key results concerning Chow rings.

\begin{notation}
We fix the following notation and conventions throughout the article.
    \begin{itemize}
    \item We work over an algebraically closed field of arbitrary characteristic. We use the terminology of \cite{Hartshorne} for algebraic geometry. In particular, varieties are reduced and irreducible. By a subvariety we mean a closed subvariety and by a point of a variety we mean a closed point. 
    \item Unless specified otherwise, when we refer to a maximal or minimal element of a set it is with respect to containment. 
    \item Unless specified otherwise, intersections of schemes (e.g. of varieties) are scheme-theoretic. 
    \item The tangent space of the smooth variety $X$ at a point $x \in X$ is denoted by $T_{X,x}$ or simply by $T_x$ when the ambient variety is clear.  
    For a smooth subvariety $Y$ of $X$ and a point $x \in Y$, we consider the tangent space  $T_{Y,x}$ as a subspace of $T_{X,x}$. 
    For $x\notin Y$, we set $T_{Y,x}:=T_{X,x}$ by convention. 
    \end{itemize}
\end{notation}

\subsection{On Li's theory of wonderful compactifications}
\label{subsection.lili}

First we review the necessary tools from Li's theory of wonderful compactifications \cite{LiLi}. Throughout this subsection we fix $X$ to be a smooth variety.

\begin{definition}[Transverse intersections]
\label{definition.transversal.intersection} Let $X$ be a smooth variety and let $Y_1,\ldots,Y_r$ be smooth subvarieties of $X$.
One says that $Y_1,\ldots,Y_r$ \emph{intersect transversally} (or that they are \emph{transverse}) on $X$ at a point $p\in X$ if there exists a regular system of parameters $x_1,\ldots,x_n\in\mathcal{O}_{X,p}$ and pairwise disjoint
subsets $S_1,\ldots,S_r\subseteq\{1,\ldots,n\}$ such that $I(Y_i)_p=\langle x_j\mid j\in S_i\rangle$ for every $i$ with $p\in Y_i$. 
For $p\notin Y_i$, we set $S_i=\emptyset$; in this case $I(Y_i)_p=\mathcal{O}_{X,p}$, and $Y_i$ imposes no condition on transversality at $p$.
One says that $Y_1,\ldots,Y_r$ \emph{intersect transversally} (or that they are \emph{transverse}) on $X$ if they do so at every point of $X$.
\end{definition}

\begin{remark}
    The definition of a transversal intersection in Definition~\ref{definition.transversal.intersection} is equivalent to that in \cite[Section~5.1.2]{LiLi}; see Section~\ref{subsection.course.transversal.intersections} for additional useful characterizations and basic properties. 
\end{remark}

\begin{definition}\label{bg: li arrangement def}
An \emph{arrangement} of subvarieties of $X$ is a finite collection $\mathcal{A}=\{Y_i\}$ of smooth subvarieties $Y_i\subsetneq X$ such that, for all $i\neq j$, the intersection $Y_i\cap Y_j$ is either empty or equal to some $Y_k$. In particular, $\mathcal{A}$ is closed under nonempty intersections.
\end{definition}

\begin{remark}
The definition of an arrangement in Definition~\ref{bg: li arrangement def} is equivalent to the definition of a (simple) arrangement in \cite[Definition~2.1]{LiLi} (see \cite[Section~5.1.1]{LiLi} for details).
\end{remark}

\begin{definition}\label{bg: li def building set}
    Let $\mathcal{A}$ be an arrangement of subvarieties of $X$. 
    \begin{itemize}
        \item A subset $\mathcal{G}\subseteq\mathcal{A}$ is called a \emph{building set of $\mathcal{A}$} if, for all $Y\in\mathcal{A}$, the minimal elements of $\{G\in\mathcal{G}\,\vert\, Y\subseteq G\}$ intersect transversally, and their intersection is $Y$. In particular, this condition is automatically satisfied if $Y\in\mathcal{G}$. 
        \item Given a building set $\mathcal{G}$ of $\mathcal{A}$ and an element $Y\in\mathcal{A}$, the minimal elements of $\mathcal{G}$ containing $Y$ are called the \emph{$\mathcal{G}$-factors} of $Y$. We denote the set of $\mathcal{G}$-factors of $Y$ by $\operatorname{Fact}(\mathcal{G},Y)$.
        \item A collection $\mathcal{G}$ of subvarieties of $X$ is called a \emph{building set} if the set $\mathcal{A}$ of all possible nonempty intersections of elements of $\mathcal{G}$ is an arrangement, and $\mathcal{G}$ is a building set of $\mathcal{A}$. In this situation, $\mathcal{A}$ is called the \emph{arrangement induced} by $\mathcal{G}$. 
    \end{itemize}
\end{definition}
We arrive at the main definition of this subsection.

\begin{definition}\label{bg: def wonderful}
Let $\mathcal{G}$ be a building set on $X$.
If $\mathcal{G}$ is nonempty, set $X^\circ:=X\smallsetminus\bigcup_{G\in\mathcal{G}}G$. The closure of the image of the locally closed embedding
\[
X^\circ\hookrightarrow\prod_{G\in\mathcal{G}}\operatorname{Bl}_G X
\]
is called the \emph{wonderful compactification} of $\mathcal{G}$
on $X$ and is denoted by $X_\mathcal{G}$.
The projections $X_\mathcal{G}\to\operatorname{Bl}_G X\to X$ all agree, since they restrict to the identity on the dense open subset $X^\circ$. This defines a natural birational morphism $X_\mathcal{G}\to X$. 
If $\mathcal{G}=\emptyset$, we define $X_\mathcal{G}:=X$ and the map to $X$ is the identity.
\end{definition}

The main theorem in Li's article \cite{LiLi} is the following.

\begin{theorem}[{\cite[Theorem 1.2]{LiLi}}]
    Let $X$ be a smooth variety and let $\mathcal{G}$ be a nonempty building set. Then, the wonderful compactification $X_\mathcal{G}$ is a smooth variety. Moreover, for each $Y\in\mathcal{G}$ there is a smooth divisor $D_Y\subseteq X_\mathcal{G}$ mapping onto $Y$ under the natural morphism $X_\mathcal{G} \rightarrow X$ such that:
    \begin{enumerate}[label=(\roman*),ref=(\roman*)] 
        \item the union of these divisors is $X_\mathcal{G}\smallsetminus X^\circ$;
        \item any set of these divisors meets transversally, and an intersection of divisors $D_{Y_1}\cap\cdots\cap D_{Y_r}$ is nonempty if and only if $\{Y_1,\dots,Y_r\}$ is $\mathcal{G}$-nested (see the following definition).
    \end{enumerate}
\end{theorem}

\begin{definition}\label{definition.G.nested}
A subset $\mathcal{T}\subseteq\mathcal{G}$ is said to be \emph{$\mathcal{G}$-nested} if there is a flag $Z_1\subseteq\cdots\subseteq Z_r$ of elements of the arrangement $\mathcal{A}$ such that
\[
\mathcal{T}=\bigcup_{i=1}^r \, 
\{Y\in\mathcal{G}\mid Y\text{ is a $\mathcal{G}$-factor of }Z_i\}.
\]
In this case, we say that $\mathcal{T}$ is induced by the flag. 
We allow the empty flag, which induces $\mathcal{T}=\emptyset$.
\end{definition}

To keep track of subvarieties contained in the blowup center,
whose strict transforms are empty, we use the following notion.

\begin{definition}[Dominant transform]
    Let $Z$ be a smooth subvariety of a smooth variety $X$, and $\pi: \Bl_Z X \to X$ be the blow-up of $X$ along $Z$. For any subvariety $V$ of $X$, we define the \emph{dominant transform of $V$} to be the strict transform $\overline{\pi^{-1}(V \smallsetminus Z)}$ of $V$ if $V\nsubseteq Z$, and to be the scheme-theoretic inverse image $\pi^{-1}(V)$ if $V\subseteq Z$. We denote the dominant transform by $\widetilde{V}$.
\end{definition}

Wonderful compactifications can be constructed by successive blowups along smooth centers. The following theorem gives a sufficient condition on the order of the building-set elements for their successive dominant transforms to yield such a construction.

\begin{theorem}[{\cite[Theorem~1.3]{LiLi}}] 
\label{thm: li blow-up building set}
    Let $X$ be a smooth variety and let $\mathcal{G} = \{Y_1 , Y_2 ,\dots , Y_n\}$ be a building set on $X$. Suppose the elements of $\mathcal{G}$ are ordered in such a way that the first $k$ terms $Y_1,Y_2,\dots, Y_k$ form a building set for all $1\leq k\leq n$. Then, the wonderful compactification $X_\mathcal{G}$ is isomorphic to the iterated blow-up $\Bl_\mathcal{G} X$ of $X$ along the dominant transforms of the elements of $\mathcal{G}$ in this chosen order, that is,
    \[
        X_\mathcal{G} = \Bl_{\widetilde{Y}_n}\cdots\Bl_{\widetilde{Y}_2}\Bl_{Y_1} X,
    \]
    where each blow-up is along a smooth subvariety.
\end{theorem}

For orders compatible with containments, we use the following inductive construction.

\begin{definition}[{\cite[Definition~2.12]{LiLi}}] \label{bg.def.transform.of.building.set}
Suppose $\mathcal{G}=\{Y_1,\dots,Y_n\}$ is a building set on $X$ with induced arrangement $\mathcal{A}$, ordered compatibly with
containments (i.e. $Y_i \subseteq Y_j$ implies $i \leq j$). For $0\leq m\leq n$, define the triple $\left(X^{(m)},\mathcal{G}^{(m)},\mathcal{A}^{(m)}\right)$ inductively. For $m=0$, set $X^{(0)}:=X$, $\mathcal{G}^{(0)}:=\mathcal{G}$, $\mathcal{A}^{(0)}:=\mathcal{A}$, and $Y_i^{(0)}:=Y_i$ for $1\leq i\leq n$.
For $1\leq m\leq n$, let $X^{(m)}$ be the blow-up of $X^{(m-1)}$ along $Y_m^{(m-1)}$, and set $Y_i^{(m)}:=\widetilde{Y_i^{(m-1)}}$ for $1\leq i\leq n$. 
Set $\mathcal{G}^{(m)}:=\{Y_1^{(m)},\dots,Y_n^{(m)}\}$, and let $\mathcal{A}^{(m)}$ be the set of all nonempty
intersections of members of $\mathcal{G}^{(m)}$.
By \cite[Proposition~2.8]{LiLi}, each $\mathcal{G}^{(m)}$ is a building set on $X^{(m)}$ with induced arrangement $\mathcal{A}^{(m)}$.
\end{definition}

\begin{proposition}[{\cite[Proposition~2.13]{LiLi}}]
\label{bg.prop.transform.of.building.set}
    In the setting of the previous definition, $X_\mathcal{G} \cong X^{(n)}$.
\end{proposition}

In light of this construction of $X_\mathcal{G}$ we introduce some useful terminology:

\begin{definition}[Strict transforms in a wonderful compactification] 
    Given a smooth variety $X$, a building set $\mathcal{G}$, and a subvariety $Y$ of $X$ that is not contained in $\bigcup_{Z\in\mathcal{G}} Z$, we have the inclusions
    \[
    Y\smallsetminus\bigcup_{Z\in\mathcal{G}}Z \ \subseteq \ X\smallsetminus\bigcup_{Z\in\mathcal{G}}Z \ \subseteq X_{\mathcal{G}}. 
    \]
    We define the \emph{strict transform} of $Y$ in $X_{\mathcal{G}}$ to be the closure of $Y^\circ:=Y\smallsetminus\bigcup_{Z\in\mathcal{G}}Z$ inside $X_{\mathcal{G}}$. This terminology agrees with the usual notion of strict transform when $X_{\mathcal{G}}$ is considered as an iterated blowup of $X$ along the dominant transforms of the elements of $\mathcal{G}$, but does not depend on a particular iterated blow up presentation of $X_{\mathcal{G}}$. 
\end{definition}

To conclude this brief overview we recall some results that will be useful.

\begin{lemma}\label{lemma: useful lemmas from li's}
    Let $X$ be a smooth variety.
    \begin{enumerate}[label=(\roman*),ref=(\roman*)]
        \item \label{lemma: useful lemmas from li's item 1} (\cite[Lemma~2.6(i)]{LiLi}) Let $\mathcal{G}$ be a building set on $X$, and $F\in\mathcal{G}$ be minimal. Then, any other $G\in\mathcal{G}$ either contains $F$ or intersects $F$ transversally.
        \item \label{lemma: useful lemmas from li's item 2} (\cite[Lemma~2.9(ii)]{LiLi}) Let $V$ and $W$ be smooth subvarieties of $X$. If the intersection of $V$ and $W$ is a smooth variety $Z$ properly contained in $V$ and $W$, then $\widetilde{V}\cap\widetilde{W}=\emptyset$ in $\operatorname{Bl}_Z X$ (see also \cite[Ex. II.7.12]{Hartshorne} for a version without smoothness). 
    \end{enumerate} 
\end{lemma}

%%%%%%%%%%%%%%%%%%%%%%%%%%%%%%%%%%%%%%%%%%%%%%%%%%%%%%%%%%%%%%%%%%%%%%%

\subsection{Crash course on transverse intersections} \label{subsection.course.transversal.intersections}

The notion of a transversal intersection (Definition~\ref{definition.transversal.intersection}) is fundamental to our article. In this subsection we review the basic properties that we will use throughout this work. We start with the following characterizations, which appear without a proof in \cite[Section 5.1.2]{LiLi}. 

\begin{lemma} \label{characterizations.transverse.intersection}
Let $X$ be a smooth variety and let $Y_1,\ldots,Y_r$ be smooth subvarieties of $X$. 
For each point $p \in X$ the following are equivalent: 
\begin{enumerate}[label=(\roman*),ref=(\roman*)]
    \item \label{equivalence.item.i} $Y_1,\ldots,Y_r$ intersect transversally on $X$ at $p$, in the sense of Definition~\ref{definition.transversal.intersection}. 
    \item \label{equivalence.item.ii} The sum $T_{Y_1, p}^{\perp}+T_{Y_2, p}^{\perp}+\cdots+T_{Y_r, p}^{\perp}$ is direct in the dual space $T_{X,p}^*$ of $T_{X,p}$.
    \item \label{equivalence.item.iii} $\operatorname{codim}\left(\bigcap_{i=1}^r T_{Y_i, p}, T_{X,p} \right)=\sum_{i=1}^r \operatorname{codim}_p\left(Y_i, X\right)$. 
\end{enumerate}
In \ref{equivalence.item.iii}, $\operatorname{codim}_p\!\left(Y_i, X\right) :=  \operatorname{codim}\!\left(Y_i, X\right)$ if $p \in Y_i$ and $\operatorname{codim}_p\!\left(Y_i, X\right) := 0$ if $p \notin Y_i$, and in \ref{equivalence.item.ii}  
for any subspace $W$ of a vector space $V$ we define
\[
W^{\perp} 
:= \{ \phi \in  V^* \ | \  \phi(w) = 0 \textnormal{ for all } w \in W  \} \subseteq V^*.
\]
\end{lemma}

\begin{proof}
We fix a point $p$ in $X$. Since the conditions in \ref{equivalence.item.i}, \ref{equivalence.item.ii}  and \ref{equivalence.item.iii} depend only on the $Y_i$ that contain $p$, we may assume that $p$ is in $Y_i$ for every $i$.   

\emph{\ref{equivalence.item.i}$\Rightarrow$\ref{equivalence.item.ii}.}
Let $\mathfrak m_p$ denote the maximal ideal of $\mathcal O_{X,p}$.  
By Definition~\ref{definition.transversal.intersection}, there is a regular system of parameters $x_1,\dots,x_n \in \mathcal O_{X,p}$ and pairwise disjoint subsets $S_i\subseteq\{1,\dots,n\}$ such that $I(Y_i)_p=\langle x_j\mid j\in S_i\rangle \subseteq \mathcal O_{X,p}$. 
The classes of $x_1,\dots,x_n$ form a basis of $\mathfrak m_p/\mathfrak m_p^2$ as a $\mathcal{O}_{X,p}/\mathfrak m_p$-vector space.
Via the identification $T_{X,p}^*\cong \mathfrak m_p/\mathfrak m_p^2$, the space $T^\perp_{Y_i,p}$ is the span of the classes of $\{x_j \, | \, j\in S_i\}$. The disjointness of the $S_i$ makes these spans linearly independent, hence $\sum_i T^\perp_{Y_i,p}$ is a direct sum. 

\smallskip

\emph{\ref{equivalence.item.ii}$\Rightarrow$\ref{equivalence.item.iii}.}
For subspaces $U_i$ of a finite dimensional vector space we have $(U_i^{\perp})^{\perp}=U_i$ and $\left(U_1+\cdots+U_r\right)^\perp=\bigcap_i U_i^\perp$. 
Taking $U_i:=T^\perp_{Y_i,p}\subseteq T_{X,p}^*$, we get
\[
\operatorname{codim}\left(\bigcap_{i=1}^r T_{Y_i,p},T_{X,p}\right)
= \operatorname{codim}\left( \left(\sum_{i=1}^r T_{Y_i, p}^{\perp}\right)^{\perp} , T_{X,p} \right)
=\dim\left(\sum_{i=1}^r T^\perp_{Y_i,p}\right). 
\] 
By \ref{equivalence.item.ii} the sum on the right is direct, so its dimension is $\sum_{i=1}^r \dim T^\perp_{Y_i,p}$. Since $\dim T^\perp_{Y_i,p}=\operatorname{codim}_p(Y_i,X)$, we obtain \ref{equivalence.item.iii}.

\smallskip

\emph{\ref{equivalence.item.iii}$\Rightarrow$\ref{equivalence.item.i}.}
Using the same orthogonality relation as before, condition \ref{equivalence.item.iii} gives
\[
\dim\!\left(\sum_{i=1}^r T^\perp_{Y_i,p}\right)
=
\operatorname{codim}\!\left(\bigcap_{i=1}^r T_{Y_i,p},T_{X,p}\right)
=
\sum_{i=1}^r \operatorname{codim}_p(Y_i,X)
=
\sum_{i=1}^r \dim T^\perp_{Y_i,p}.
\]
Thus the sum $\sum_i T^\perp_{Y_i,p}\subseteq T_{X,p}^*$ is a direct sum. 
For each $i$, choose a basis of $T^\perp_{Y_i,p}$. Since the sum is a direct sum, the union of these bases is linearly independent in
\[
T_{X,p}^*\cong \mathfrak m_p/\mathfrak m_p^2.
\]
Extend this union to a basis of $T_{X,p}^*$, and choose lifts $x_1,\ldots,x_n\in \mathfrak m_p$ of the resulting basis. We choose the lifts corresponding to the basis of $T^\perp_{Y_i,p}$ inside $I_i:=I(Y_i)_p$. 
This is possible because $T^\perp_{Y_i,p}$ is the image of $I_i$ in
$\mathfrak m_p/\mathfrak m_p^2$. 
The classes of $x_1,\ldots,x_n$ form a basis of
$\mathfrak m_p/\mathfrak m_p^2$, so Nakayama's lemma implies that they generate
$\mathfrak m_p$. Since $\mathcal O_{X,p}$ is regular of dimension $n$, they form a regular system of parameters. 

It remains to show that the ideals $I_i$ are generated by the corresponding subsets of these parameters.  
Let $S_i\subseteq\{1,\ldots,n\}$ be the set of indices corresponding to the chosen basis of $T^\perp_{Y_i,p}$, and set $J_i:=\langle x_j\mid j\in S_i\rangle \subseteq I_i$. 
Since $X$ and $Y_i$ are smooth at $p$, the conormal sequence at $p$,  
\[
0\longrightarrow I_i/I_i^2
\longrightarrow \Omega_{X,p}\otimes_{\mathcal O_{X,p}}\mathcal O_{Y_i,p}
\longrightarrow \Omega_{Y_i,p}
\longrightarrow 0
\]
is exact, with free terms, and hence split. After tensoring over $\mathcal O_{Y_i,p}$ with 
$k(p)=\mathcal O_{X,p}/\mathfrak m_p=\mathcal O_{Y_i,p}/\mathfrak m_{Y_i,p}$, 
and using the standard identifications with cotangent spaces, we get the exact sequence 
\[
0\longrightarrow I_i/\mathfrak m_p I_i
\longrightarrow \mathfrak m_p/\mathfrak m_p^2
\longrightarrow \mathfrak m_{Y_i,p}/\mathfrak m_{Y_i,p}^2
\longrightarrow 0.
\]
Under the identification $\mathfrak m_p/\mathfrak m_p^2\simeq T_{X,p}^*$, the image of
$I_i/\mathfrak m_p I_i$ is $T_{Y_i,p}^{\perp}$. Hence the chosen elements
$x_j$ for $j\in S_i$ form a basis of $I_i/\mathfrak m_p I_i$. Therefore
\[
I_i=J_i+\mathfrak m_p I_i,
\]
and Nakayama's lemma gives $I_i=J_i$.
Therefore each $I(Y_i)_p$ is generated by the subset $\{x_j\mid j\in S_i\}$ of the regular system of parameters $x_1,\ldots,x_n$. Since the subspaces $T^\perp_{Y_i,p}$ form a direct sum, the index sets $S_1,\ldots,S_r$ are pairwise disjoint. Hence $Y_1,\ldots,Y_r$ intersect transversally at $p$, in the sense of Definition~\ref{definition.transversal.intersection}.  
\end{proof}

\begin{example}[Transversality for two subvarieties]\label{example.transversality.two.subvarieties}
    By Lemma~\ref{characterizations.transverse.intersection}, two smooth subvarieties $Y_1$ and $Y_2$ of a smooth variety $X$ intersect transversally on $X$ if and only if $T_{Y_1, p}^{\perp}+T_{Y_2, p}^{\perp}$ is a direct sum in the dual space $T_{X,p}^* \cong \mathfrak{m}_p / \mathfrak{m}_p^2$ of $T_{X,p}$ for each point $p \in X$, or, equivalently, if and only if 
    $$
    T_{X,p}=T_{Y_1, p}+T_{Y_2, p} \quad \textnormal{for all } p \in X. 
    $$
\end{example}

In the following lemma, we use the same definition of transversality for smooth closed subschemes, allowing empty or disconnected intersections. For members of an arrangement, all these intersections are empty or smooth varieties. 

\begin{lemma} \label{lemma.properties.transversal.intersections} 
Let $X$ be a smooth variety and let $Y_1,\ldots,Y_r$ be smooth subvarieties of $X$ that intersect transversally. Then,
\begin{enumerate}[label=(\roman*),ref=(\roman*)]
\item \label{lemma.properties.transversal.intersections.item.one.subset} Any nonempty subcollection of $Y_1,\ldots,Y_r$ intersects transversally. 
\item \label{lemma.properties.transversal.intersections.item.disjoint.subsets} If $I_1,\ldots,I_m$ are nonempty pairwise disjoint subsets of $\{1,\ldots,r\}$ and we set $Y_{I_j} := \bigcap_{i \in I_j} Y_i$, then $Y_{I_1},\ldots,Y_{I_m}$ intersect transversally. In particular, each $Y_{I_j}$ is smooth.
\item \label{lemma.properties.transversal.intersections.item.restrictions} The intersections $Y_1 \cap Y_r, \ldots, Y_{r-1} \cap Y_r$ are smooth and intersect transversally on $Y_r$.
\end{enumerate}
\end{lemma}

\begin{proof}
Fix $p\in X$. By Definition~\ref{definition.transversal.intersection}, there is a regular system of parameters $x_1,\ldots,x_n\in\mathcal{O}_{X,p}$ and pairwise disjoint index sets $S_1,\ldots,S_r\subseteq\{1,\ldots,n\}$ such that $I(Y_i)_p=\langle x_j\mid j\in S_i\rangle$ for every $i$ with $p\in Y_i$. For $p\notin Y_i$, we have $I(Y_i)_p=\mathcal{O}_{X,p}$ and set $S_i=\emptyset$.

\begin{enumerate}[label=(\roman*),ref=(\roman*)]
\item For any nonempty $J\subseteq\{1,\ldots,r\}$, the same regular system of parameters together with the disjoint family $\{S_i\}_{i\in J}$ proves the transversality of $\{Y_i\}_{i\in J}$ at $p$.

\item For each $j$ with $p\in Y_{I_j}$, all the $Y_i$ with $i\in I_j$ contain $p$. Set $S_{I_j}:=\bigsqcup_{i\in I_j}S_i$. Then
\[
I(Y_{I_j})_p=\sum_{i\in I_j}I(Y_i)_p=\langle x_\ell\mid\ell\in S_{I_j}\rangle.
\]
Thus each such $Y_{I_j}$ is smooth at $p$, and these parameter sets are pairwise disjoint. The intersections $Y_{I_j}$ not containing $p$ impose no condition on transversality at $p$. Hence $Y_{I_1},\ldots,Y_{I_m}$ intersect transversally at $p$.

\item For this part, take $p\in Y_r$ and write $S:=S_r$. In $\mathcal{O}_{Y_r,p}\cong\mathcal{O}_{X,p}/\langle x_j\mid j\in S\rangle$, the images $\{\overline{x}_k\mid k\notin S\}$ form a regular system of parameters. For each $i<r$ with $p\in Y_i$,
\[
I\!\left((Y_i\cap Y_r)\subseteq Y_r\right)_p
=\frac{I(Y_i)_p+I(Y_r)_p}{I(Y_r)_p}
=\langle\overline{x}_\ell\mid\ell\in S_i\rangle.
\]
Thus each such $Y_i\cap Y_r$ is smooth at $p$, and their defining parameter sets are pairwise disjoint. The remaining intersections do not contain $p$ and impose no condition on transversality there. Hence $Y_1\cap Y_r,\ldots,Y_{r-1}\cap Y_r$ intersect transversally on $Y_r$ at $p$.
\end{enumerate}
Since $p$ was arbitrary, the claims hold globally.
\end{proof}

%%%%%%%%%%%%%%%%%%%%%%%%%%%%%%%%%%%%%%%%%%%%%%%%%%%%%%%%%%%%%%%%%%%%%%%

\subsection{Transversal intersections and additivity of codimensions} 

In this subsection we review the relation between additivity of codimensions and transversality of intersections, and then deduce some consequences for the intersection of elements of a building set.

\begin{lemma}
    \label{subadditivity.codimension.smooth.ambient}
    Let $X$ be a smooth variety and let $Y_1,\dots,Y_r$ be subvarieties of $X$. For every irreducible component $W$ of the scheme-theoretic intersection $Y_1\cap Y_2\cap\cdots\cap Y_r$, one has
    \[
        \codim_X(W) \leq \sum_{i=1}^r\codim_{X}(Y_i).
    \]
\end{lemma}

\begin{proof}
    Fix an irreducible component $W$ of $Y_1\cap\cdots\cap Y_r$, and let $w \in W$ be a closed point of $X$ that is not in any irreducible component $W'$ of $Y_1\cap\cdots\cap Y_r$ except for $W'=W$. 
    Set $R:=\mathcal{O}_{X,w}$.  
    Since $X$ is smooth, $R$ is a regular local ring.

    Let $\mathfrak{p}_i\subseteq R$ be the prime ideal defining $Y_i$ at $w$, for all $i=1,\dots,r$. 
    We have that $W$ and its reduced structure $W_{\operatorname{red}}$ have the same dimension. 
    Let $\mathfrak{p}\subseteq R$ be the prime ideal defining $W_{\operatorname{red}}$ at $w$. 
    Then, $\mathfrak{p}$ is a minimal prime over $\sum_{i=1}^r\mathfrak{p}_i$. 
    Recall that the codimension of a subvariety of $X$ is equal to the height of its prime ideal in the local ring at any of its closed points (no smoothness needed: reducing to the affine case, using that heights are preserved under localization, and invoking 
    \cite[Theorem~I.1.8A(b), Theorem~I.3.2(c)]{Hartshorne}, we obtain: for any subvariety $Z\subseteq X$ and any closed point $w\in Z$, if $\mathfrak p_Z\subseteq \mathcal{O}_{X,w}$ denotes the prime defining $Z$ at $w$, then $\operatorname{ht}_{\mathcal{O}_{X,w}}(\mathfrak p_Z)=\codim_X(Z)$).  
    Therefore, the inequality in the statement is equivalent to
    \[
        \operatorname{ht}_R(\mathfrak{p})\leq\sum_{i=1}^r\operatorname{ht}_R(\mathfrak{p}_i).
    \]
    We prove this inequality by induction on $r$.    For $r=2$, this is the content of Serre's dimension theorem \cite[Chapter~V, Theorem~3]{Serre-LocalAlgebra}: in a regular ring $R$, if $\mathfrak{p}_1$ and $\mathfrak{p}_2$ are primes, and  $\mathfrak{p}$ is a minimal prime over $\mathfrak{p}_1+\mathfrak{p}_2$, then
    \[
        \operatorname{ht}_R(\mathfrak{p}) \leq \operatorname{ht}_R(\mathfrak{p}_1)+\operatorname{ht}_R(\mathfrak{p}_2).
    \]
    Assume $r\geq 3$ and set $J=\sum_{i=1}^{r-1}\mathfrak{p}_i$. Since $\mathfrak{p}\supseteq J$, we may choose a minimal prime $\mathfrak{q}$ over $J$ such that $\mathfrak{p}\supseteq\mathfrak{q}\supseteq J$. We have that $\mathfrak{p}\supseteq\mathfrak{q}+\mathfrak{p}_r\supseteq J +\mathfrak{p}_r$, and $\mathfrak{p}$ is minimal over $J+\mathfrak{p}_r = \sum_{i=1}^r\mathfrak{p}_i$. Then, $\mathfrak{p}$ is minimal over $\mathfrak{q}+\mathfrak{p}_r$. Applying Serre's Theorem \cite[Chapter~V, Theorem~3]{Serre-LocalAlgebra} once again, we obtain that
    \begin{equation}\label{height.eq.1}
\operatorname{ht}_R(\mathfrak{p})\leq\operatorname{ht}_R(\mathfrak{q})+\operatorname{ht}_R(\mathfrak{p}_r).
    \end{equation}
    By the induction hypothesis applied to the family $\mathfrak{p}_1,\dots,\mathfrak{p}_{r-1}$, for which $\mathfrak{q}$ is minimal over their sum, we have
    \begin{equation}\label{height.eq.2}
        \operatorname{ht}_R(\mathfrak{q}) \leq \sum_{i=1}^{r-1}\operatorname{ht}_R(\mathfrak{p}_i).
    \end{equation}
    Combining Equations~\eqref{height.eq.1} and \eqref{height.eq.2} yields $\operatorname{ht}_R(\mathfrak{p})\leq\sum_{i=1}^r\operatorname{ht}_R(\mathfrak{p}_i)$. This completes the proof. 
\end{proof}

%%%%%%%%%%%%%%%%%%%%%%%%%%%%%%%%%%%%%%%%%%%%%%%%%%%%%%%%%%%%%%%%%%%%%%%

In the following lemma, we use the same tangent-space notation as before for schemes locally of finite type over an algebraically closed field $k$. 

\begin{lemma}[Tangent space of a scheme-theoretic intersection]\label{tangent.intersection} 
Let $X$ be a scheme locally of finite type over an algebraically closed field $k$.
Let $Y_1,\dots,Y_r\subseteq X$ be closed subschemes with scheme-theoretic intersection $Y$.  
Then, for every closed point $p \in Y$ we have
\[
    T_{Y,p}=\bigcap_{i=1}^r T_{Y_i,p}\subseteq T_{X,p}.
\]
\end{lemma}

\begin{proof} 
    Let $Z\subseteq X$ be a closed subscheme containing $p$. 
    In this setting the derivation description of tangent spaces says that 
    \[
        T_{X,p} = \operatorname{Der}_k(\mathcal{O}_{X,p},k)\quad\text{and}\quad T_{Z,p} = \operatorname{Der}_k(\mathcal{O}_{Z,p},k).
    \]
    Denote the ideal sheaf of any closed subscheme $W$ of $X$ by $\mathcal{I}_W$. 
    We have that $\mathcal{I}_{Y} = \mathcal{I}_{Y_1} + \cdots  + \mathcal{I}_{Y_r}$. 
    A derivation $\mathcal{O}_{Z,p}=\mathcal{O}_{X,p}/\mathcal{I}_{Z,p}\to k$ corresponds to a derivation $\mathcal{O}_{X,p}\to k$ that vanishes on $\mathcal{I}_{Z,p}\subseteq\mathcal{O}_{X,p}$. This identifies $T_{Z,p}$ with the $k$-subspace 
    \[
        \{D\in T_{X,p}\,\vert\, D(\mathcal{I}_{Z,p}) =0\}\subseteq T_{X,p}.
    \]
    Taking $Z=Y$, we have that $\mathcal{I}_{Y,p} = \mathcal{I}_{Y_1,p} + \cdots  + \mathcal{I}_{Y_r,p}$, which implies that
    \begin{align*}
        T_{Y,p} &= \{D\in T_{X,p}\,\vert\, D(\mathcal{I}_{Y,p}) =0\}
                = \{D\in T_{X,p}\,\vert\, D(\mathcal{I}_{Y_1,p} + \cdots + \mathcal{I}_{Y_r,p}) =0\}\\
                &= \bigcap_{i=1}^r \{D\in T_{X,p}\,\vert\, D(\mathcal{I}_{Y_i,p}) =0\} 
                = \bigcap_{i=1}^r T_{Y_i,p} \subseteq T_{X,p},
    \end{align*}
    as claimed.
\end{proof}

\begin{lemma}\label{lemma.codimension.transversality.two.subvarieties}
Let $X$ be a smooth variety and let $Y_1,Y_2$ be smooth subvarieties of $X$ such that $Y_1\cap Y_2$ is a nonempty smooth subvariety. Then, $Y_1$ and $Y_2$ intersect transversally on $X$ if and only if
\[
    \codim_X(Y_1\cap Y_2)
    = \codim_X(Y_1)+\codim_X(Y_2).
\]
\end{lemma}

\begin{proof}
For every $p\in Y_1\cap Y_2$, Lemma~\ref{tangent.intersection} gives $T_{Y_1\cap Y_2,p}=T_{Y_1,p}\cap T_{Y_2,p}$. Since $X$ and $Y_1\cap Y_2$ are smooth, the codimension equality in the statement is equivalent to 
\[
\codim(T_{Y_1,p}\cap T_{Y_2,p},T_{X,p}) = \codim_X(Y_1)+\codim_X(Y_2).
\]
By Lemma~\ref{characterizations.transverse.intersection}, this is equivalent to transversality at $p$. At every point outside $Y_1\cap Y_2$, at most one of the two smooth subvarieties contains the point, so transversality is automatic. Since $Y_1\cap Y_2$ is nonempty, the result follows.
\end{proof}

The analogue of Lemma~\ref{lemma.codimension.transversality.two.subvarieties} need not hold for more than two smooth subvarieties. 
For example, the hypersurfaces defined by $y=0$, $x=1$, and $y=xz$ in $\mathbb{A}^3$ have the reduced point $(1,0,0)$ as their intersection, but do not intersect transversally at the origin $(0,0,0)$.
However, the result holds for members of an arrangement, as the following lemma shows.

\begin{lemma}\label{lemma.codimension.characterization.transversal.intersections}
Let $X$ be a smooth variety and let $Y_1,\dots,Y_r$ be members of an arrangement on $X$ with nonempty intersection. Then, $Y_1,\dots,Y_r$ intersect transversally on $X$ if and only if
\begin{equation}\label{equation.codimension.criterion.arrangement}
    \codim_X(Y_1\cap\cdots\cap Y_r)=\sum_{i=1}^r\codim_X(Y_i).
\end{equation}
\end{lemma}

\begin{proof}
Set $Z:=Y_1\cap\cdots\cap Y_r$. For every $p\in Z$, Lemma~\ref{tangent.intersection} gives $T_{Z,p}=\bigcap_{i=1}^r T_{Y_i,p}$. Since $X$ and $Z$ are smooth, we have
\[
    \codim\left(\bigcap_{i=1}^r T_{Y_i,p},T_{X,p}\right)=\codim(T_{Z,p},T_{X,p})=\codim_X(Z).
\]
Thus, by Lemma~\ref{characterizations.transverse.intersection}, transversality at $p$ is equivalent to Equation~\eqref{equation.codimension.criterion.arrangement}.
In particular, transversality on $X$ implies Equation~\eqref{equation.codimension.criterion.arrangement}.

Conversely, assume Equation~\eqref{equation.codimension.criterion.arrangement} holds, and fix $p\in Z$. By the preceding argument, $Y_1,\dots,Y_r$ intersect transversally at $p$. For any nonempty subset $I\subseteq\{1,\dots,r\}$, set $Z_I:=\bigcap_{i\in I}Y_i$. Since $Z\subseteq Z_I$, the intersection $Z_I$ is a nonempty member of the arrangement and hence is smooth and irreducible. The subfamily $\{Y_i\}_{i\in I}$ intersects transversally at $p$, so the same tangent-space calculation gives
\[
    \codim_X(Z_I)=\sum_{i\in I}\codim_X(Y_i).
\]
Applying that calculation at any point of $Z_I$ now shows that $\{Y_i\}_{i\in I}$ intersects transversally at every point of $Z_I$.

Finally, given $q\in X$, let $I(q):=\{i\mid q\in Y_i\}$. If $I(q)$ is nonempty, then $q\in Z_{I(q)}$, and the preceding argument gives transversality at $q$ of the members containing $q$. The remaining members impose no condition on transversality there, by Definition~\ref{definition.transversal.intersection}. If $I(q)$ is empty, transversality at $q$ is automatic. Therefore $Y_1,\dots,Y_r$ intersect transversally on $X$.
\end{proof}

%%%%%%%%%%%%%%%%%%%%%%%%%%%%%%%%%%%%%%%%%%%%%%%%%%%%%%%%%%%%%%%%%%%%%%%

\begin{lemma}
\label{transversality.within.building.set}
Let $X$ be a smooth variety and let $\mathcal{G}$ be a building set on $X$. If $Y_1$ and $Y_2$ are in $\mathcal{G}$, but $Y_1\cap Y_2 \neq \emptyset$ is not in $\mathcal{G}$, then $Y_1$ and $Y_2$ are the $\mathcal{G}$-factors of $Y_1\cap Y_2$. 
In particular, $Y_1$ and $Y_2$ intersect transversally.    
\end{lemma}
\begin{proof}
    Since both $Y_1,Y_2\in\mathcal{G}$ contain $Y_1\cap Y_2$, then there exist $Y_1',Y_2'\in\mathcal{G}$ minimal among the elements of $\mathcal{G}$ that contain $Y_1\cap Y_2$, such that $Y_1'\subseteq Y_1$ and $Y_2'\subseteq Y_2$.       
    Notice that $Y_1'\cap Y_2'\supseteq Y_1\cap Y_2$, then we have 
    \[
        Y_1\cap Y_2 = Y_1'\cap Y_2' = Y_1'\cap Y_2 = Y_1\cap Y_2'.
    \]
    Notice also that $Y_1'\neq Y_2'$ because otherwise $Y_1' = Y_2' = Y_1'\cap Y_2' = Y_1\cap Y_2 \notin \mathcal{G}$ which is a contradiction. 
    Hence, $Y_1'$ and $Y_2'$ are transverse because they are distinct $\mathcal{G}$-factors of $Y_1\cap Y_2$. Then, $\codim(Y_1'\cap Y_2') = \codim(Y_1')+\codim(Y_2')$ by Lemma~\ref{lemma.codimension.transversality.two.subvarieties}. It follows from Lemma~\ref{subadditivity.codimension.smooth.ambient} that
    \begin{align*}
        \codim(Y_1') + \codim(Y_2') &= \codim(Y_1'\cap Y_2')\\
                                    &= \codim(Y_1'\cap Y_2) \leq \codim(Y_1')+\codim(Y_2).
    \end{align*}
    Then, $\codim(Y_2')\leq\codim(Y_2)$, which implies that $Y_2'=Y_2$. Similarly, we have that $Y_1'=Y_1$. 
    Therefore, $Y_1$ and $Y_2$ are $\mathcal{G}$-factors of $Y_1\cap Y_2$. 
    If there were additional $\mathcal{G}$-factors of $Y_1\cap Y_2$, by Lemma~\ref{lemma.codimension.characterization.transversal.intersections} the sum of the codimensions of the additional $\mathcal{G}$-factors would have to be zero, and since such codimensions are positive integers, we conclude that $Y_1$ and $Y_2$ are the only $\mathcal{G}$-factors of $Y_1\cap Y_2$.
    We deduce additionally that $Y_1$ and $Y_2$ are transverse by the definition of a building set. 
\end{proof} 

\begin{corollary}
    \label{injectivity.of.restrictions.within.building.set}
Let $X$ be a smooth variety and let $\mathcal{G}$ be a building set on $X$. If $Y$, $Y_1$ and $Y_2$ in $\mathcal{G}$ are such that $Y \cap Y_1 = Y \cap Y_2 \neq \emptyset$ and this intersection is not in $\mathcal{G}$, then $Y_1 = Y_2$.
\end{corollary}
\begin{proof}
    By Lemma~\ref{transversality.within.building.set}, $Y$ and $Y_1$ 
    are the $\mathcal{G}$-factors of $Y \cap Y_1$.
    Similarly, $Y$ and $Y_2$  
    are the $\mathcal{G}$-factors of $Y \cap Y_2$.
    Since $Y \cap Y_1 = Y \cap Y_2$, then $\{Y,Y_1\} = \{Y,Y_2\}$. 
    Since $Y \cap Y_1 = Y \cap Y_2$ is not in $\mathcal{G}$, then $Y \neq Y_1$ and $Y \neq Y_2$, and hence $Y_1 = Y_2$ as desired. 
\end{proof}

\begin{example}
    In general, elements of a building set do not need to intersect transversally. For example, let $X=\mathbb{A}^3$ with coordinates $x,y,z$. Let $L_1=V(x,y)$, $L_2=V(x,z)$ and $p=V(x,y,z)$. Then, $\mathcal{G}=\{p,L_1,L_2\}$ is a building set on $X$, but $L_1$ and $L_2$ do not intersect transversally on $X$, since $T_{L_1,p}^\perp + T_{L_2,p}^\perp$ is not a direct sum in $T_{X,p}^*$.
\end{example}

%%%%%%%%%%%%%%%%%%%%%%%%%%%%%%%%%%%%%%%%%%%%%%%%%%%%%%%%%%%%%%%%%%%%%%%

\subsection{Chow rings and blow-ups}

We recall some results on the Chow rings of blow-ups that we will use to study wonderful compactifications through their iterated blow-up constructions. We begin by introducing the notation. 

We denote by $A^\bullet(X)$ the Chow ring of a smooth variety $X$. For a morphism $f:X'\to X$ between smooth varieties, we write $f^*:A^\bullet(X)\to A^\bullet(X')$ for the pullback. We also denote by $f^*$ its coefficientwise extension to polynomial rings $f^*:A^\bullet(X)[t]\to A^\bullet(X')[t]$, fixing $t$.

For a blowup $\pi:\widetilde X\to X$ of a smooth variety along a smooth subvariety, the pullback $\pi^*:A^\bullet(X)\to A^\bullet(\widetilde X)$ is injective by \cite[Proposition~6.7(b)]{Fulton-book}, and we identify classes in $A^\bullet(X)$ with their pullbacks. Likewise, we identify polynomials in $A^\bullet(X)[t]$ with their images in $A^\bullet(\widetilde X)[t]$ under the coefficientwise extension of $\pi^*$.

For an inclusion $i:Z\hookrightarrow X$ of  smooth varieties, we write $\alpha|_Z:=i^*(\alpha)$ for $\alpha\in A^\bullet(X)$ or $\alpha\in A^\bullet(X)[t]$. We define
\[
J_Z:=\ker\big(i^*:A^\bullet(X)\to A^\bullet(Z)\big).
\]
We use the notation $J_{X,Z}=J_Z$ when we want to make the ambient space explicit.

\begin{definition}[Chern polynomials] \label{definition.Chern.polynomial}
    Let $X$ be a smooth variety and $Z$ be a smooth subvariety of codimension $d$.
    A \emph{Chern polynomial for $Z$ in $X$} is a monic degree-$d$ polynomial
    \[
      P_{X,Z}(t) = t^d + a_1 t^{d-1} + \cdots + a_d \in A^{\bullet}(X)[t]
    \]
    with $a_i\in A^i(X)$, constant term $a_d=[Z]$, and $a_i|_Z=c_i(N_{Z/X})$ for $1\leq i\leq d$, where $N_{Z/X}$ is the normal bundle of $Z$ in $X$. Chern polynomials are guaranteed to exist if the restriction map $i^*$ is surjective, but need not be unique. The polynomial $P_{X,Z}(t)$ is related to the Chern polynomial $c_t(N_{Z/X})$ of the normal bundle $N_{Z/X}$ by 
\[
t^d\,P_{X,Z}(1/t)\big|_Z=c_t(N_{Z/X}).
\]
  \end{definition}

\begin{example}
    If $Z = D$ is a smooth divisor on $X$, then $d=1$ and we may take $P_{X,D}(t) = t + D$.
\end{example}

This notation for Chern polynomials is standard. However, we introduce new notation as well.

\begin{definition} \label{definition.Chern.polynomial.triple}
    Let $X$ be a smooth variety and $Z\subsetneq Y$ be smooth subvarieties of $X$. 
    A \emph{Chern polynomial for} $Z \subseteq Y$ in $X$ is any polynomial $P_{X,Y,Z}(t) \in A^\bullet(X)[t]$ that maps to a Chern polynomial $P_{Y,Z}(t)\in A^\bullet(Y)[t]$ under the homomorphism induced by the restriction $A^\bullet(X)\to A^\bullet(Y)$.
    A Chern polynomial $P_{X,Y,Z}(t)$ is guaranteed to exist if both restriction maps $A^\bullet(X)\to A^\bullet(Y)$ and $A^\bullet(Y)\to A^\bullet(Z)$ are surjective, but need not be unique.  
\end{definition}

With this terminology in place, we now state Keel's result describing the Chow ring of a blow-up.

\begin{lemma}[Keel's blow-up formula {\cite[Appendix, Theorem~1]{Keel-thesis}}] \label{lemma: keel's blowup formula}
Let $X$ be a smooth variety and $Z \subseteq X$ a smooth subvariety, and suppose that $A^{\bullet}(X) \to A^{\bullet}(Z)$ is surjective. Let $\pi : \widetilde{X} = \Bl_Z X \to X$ be the blow-up of $X$ along $Z$ with exceptional divisor $E$, and fix a Chern polynomial $P_{X,Z}(t)$ for $Z$ in $X$. 
Then, pullback by $\pi$ on coefficients and sending the variable $E$ to the class of the exceptional divisor induce an isomorphism of $A^{\bullet}(X)$-algebras 
\[
  A^{\bullet}(\widetilde{X})\cong \frac{A^{\bullet}(X)[E]}{\langle J_Z \cdot E,\ P_{X,Z}(-E)\rangle}.
\]
\end{lemma}

Keel's formula expresses the Chow ring of a blow-up in terms of a Chern polynomial and the kernel of a restriction map. For the iterated blowup construction of wonderful compactifications, we therefore need to understand how these data change under blow-ups. Before addressing this, we record two results on Chern polynomials for transverse intersections.

\begin{lemma}[{\cite[Lemma~5.1(b)]{Fulton-MacPherson}}]\label{lemma.chern.polynomial.identities}
Let $X$ be a smooth variety, and let $Z,W\subseteq X$ be smooth subvarieties intersecting transversally, with $Z \cap W$ nonempty and irreducible, and having Chern polynomials $P_{X,Z}(t)$ and $P_{X,W}(t)$. 
Then, 
\begin{enumerate}[label=(\roman*),ref=(\roman*)]
  \item \label{lemma.chern.polynomial.identities.i} $P_{X,Z}(t)\cdot P_{X,W}(t)$ is a Chern polynomial of $Z\cap W$ in $X$. 
  \item \label{lemma.chern.polynomial.identities.ii} The image of $P_{X,Z}(t)$ under the restriction map $A^\bullet(X)[t]\to A^\bullet(W)[t]$ is a Chern polynomial of $Z\cap W$ in $W$. 
\end{enumerate}
\end{lemma}

\begin{proof}
Set $K:=Z\cap W$, 
$d_Z:=\codim_X(Z)$, $d_W:=\codim_X(W)$, and $d_K:=\codim_X(K)=d_Z+d_W$.

Recall that $P_{X,Z}(t)$ is characterized by: it is monic of degree $d_Z$ with constant term $[Z]$,
and upon restriction to $Z$ one has the reversal identity
\begin{equation}\label{eq:reversal}
t^{d_Z}\,P_{X,Z}(1/t)\big|_{Z}=c_t(N_{Z/X}).
\end{equation}
The same applies to $W$ and $K$.

\begin{enumerate}[label=(\roman*),ref=(\roman*)]
\item Let $Q(t):=P_{X,Z}(t)\,P_{X,W}(t)$. Then $Q(t)$ is monic of degree $d_K$, and its constant term is
$[Z]\cdot[W]=[K]$ by transversality. Restricting \eqref{eq:reversal} to $K\subseteq Z$ and $K\subseteq W$ gives
\[
t^{d_Z}P_{X,Z}(1/t)\big|_{K}=c_t(N_{Z/X}|_{K}),
\qquad
t^{d_W}P_{X,W}(1/t)\big|_{K}=c_t(N_{W/X}|_{K}).
\]
Multiplying and using $d_K=d_Z+d_W$ yields
\[
t^{d_K}Q(1/t)\big|_{K}
=
c_t(N_{Z/X}|_{K})\;c_t(N_{W/X}|_{K}).
\]
For a transversal intersection, $N_{K/X}\cong N_{Z/X}|_{K}\oplus N_{W/X}|_{K}$, so by the Whitney sum formula
\[
c_t(N_{K/X}) 
=c_t(N_{Z/X}|_{K})\;c_t(N_{W/X}|_{K})
=t^{d_K}Q(1/t)\big|_{K}.
\]
Thus $Q(t)$ satisfies the defining properties of a Chern polynomial of $K$ in $X$.

\item Consider $\big(P_{X,Z}(t)\big)\big|_{W}\in A^\bullet(W)[t]$. It is monic of degree
$d_Z=\codim_W(K)$ and its constant term is $[Z]|_W=[K]$ in $A^\bullet(W)$ by transversality. 
Moreover, restricting \eqref{eq:reversal} further to $K\subseteq Z$ gives
\[
t^{d_Z}\,\big(P_{X,Z}(1/t)\big)\big|_{K}
=c_t(N_{Z/X})|_{K}
=c_t(N_{Z/X}|_{K}).
\]
By transversality $N_{Z/X}|_{K} \cong N_{K/W}$, hence
\[
t^{d_Z}\,\big(P_{X,Z}(1/t)\big)\big|_{K}=c_t(N_{K/W}).
\]
Therefore $\big(P_{X,Z}(t)\big)\big|_{W}$ is a Chern polynomial for $K = Z\cap W$ in $W$. \qedhere
\end{enumerate}
\end{proof}
We obtain the following useful corollary.
\begin{corollary}\label{corollary.chern.polynomial.remove.transversal.factor}
Let $V,W$ and $Y$ be smooth subvarieties of a smooth variety $X$ such that $W\subsetneq V$; $Y$ is transversal to both $V$ and $W$; and $V\cap Y$ and $W\cap Y$ are irreducible and nonempty. 
Then, the image of a Chern polynomial $P_{V,W}(t)$ under the restriction map $A^\bullet(V)[t]\to A^\bullet(V\cap Y)[t]$ is a Chern polynomial for $W\cap Y$ in $V\cap Y$.
\end{corollary}
\begin{proof}
    This is a direct consequence of Lemma~\ref{lemma.chern.polynomial.identities}\ref{lemma.chern.polynomial.identities.ii}, provided we can prove that $Y\cap V$ and $W$ intersect transversally as subvarieties of $V$.
    
    By Example~\ref{example.transversality.two.subvarieties}, it is enough to prove that $T_p W + T_p(V\cap Y) = T_p V$ for every $p\in W\cap Y$. It is clear that the left hand side is contained in $T_p V$. Conversely, let $v\in T_p V$. Since $W$ and $Y$ intersect transversally on $X$, then $T_p V \subseteq T_p X = T_p W + T_p Y$, and hence there are $w\in T_pW$ and $y\in T_pY$ such that $v=w+y$. This implies that $y = v - w \in T_pV$, as $w\in T_pW\subseteq T_pV$. In particular, $y\in T_pV \cap T_pY = T_p(V\cap Y)$.
    Hence $v = w + y \in T_pW+T_p(V\cap Y)$, as desired. 
\end{proof}
We can now describe how Chern polynomials behave under blow-ups.
\begin{lemma}[{Chern polynomials under blow-up \cite[Lemma~5.2]{Fulton-MacPherson}}]\label{lemma.Chern.polynomials.blowup} 
    Let $X$ be smooth, $Z \subseteq X$ a smooth subvariety, and $\widetilde{X} = \Bl_Z X$ with exceptional divisor $E$. Let $W \subseteq X$ be a smooth subvariety not contained in $Z$, and $\widetilde{W}$ its strict transform. 
    Suppose that $P_{X,W}(t)$ is a Chern polynomial for $W$ in $X$.
    \begin{enumerate}[label=(\roman*),ref=(\roman*)] 
        \item If $W$ and $Z$ intersect transversally, then $P_{X,W}(t)$ is a Chern polynomial for $\widetilde{W}$ in $\widetilde{X}$.
        \item If $Z \subsetneq W$, then $P_{X,W}(t-E)$ is a Chern polynomial for $\widetilde{W}$ in $\widetilde{X}$.
    \end{enumerate}
\end{lemma}

The following lemma presents a strict transform as a blow-up. Its proof simply reconciles this version of the statement with versions available in the literature. 

\begin{lemma}[Strict transforms as blow ups]\label{lemma.strict.transforms.as.blowups}
Let $X$ be a variety and let $Y$ and $Z$ be subvarieties of $X$. 
Assume that $Y$ is not contained in $Z$. 
Then, the strict transform of $Y$ in the blow up of $X$ along $Z$ is canonically isomorphic to the blow up of $Y$ along the scheme-theoretic intersection $Y \cap Z$.
\end{lemma}

\begin{proof}
Let $i:Y\hookrightarrow X$ be the closed immersion and let $\mathcal{I} \subseteq \mathcal{O}_X$ be the ideal sheaf of $Z$.
Set $\mathcal{J}:=i^{-1}\mathcal{I}\cdot\mathcal{O}_Y$, which is the ideal sheaf of the scheme-theoretic intersection $Y\cap Z$ in $Y$.
By \cite[II, Cor.~7.15]{Hartshorne}, there is a unique morphism $\widetilde{i}:\ \operatorname{Bl}_{\mathcal{J}}Y\ \longrightarrow\ \operatorname{Bl}_{\mathcal{I}}X$ 
making the natural diagram commute, and since $i$ is a closed immersion, $\widetilde{i}$ is a closed immersion.
Hartshorne (loc.\ cit.) \emph{defines} the strict transform of $Y$ under the blow-up $\operatorname{Bl}_{\mathcal{I}}X\to X$ to be the image $\widetilde{i}\big(\operatorname{Bl}_{\mathcal{J}}Y\big)\subseteq \operatorname{Bl}_{\mathcal{I}}X$.
Therefore the strict transform of $Y$ is canonically isomorphic to $\operatorname{Bl}_{\mathcal{J}}Y=\operatorname{Bl}_{Y\cap Z}Y$, as claimed.
\end{proof}

While Lemma~\ref{lemma.Chern.polynomials.blowup} explains how Chern polynomials change upon blowing up (at least in the two cases that are relevant for our present work), the next result describes the ideal $J_{\widetilde{X},\widetilde{W}}\subseteq A^\bullet(\widetilde{X})$ corresponding to the strict transform $\widetilde{W}$ of a subvariety $W$ of $X$. 

\begin{lemma}\label{lemma.Jideal.blowup} 
    Let $X$ be smooth, $Z \subseteq X$ a smooth subvariety, and $\widetilde{X} = \Bl_Z X$ with exceptional divisor $E$. 
    Let $W \subseteq X$ be a smooth subvariety not contained in $Z$, 
    such that $W\cap Z$ is a smooth variety, 
    and let $\widetilde{W}$ be its strict transform. 
    Suppose that the pullback maps $A^{\bullet}(X) \to A^{\bullet}(Z)$, $A^{\bullet}(X) \to A^{\bullet}(W)$, and $A^\bullet(W)\to A^\bullet(W\cap Z)$ are surjective. 
    Let $P_{X,W,W\cap Z}(t) \in A^\bullet(X)[t]$ be a Chern polynomial for $W\cap Z \subseteq W$ in $X$. 
    Then, the restriction map $A^\bullet(\widetilde{X})\to A^\bullet(\widetilde{W})$ is surjective and its kernel $J_{\widetilde{X},\widetilde{W}}$ is given by 
    \[
        J_{\widetilde{X},\widetilde{W}} = \langle \, J_{X,W}  \, , \, J_{X,W\cap Z}\cdot E \, , \, P_{X,W,W\cap Z}(-E) \, \rangle.
    \]
\end{lemma}

\begin{proof}
    Consider the following commutative diagram
    \[
    \begin{tikzcd}
      A^\bullet(X)[E] \ar[d,"\pi_X",swap] \ar[r, "i^*"] & A^\bullet(W)[e] \ar[r,"\pi_W"] & A^\bullet(\widetilde{W}) \\
      A^\bullet(\widetilde{X}),\ar[rru,"j^*", swap]
    \end{tikzcd}
    \]
    where the top arrow extends the restriction map induced by $i:W\hookrightarrow X$ by sending $E$ to $e$, and $j^*$ is induced by the inclusion $j:\widetilde{W}\hookrightarrow\widetilde{X}$. By Lemma~\ref{lemma.strict.transforms.as.blowups}, $\widetilde{W}\cong\operatorname{Bl}_{W\cap Z}W$, and $\pi_X$ and $\pi_W$ are the quotient maps given by Keel's blow-up formula: 
    \[
        A^\bullet(\widetilde{X}) = \frac{A^\bullet(X)[E]}{\langle \, J_{X,Z}\cdot E, P_{X,Z}(-E) \, \rangle },
        \quad\text{and}\quad
        A^\bullet(\widetilde{W}) = \frac{A^\bullet(W)[e]}{\langle \, J_{W,W\cap Z}\cdot e,P_{W,W\cap Z}(-e) \, \rangle },
    \]
    where $P_{X,Z}(t) \in A^\bullet(X)[t]$ and $P_{W,W\cap Z}(t) \in A^\bullet(W)[t]$ are the Chern polynomials, and $E$ and $e =\widetilde{W}\cap E$ are the classes of the exceptional divisors in $\widetilde{X}$ and $\widetilde{W}$, respectively.
    Choose $P_{W,W\cap Z}(t)$ to be the restriction of $P_{X,W,W\cap Z}(t)$ to $A^\bullet(W)[t]$. 
    The surjectivity of $A^\bullet(\widetilde{X})\to A^\bullet(\widetilde{W})$ follows from the previous discussion since $i^*$ is surjective by hypothesis and $\pi_W$ by construction.

    Since $\pi_{X}$ is surjective it follows that $J_{\widetilde{X},\widetilde{W}} := \ker(j^*) = \pi_X\left(\ker(\pi_W\circ i^*)\right)$. We have 
    \[
        \ker(\pi_W\circ i^*) = (i^*)^{-1}(\langle J_{W,W\cap Z}\cdot e,P_{W,W\cap Z}(-e)\rangle) = (i^*)^{-1}(\langle J_{W,W\cap Z}\cdot e \rangle )+(i^*)^{-1}( \langle P_{W,W\cap Z}(-e) \rangle ),
    \]
    where the last equality follows from $i^*$ being surjective.

    Since $i^*$ is surjective, note that $(i^*)^{-1}(\langle P_{W,W\cap Z}(-e) \rangle ) = \langle J_{X,W} \rangle + \langle P_{X,W,W\cap Z}(-E) \rangle $. On the other hand, by functoriality of pullbacks applied to the composition $A^\bullet(X)\to A^\bullet(W)\to  A^\bullet(W\cap Z)$, 
    \[
        (i^*)^{-1}(J_{W,W\cap Z})= \ker(A^\bullet(X)\to A^\bullet(W\cap Z)) = J_{X,W\cap Z}.
    \] 
    Since $i^*$ is surjective we deduce 
    \[
    (i^*)^{-1}( \langle J_{W,W\cap Z}\cdot e \rangle )= \langle J_{X,W\cap Z}\cdot E \rangle + \operatorname{ker}(i^*) 
    = \langle J_{X,W\cap Z}\cdot E \rangle +  \langle J_{X,W} \rangle. 
    \]
    Together, these equalities give the desired conclusion:   
    \[
        J_{\widetilde{X},\widetilde{W}} = 
            \langle \, J_{X,W}  \, , \, J_{X,W\cap Z}\cdot E \, , \, P_{X,W,W\cap Z}(-E) \, \rangle. \qedhere 
    \]
\end{proof}

Lemma~\ref{lemma.Jideal.blowup} is a generalization of \cite[Lemma~5.5(1)]{Fulton-MacPherson} and \cite[Lemma~2.2]{petersen2017chow} (itself a correction of \cite[Lemma~5.4]{Fulton-MacPherson}); for the latter, we omit the elementary case $W\cap Z=\emptyset$ for simplicity and state explicitly the assumption that $A^\bullet(W)\to A^\bullet(W\cap Z)$ is surjective.

\begin{corollary}\label{lemma.Jideal.known.cases}
    In the setting of Lemma~\ref{lemma.Jideal.blowup}:
    \begin{enumerate}[label=(\roman*),ref=(\roman*)]
        \item \label{lemma.Jideal.known.cases.item.i} \cite[Lemma~2.2]{petersen2017chow} If $W$ and $Z$ intersect transversally, then 
        \[ J_{\widetilde{X},\widetilde{W}} = \langle \, J_{X,W} \, , \, J_{X,W\cap Z}\cdot E \, \rangle. \]
        \item \label{lemma.Jideal.known.cases.item.ii} \cite[Lemma~5.5(1)]{Fulton-MacPherson} If $Z \subsetneq W$, then 
        \[ J_{\widetilde{X},\widetilde{W}} = \langle \, J_{X,W} \, , \, P_{X, W,Z}(-E) \, \rangle. \]
    \end{enumerate}
\end{corollary}

\begin{proof}
    By Lemma~\ref{lemma.Jideal.blowup}, the ideal $J_{\widetilde{X},\widetilde{W}}$ is generated by $J_{X,W}$, $J_{X,W\cap Z}\cdot E$, and $P_{X,W,W\cap Z}(-E)$.
    \begin{enumerate}[label=(\roman*),ref=(\roman*)]
    \item Since $W$ and $Z$ intersect transversally, by Lemma~\ref{lemma.chern.polynomial.identities}\ref{lemma.chern.polynomial.identities.ii} we may choose $P_{X,W,W\cap Z}(t):=P_{X,Z}(t)$.  
    By Keel's formula $P_{X,Z}(-E) = 0$ in $A^\bullet(\widetilde{X})$, and then $J_{\widetilde{X},\widetilde{W}}$ is generated by $J_{X,W}$ and $J_{X,W\cap Z}\cdot E$. 
    \item Since $W \cap Z = Z$, we have   
    $J_{X,W\cap Z}\cdot E = J_{X, Z}\cdot E$ and $P_{X,W,W\cap Z}(-E) = P_{X, W, Z}(-E)$. 
    By Keel's formula $J_{X,Z} \cdot E = 0$ in $A^\bullet(\widetilde{X})$, and then $J_{\widetilde{X},\widetilde{W}}$ is generated by $J_{X,W}$ and $P_{X,W, Z}(-E)$. \qedhere 
    \end{enumerate}
\end{proof}

%%%%%%%%%%%%%%%%%%%%%%%%%%%%%%%%%%%%%%%%%%%%%%%%%%%%%%%%%%%%%%%%
%%%%%%%%%%%%%%%%%%%%%%%%%%%%%%%%%%%%%%%%%%%%%%%%%%%%%%%%%%%%%%%%
%%%%                                                        %%%%
%%%%        3. Restriction of building sets    
%%%%                                                        %%%%
%%%%%%%%%%%%%%%%%%%%%%%%%%%%%%%%%%%%%%%%%%%%%%%%%%%%%%%%%%%%%%%%
%%%%%%%%%%%%%%%%%%%%%%%%%%%%%%%%%%%%%%%%%%%%%%%%%%%%%%%%%%%%%%%%

\section{Restriction of building sets} \label{section.restriction}

In this section we introduce the restriction of a building set and study its main properties. 
Given a smooth variety $X$, a building set $\mathcal{G}$ on $X$, and an element $Y$ in $\mathcal{G}$, in Section~\ref{subsection.restriction.building.sets} we show that the restriction $\mathcal{G}|_Y$ of $\mathcal{G}$ to $Y$ is a building set on $Y$ (see Theorem~\ref{theorem.restriction.of.building.set.is.building.set}). 
Next, in Section~\ref{subsection.strict.transform.is.wonderful} we show that if $Y$ is maximal in $\mathcal{G}$ then $\mathcal{G} \smallsetminus \{Y\}$ is a building set on $X$, and the wonderful compactification of $\mathcal{G}|_Y$ on $Y$ is canonically isomorphic to the strict transform of $Y$ in the wonderful compactification of $\mathcal{G} \smallsetminus \{Y\}$ on $X$ (see Theorem~\ref{strict.transforms.are.wonderful}). 
Lastly, in Section~\ref{subsection.restricted.building.set.iterated.blowup} we show that, when $Y$ is maximal and $\mathcal{G}$ is ordered compatibly with containments, the ambient iterated blowup construction induces the wonderful compactification construction on $Y$ in the corresponding order.

%%%%%%%%%%%%%%%%%%%%%%%%%%%%%%%%%%%%%%%%%%%%%%%%%%%%%%%%%%%%%%%%

\subsection{The restriction of a building set} \label{subsection.restriction.building.sets}

\begin{definition}[Restriction of a building set] \label{definition.restriction.building.set}
    Let $X$ be a smooth variety, $\mathcal{G}$ be a building set on $X$, and $Y$ be an element of $\mathcal{G}$. 
    The \emph{restriction} $\mathcal{G}|_{Y}$ of the building set $\mathcal{G}$ to $Y$ is defined to be the following (possibly empty) collection of proper closed subsets of $Y$,
    \[
    \mathcal{G}|_{Y} := \left\{ Y \cap Z \ | \ Z \in \mathcal{G} \textnormal{ such that } \emptyset \neq Y \cap Z \subsetneq Y \right\}. 
    \]
\end{definition}

Although one can formally restrict a building set $\mathcal{G}$ to any subvariety $Y$ by the rule of Definition~\ref{definition.restriction.building.set}, the resulting family $\mathcal{G}|_Y$ need not be a building set in general. 
For this reason, throughout this article we only consider restrictions to members $Y\in\mathcal{G}$ as in Definition~\ref{definition.restriction.building.set}, and with this convention the restriction is again a building set, as the next theorem shows.

\begin{theorem}[Restriction of a building set is a building set] \label{theorem.restriction.of.building.set.is.building.set}
Let $X$ be a smooth variety and let $\mathcal{G}$ be a building set on $X$. Then, for every $Y\in\mathcal{G}$, the restriction $\mathcal{G}|_Y$ is a building set on $Y$. 
\end{theorem}

\begin{proof} 
Let $Y_1,\dots,Y_{n-1},Y_n$ be the distinct elements of $\mathcal{G}$, and without loss of generality we take $Y=Y_n$.
Let $\mathcal{A}$ be the arrangement on $X$ induced by $\mathcal{G}$, and $\mathcal{A}_Y$ be the collection of all nonempty intersections of elements in $\mathcal{G}|_Y$. Notice that $\mathcal{A}_Y=\{Z\in\mathcal{A}\,\vert\,Z\subsetneq Y\}$. Hence, $\mathcal{A}_Y$ is an arrangement on $Y$, since it is a collection of smooth subvarieties of $Y$ that is closed under nonempty scheme-theoretic intersection. 

To complete the argument we need to show that every $S\in\mathcal{A}_Y$ is equal to the transverse intersection in $Y$ of the minimal elements of $\mathcal{G}|_Y$ containing $S$. 

Fix $S\in\mathcal{A}_Y$. Then, $S\in\mathcal{A}$ as well, so it is the transverse intersection of the minimal elements in $\mathcal{G}$ containing $S$, say $Y_{i_1},\dots,Y_{i_r}$. In particular,
\begin{equation} \label{equation.codimension.sum}
    \codim_X(S) = \sum_{\ell=1}^r\codim_X(Y_{i_\ell}). 
\end{equation}
We consider two cases depending on whether $Y_n$ is minimal over $S$:
\smallskip

\noindent\textbf{Case 1:} Suppose that $Y\in\{Y_{i_1},\dots,Y_{i_r}\}$, say $Y = Y_{i_r}$. 
Thus, $Y_{i_1}\cap Y,\dots,Y_{i_{r-1}}\cap Y$ are in $\mathcal{G}|_Y$. We have 
\begin{equation}\label{equation.S.two.presentations}
    S \ = \  Y_{i_1}\cap Y_{i_2}\cap\cdots\cap Y_{i_{r-1}}\cap Y 
      \ = \ (Y_{i_1}\cap Y) \cap \cdots\cap (Y_{i_{r-1}}\cap Y),
\end{equation}
and $Y_{i_1}\cap Y,\dots,Y_{i_{r-1}}\cap Y$ intersect transversally in $Y$, by Lemma~\ref{lemma.properties.transversal.intersections}\ref{lemma.properties.transversal.intersections.item.restrictions}. We are then reduced to showing that $Y_{i_1}\cap Y,\dots,Y_{i_{r-1}}\cap Y$ are precisely the minimal elements of $\mathcal{G}|_Y$ containing $S$.

Let $Y'\cap Y$ be a minimal element of $\mathcal{G}|_Y$ containing $S$, where $Y'\in\mathcal{G}$. Since $S\subseteq Y'$, then there exists $j\in\{1,\dots,r\}$ such that $Y_{i_j}\subseteq Y'$. Then, $S\subseteq Y_{i_j}\cap Y\subseteq Y'\cap Y$, and by minimality $Y'\cap Y = Y_{i_j}\cap Y$. We see that $j\neq r$, because $Y_{i_j}\subseteq Y'$, but $Y_{i_r}=Y\not\subseteq Y'$ as $Y'\cap Y\in\mathcal{G}|_Y$. It follows that the minimal elements of $\mathcal{G}|_Y$ containing $S$ are in $\{Y_{i_1}\cap Y,\dots,Y_{i_{r-1}}\cap Y\}$, but these are incomparable by containment, since they intersect transversally. It follows that these are precisely the minimal elements containing $S$, completing the proof of the first case.
\smallskip

\textbf{Case 2:} Let us now suppose that $Y\not\in\{Y_{i_1},\dots,Y_{i_r}\}$. Since $S\subseteq Y$, at least one of the $Y_{i_1},\dots,Y_{i_r}$ is contained in $Y$.  By relabeling if necessary, we may assume that $Y_{i_r} \subsetneq Y$. If $r=1$, then $S=Y_{i_r}$ and the desired conclusion holds, so we now consider the case $r \geq 2$.  

Given $1\leq j\leq r-1$, let us show that $Y_{i_j}$ and $Y$ meet transversally on $X$. By contradiction assume that they do not meet transversally. By Lemma~\ref{transversality.within.building.set}, since $Y_{i_j},Y\in\mathcal{G}$, then $Y_{i_j}\cap Y\in\mathcal{G}$. Given that $S\subseteq Y_{i_j}\cap Y\subseteq Y_{i_j}$, the minimality of $Y_{i_j}$ among the elements of $\mathcal{G}$ containing $S$ implies that $Y_{i_j}\cap Y= Y_{i_j}$, so $Y_{i_j}\subseteq Y$. 

On the other hand, we have that $Y_{i_j}$ and $Y_{i_r}$ intersect transversally on $X$, as they are distinct $\mathcal{G}$-factors of $S$. Then, given any point $p\in S$, we have that
\[
    T_{Y_{i_j},p} + T_{Y_{i_r},p} = T_{X,p}.
\]
Moreover, $T_{Y_{i_j},p}\subseteq T_{Y,p}$ and $T_{Y_{i_r},p}\subseteq T_{Y,p}$, which would imply that $T_{X,p}\subseteq T_{Y,p}\subsetneq T_{X,p}$. This is a contradiction, and we conclude that $Y_{i_j}$ and $Y$ intersect transversally for each $1\leq j\leq r-1$. 

Notice that $Y_{i_1} \cap Y,\dots,Y_{i_{r-1}} \cap Y,Y_{i_r}=Y_{i_r} \cap Y$ are proper subsets of $Y$ and hence they are in $\mathcal{G}|_Y$.  
Next we show that $Y_{i_1} \cap Y,\dots,Y_{i_{r-1}} \cap Y,Y_{i_r} \cap Y$ intersect transversally in $Y$. First, note that
\begin{equation} \label{equation.codimension.ell}
    \codim_Y(Y_{i_\ell}\cap Y)=\codim_X(Y_{i_\ell})\quad\text{for all}\quad 1\leq\ell \leq r-1,
\end{equation}
since $Y_{i_\ell}$ and $Y$ intersect transversally on $X$. 

Similarly, since $Y_{i_r}\subseteq Y$ we have that
\begin{equation}\label{equation.codimension.r}
\begin{aligned} 
    \codim_Y(Y_{i_r}\cap Y) = \codim_Y(Y_{i_r})
    =\dim Y - \dim Y_{i_r}
    &=\dim X - \dim Y_{i_r} - (\dim X - \dim Y)\\
    &=\codim_X(Y_{i_r}) - (\dim X - \dim Y).
\end{aligned}
\end{equation}

Lastly, note that since $S\subseteq Y$, we have
\begin{equation}\label{equation.codimension.S}
\begin{aligned} 
    \codim_Y(S) = \dim Y - \dim S 
    &=\dim X - \dim S - (\dim X - \dim Y) \\
    &= \codim_X(S) -  (\dim X - \dim Y).
\end{aligned}
\end{equation}

Putting together \eqref{equation.codimension.ell}, \eqref{equation.codimension.r} and \eqref{equation.codimension.S}, and using \eqref{equation.codimension.sum}  we obtain 
\begin{equation}\label{equation.codim.S}
    \codim_Y(S) - \sum_{\ell=1}^r\codim_Y(Y_{i_\ell}\cap Y)  
    = \codim_X(S) - \sum_{\ell=1}^r\codim_X(Y_{i_\ell}) = 0. 
\end{equation}

The subvarieties $Y_{i_1}\cap Y,\dots,Y_{i_r}\cap Y$ belong to the arrangement $\mathcal{A}_Y$ on $Y$, and their intersection is $S$. Thus, Equation~\eqref{equation.codim.S} and Lemma~\ref{lemma.codimension.characterization.transversal.intersections} imply that they intersect transversally on $Y$. In particular, they are distinct and pairwise incomparable. Moreover, every member $Y'\cap Y$ of $\mathcal{G}|_Y$ containing $S$, with $Y'\in\mathcal{G}$, contains one of these subvarieties. Indeed, $Y'\supseteq S$, so $Y'$ contains some $\mathcal{G}$-factor $Y_{i_\ell}$ of $S$, and hence $Y'\cap Y\supseteq Y_{i_\ell}\cap Y$. Together with pairwise incomparability, this shows that $Y_{i_1}\cap Y,\dots,Y_{i_r}\cap Y$ are precisely the minimal elements of $\mathcal{G}|_Y$ containing $S$. 
This completes the proof of Case 2, thus completing the argument. 
\end{proof}

\begin{remark}[$\mathcal{G}|_Y$-factors from $\mathcal{G}$-factors]\label{remark.restriction.of.g.factors}
    In the setting of the previous theorem, the arrangement of the  building set $\mathcal{G}|_Y$ on $Y$ is $\mathcal{A}|_Y := \{S \in\mathcal{A}\,\vert\,S \subsetneq Y\}$. 
    The $\mathcal{G}|_Y$-factors of any $S \in \mathcal{A}|_Y \subseteq \mathcal{A}$ are its $\mathcal{G}$-factors different from $Y$, intersected with $Y$.  
    More precisely, we showed within the proof of Theorem~\ref{theorem.restriction.of.building.set.is.building.set} that if $Y_{i_1},\ldots,Y_{i_r}$ are the distinct $\mathcal{G}$-factors of $S$ that are not equal to $Y$ (i.e. remove $Y$ from the list if $Y$ is a $\mathcal{G}$-factor of $S$), then  $Y_{i_1} \cap Y,\ldots,Y_{i_r} \cap Y$ are precisely the distinct $\mathcal{G}|_Y$-factors of $S$. 
    Moreover, within the proof we also showed that for $S \in \mathcal{A}|_Y$ there is exactly one $\mathcal{G}$-factor of $S$ contained in $Y$ (possibly equal to $Y$ itself).  
\end{remark}

%%%%%%%%%%%%%%%%%%%%%%%%%%%%%%%%%%%%%%%%%%%%%%%%%%%%%%%%%%%%%%%%

The following proposition describes the poset structure of a restricted building set.

\begin{proposition}[Poset structure of $\mathcal{G}|_Y$]
     \label{proposition.structure.restricted.building.set}
    Let $X$ be a smooth variety, $\mathcal{G}$ be a building set, and $Y$ be an element of $\mathcal{G}$. Let $\mathcal{G}_Y$ be the subset of $\mathcal{G}$ defined as $\mathcal{G}_Y := \{ Z \in \mathcal{G}  \ | \  \emptyset \neq Y \cap Z  \notin \mathcal{G}  \ \textnormal{ or } \ Z  \subsetneq Y\}$ \ (notice that the two conditions $Y \cap Z  \notin \mathcal{G}$ and $Z  \subsetneq Y$ cannot occur simultaneously for $Z \in \mathcal{G}$). Then, 
\begin{enumerate}[label=(\roman*),ref=(\roman*)]
    \item Set theoretic description: The function $\Psi_Y: \mathcal{G}_Y \rightarrow \mathcal{G}|_Y$ defined by $\Psi_{Y}(Z) := Y \cap Z$ is a bijection.  \label{proposition.structure.restricted.building.set.i}  
    \item Poset structure: For any $Z_1,Z_2 \in \mathcal{G}_Y$, we have $Z_1 \subseteq Z_2$ if and only if 
    $\Psi_{Y}(Z_1) \subseteq \Psi_{Y}(Z_2)$.  \label{proposition.structure.restricted.building.set.ii}
    \item $\mathcal{G}$-nested set description: For any $Z \in \mathcal{G}$, we have $Z \in \mathcal{G}_Y$ if and only if $Y \not\subseteq Z$ and $\{Y, Z\}$ is a $\mathcal{G}$-nested set.  \label{proposition.structure.restricted.building.set.iii}
\end{enumerate}

\end{proposition}

\begin{proof}
\begin{enumerate}[label=(\roman*),ref=(\roman*)]
    \item Given $W \in \mathcal{G}|_Y$, we know $W = Y \cap Z$ for some $Z \in \mathcal{G}$ such that $\emptyset \neq Y \cap Z \subsetneq Y$. 
    If $Y \cap Z \notin \mathcal{G}$, then $Z \in \mathcal{G}_Y$ and $\Psi_Y(Z)=Y \cap Z=W$.  
    If instead $Y \cap Z \in \mathcal{G}$, then $W = Y \cap Z \in \mathcal{G}_Y$ and $\Psi_Y(W)=W$.
    Hence, $\Psi_Y$ is surjective. 
    
    Given $Z_1,Z_2 \in \mathcal{G}_Y$ such that $\Psi_Y(Z_1)=\Psi_Y(Z_2)$, we have $Y \cap Z_1 = Y \cap Z_2 \neq \emptyset$.
    If $Y \cap Z_1  = Y \cap Z_2 \notin \mathcal{G}$, then $Z_1 =Z_2$ by Corollary~\ref{injectivity.of.restrictions.within.building.set}.
    If instead $Y \cap Z_1  = Y \cap Z_2 \in \mathcal{G}$, from the definition of $\mathcal{G}_Y$ we deduce $Z_1 \subsetneq Y$ and $Z_2 \subsetneq Y$, hence 
    $Z_1 = Y \cap Z_1 = Y \cap Z_2 = Z_2$, and then $Z_1 = Z_2$. 
    Hence, $\Psi_Y$ is also injective, and therefore $\Psi_Y$ is a bijection.
 
    \item If $Z_1,Z_2 \in \mathcal{G}_Y$ satisfy $Z_1 \subseteq Z_2$, then clearly 
    $\Psi_{Y}(Z_1) = Y \cap Z_1 \subseteq Y \cap Z_2 = \Psi_{Y}(Z_2)$.
    Conversely, assume that $Z_1,Z_2 \in \mathcal{G}_Y$ satisfy $\Psi_{Y}(Z_1) \subseteq \Psi_{Y}(Z_2)$. 
    Since $Z_1 \in \mathcal{G}_Y$, we have that either $\emptyset \neq Y \cap Z_1  \notin \mathcal{G}$ or $ Z_1  \subsetneq Y$, and we show that $Z_1 \subseteq Z_2$ considering these two cases separately.  
    
    In the first case $\emptyset \neq Y \cap Z_1  \notin \mathcal{G}$, and then by Lemma~\ref{transversality.within.building.set} we have that $Y$ and $Z_1$ are the $\mathcal{G}$-factors of $\Psi_{Y}(Z_1) = Y \cap Z_1$.
    Notice that 
    \[
    \Psi_{Y}(Z_1) = Y \cap Z_1 \subseteq Y \cap Z_2 \subseteq Z_2.
    \] 
    Since $Z_2 \in \mathcal{G}$, then $Z_2$ contains at least one of the $\mathcal{G}$-factors of $\Psi_{Y}(Z_1)$. 
    Since $Z_2$ cannot contain $Y$, we deduce that $Z_1 \subseteq Z_2$, as desired. 

    In the remaining case $Z_1 \subsetneq Y$, and then 
    \[
    Z_1 = Y \cap Z_1 = \Psi_{Y}(Z_1) \subseteq \Psi_{Y}(Z_2) = Y \cap Z_2 \subseteq Z_2,
    \]
    and hence $Z_1 \subseteq Z_2$, as desired. 
    
    \item Let $Z \in \mathcal{G}_Y$. In particular, $Y \not\subseteq Z$. We show that $\{Y, Z\}$ is $\mathcal{G}$-nested considering separately the cases $\emptyset \neq Y \cap Z  \notin \mathcal{G}$ and $Z  \subsetneq Y$.
    
    In the first case, as before, $Y$ and $Z$ are the $\mathcal{G}$-factors of $Y \cap Z$, and then $\{Y, Z\}$ is $\mathcal{G}$-nested because it is precisely the list of $\mathcal{G}$-factors of the one-step flag $Y \cap Z$ in the arrangement of $\mathcal{G}$. 

    In the second case $Z  \subsetneq Y$, and then $\{Y, Z\}$ is $\mathcal{G}$-nested because it is precisely the union of the lists of $\mathcal{G}$-factors of the elements of the flag $Z\subseteq Y$ in the arrangement of $\mathcal{G}$.

    Conversely, suppose that $Z \in \mathcal{G}$ is such that $Y \not\subseteq Z$ and $\{Y, Z\}$ is a $\mathcal{G}$-nested set.
    If $Z \subseteq Y$, then $Z \subsetneq Y$ and hence $Z \in \mathcal{G}_Y$. 
    Then, we may assume $Z \not\subseteq Y$. 
    Hence $Y$ and $Z$ are the minimal elements of the set $\{Y, Z\}$. 
    From the characterization of a $\mathcal{G}$-nested set in \cite[Definition 2.3(ii)]{LiLi}, we deduce that $Y \cap Z$ is an element of the arrangement of $\mathcal{G}$, in particular it is not empty, and $Y$ and $Z$ are its distinct $\mathcal{G}$-factors. 
    We deduce that $\emptyset \neq Y \cap Z \notin \mathcal{G}$ and hence $Z \in \mathcal{G}_Y$. This completes the proof.  \qedhere
\end{enumerate} 
\end{proof}

\begin{remark}
    In Proposition~\ref{proposition.nested.ambient.vs.restriction} we will give a characterization of $\mathcal{G}|_Y$-nested sets in terms of $\mathcal{G}$-nested sets, generalizing Proposition~\ref{proposition.structure.restricted.building.set}\ref{proposition.structure.restricted.building.set.iii} when $Y$ is maximal in $\mathcal{G}$. 
\end{remark}

%%%%%%%%%%%%%%%%%%%%%%%%%%%%%%%%%%%%%%%%%%%%%%%%%%%%%%%%%%%%%%%%

\subsection{Strict transforms under wonderful compactifications are wonderful compactifications} \label{subsection.strict.transform.is.wonderful}

%%%%%%%%%%%%%%%%%%%%%%%%%%%%%%%%%%%%%%%%%%%%%%%%%%%%%%%%%%%%%%%%

\begin{lemma}\label{lemma.LiLi.compatible.orders.on.building.sets}
Let $X$ be a smooth variety and let $\mathcal{G}=\{Y_1,Y_2,\ldots,Y_n\}$ be a building set on $X$.
Assume that the order $Y_1,\ldots,Y_n$ respects containments (i.e.\ $Y_i\subseteq Y_j$ implies $i\le j$).
Then, for each $1\le m\le n$:
\begin{enumerate}[label=(\roman*),ref=(\roman*)]
\item \label{item.building.sets.initial.segments} the initial segment $\mathcal{G}_m:=\{Y_1,\ldots,Y_m\}$ is a building set on $X$;
\item \label{item.building.sets.initial.segments.and.last} the set $\mathcal{G}_m\cup\{Y_n\}$ is a building set on $X$ for any $m\le n-1$.
\end{enumerate}
\end{lemma}

\begin{proof}

\noindent\textbf{Step 1: Arrangements induced by the subcollections.}
Let $\mathcal{A}$ be the arrangement induced by $\mathcal{G}$ (all nonempty intersections of members of $\mathcal{G}$).
For each $m$, let $\mathcal{A}_m$ be the set of all nonempty intersections of members of $\mathcal{G}_m$, and similarly
let $\mathcal{A}_m^{+}$ be the set of all nonempty intersections of members of $\mathcal{G}_m\cup\{Y_n\}$.
Both $\mathcal{A}_m$ and $\mathcal{A}_m^+$ are contained in $\mathcal{A}$ and are closed under nonempty scheme-theoretic intersections. Hence they are arrangements on $X$. 

\medskip
\noindent\textbf{Step 2: The building-set axiom for $\mathcal{G}_m$.}
Let $S\in\mathcal{A}_m$ be nonempty; fix a representation $S=\bigcap_{i=1}^r G_i$ with $G_i\in\mathcal{G}_m$.
Write $\mathcal{G}(S):=\{\,G\in\mathcal{G}:\; S\subseteq G\,\}$.
For each index $i$, consider the finite poset $\mathcal{G}(S)\cap\{\,H\in\mathcal{G}:\; H\subseteq G_i\,\}$ and choose a \emph{minimal} element $F_i$ in it. 
Then $S\subseteq F_i\subseteq G_i$, and minimality in the restricted poset forces $F_i$ to be \emph{globally} minimal in $\mathcal{G}(S)$; hence each $F_i$ is a $\mathcal{G}$-factor of $S$. Because the order on $\mathcal{G}$ respects containments, $F_i$ appears no later than $G_i$, so $F_i\in\mathcal{G}_m$.
Let $\mathcal{F}:=\{F_i\}_{i=1}^r$ with repetitions removed. Then $S=\bigcap_{F\in\mathcal{F}}F$
and the elements of $\mathcal{F}$ are $\mathcal{G}$-factors of $S$, hence they intersect transversally by Lemma~\ref{lemma.properties.transversal.intersections}\ref{lemma.properties.transversal.intersections.item.one.subset}.  
By Lemma~\ref{lemma.codimension.characterization.transversal.intersections},  
\[
\operatorname{codim}(S)=\sum_{F\in\mathcal{F}}\operatorname{codim}(F).
\]
Since the set of all $\mathcal{G}$-factors of $S$ is transverse and the sum of their codimensions is equal to $\operatorname{codim}(S)$ by Lemma~\ref{lemma.codimension.characterization.transversal.intersections}, it follows that $\mathcal{F}$ exhausts all $\mathcal{G}$-factors of $S$.
Consequently the minimal elements of $\{H\in\mathcal{G}_m:\; S\subseteq H\}$ are precisely the members of $\mathcal{F}$, since $\mathcal{F}$ consists of all the $\mathcal{G}$-factors of $S$ and is contained in $\mathcal{G}_m$. These minimal elements meet transversally, and their intersection is $S$. This is the building-set axiom for $\mathcal{G}_m$.

\medskip
\noindent\textbf{Step 3: The building-set axiom for $\mathcal{G}_m\cup\{Y_n\}$.}
Let $S\in\mathcal{A}_m^+$. The case $S=Y_n$ is immediate, so assume $S\neq Y_n$. 

\emph{First, suppose that $S$ is an intersection of elements of $\mathcal{G}_m$ only.}
By Step~2, the $\mathcal{G}$-factors of $S$, say $E_1,\ldots,E_t$, all lie in $\mathcal{G}_m$ and are exactly the minimal elements of $\{H\in\mathcal{G}_m:\, S\subseteq H\}$. Since the minimal elements of $\{H\in\mathcal{G}:\, S\subseteq H\}$ are $E_1,\ldots,E_t$, then every $H\in\mathcal{G}$ with $S\subseteq H$ contains some $E_j$. In particular, if $Y_n\supseteq S$ then $E_j\subseteq Y_n$ for some $j$, hence $Y_n$ is not minimal above $S$ in $\mathcal{G}_m\cup\{Y_n\}$. Therefore the minimal elements of $\{H\in\mathcal{G}_m\cup\{Y_n\}:\, S\subseteq H\}$ are precisely $E_1,\ldots,E_t$, and the building-set axiom follows from Step~2.

\emph{Now suppose $S$ is not an intersection of elements of $\mathcal{G}_m$ alone.}
Write
\[
S=\Big(\bigcap_{i=1}^{r} G_i\Big)\cap Y_n
\qquad\text{with }G_i\in\mathcal{G}_m.
\]
As in Step~2, choose for each $i$ a minimal element $F_i$ of $\mathcal{G}(S)\cap\{H\subseteq G_i\}$; then each $F_i$ is a $\mathcal{G}$-factor of $S$, $F_i\in\mathcal{G}_m$, and
\[
S=\Big(\bigcap_{i=1}^{r} F_i\Big)\cap Y_n.
\]
Set $F:=\bigcap_{i=1}^{r}F_i$.
Since $S$ is not an intersection of elements of $\mathcal{G}_m$, then $F\supsetneq S$. Choose a minimal element $Y_k$ of $\mathcal{G}(S)\cap\{H\subseteq Y_n\}$. Then $Y_k$ is a $\mathcal{G}$-factor of $S$ with $S\subseteq Y_k\subseteq Y_n$.
Notice that $Y_k$ is not equal to any of the $F_i$, since $Y_k\subseteq Y_n$ but no $F_i$ can be contained in $Y_n$.
Since the $F_i$ and $Y_k$ are $\mathcal{G}$-factors of $S$, after removing repetitions among the $F_i$, Lemma~\ref{lemma.properties.transversal.intersections}\ref{lemma.properties.transversal.intersections.item.disjoint.subsets} shows that $F$ and $Y_k$ intersect transversally. By Lemma~\ref{lemma.codimension.transversality.two.subvarieties}, we deduce  
\[
\operatorname{codim}(S)=\operatorname{codim}(F)+\operatorname{codim}(Y_k). 
\]
Now, by Lemma~\ref{subadditivity.codimension.smooth.ambient} applied to $F$ and $Y_n$,
\[
\operatorname{codim}(S)\le \operatorname{codim}(F)+\operatorname{codim}(Y_n).
\]
Hence $\operatorname{codim}(Y_k)\le \operatorname{codim}(Y_n)$; since $Y_k\subseteq Y_n$, we conclude $Y_k=Y_n$. Thus $Y_n$ itself is a $\mathcal{G}$-factor of $S$. 
Finally, since we have shown that $Y_n$ is a $\mathcal{G}$-factor of $S$, set
\(\mathcal{F}':=\{F_1,\ldots,F_{r},Y_n\}\) \emph{with repetitions removed.}
Each member of $\mathcal{F}'$ is a $\mathcal{G}$-factor of $S$; hence, by Lemma~\ref{lemma.properties.transversal.intersections}\ref{lemma.properties.transversal.intersections.item.one.subset}, the members of $\mathcal{F}'$ meet transversally, and since
\[
S=\Big(\bigcap_{i=1}^{r}F_i\Big)\cap Y_n=\bigcap_{H\in\mathcal{F}'} H,
\]
their intersection is $S$. By Lemma~\ref{lemma.codimension.characterization.transversal.intersections} we have
\(\operatorname{codim}(S)=\sum_{H\in\mathcal{F}'}\operatorname{codim}(H)\).
If there were any additional $\mathcal{G}$-factor of $S$ not in $\mathcal{F}'$, the codimension sum would strictly exceed $\operatorname{codim}(S)$, a contradiction. Thus $\mathcal{F}'$ exhausts all $\mathcal{G}$-factors of $S$. 
Consequently, the minimal elements of $\{H\in\mathcal{G}_m\cup\{Y_n\}:\, S\subseteq H\}$
are precisely the members of $\mathcal{F}'$, since $\mathcal{F}'$ consists of all the $\mathcal{G}$-factors of $S$ and is contained in $\mathcal{G}_m\cup\{Y_n\}$. These minimal elements meet transversally, and their intersection is $S$. This is the building-set axiom for $\mathcal{G}_m\cup\{Y_n\}$.
\end{proof}

%%%%%%%%%%%%%%%%%%%%%%%%%%%%%%%%%%%%%%%%%%%%%%%%%%%%%%%%%%%%%%%%

\begin{theorem}[Strict transforms under wonderful compactifications are wonderful compactifications]  \label{strict.transforms.are.wonderful} 
    Let $X$ be a smooth variety, $\mathcal{G}$ be a building set on $X$, and $Y$ be a maximal element of $\mathcal{G}$. 
    Then, the strict transform $\widetilde{Y}$ of $Y$ in the wonderful compactification of the building set $\mathcal{G}':=\mathcal{G} \smallsetminus \{Y\}$ on $X$ is canonically isomorphic to the wonderful compactification of the building set $\mathcal{G}|_Y$ on $Y$. 
    That is, there exists a canonical isomorphism 
    \[
    \widetilde{Y} \cong Y_{\mathcal{G}|_Y}. 
    \] 
\end{theorem}

\begin{proof}
We may assume that $\mathcal{G}|_Y$ is nonempty, since otherwise the conclusion trivially holds. Let 
\[
X^\circ := X\smallsetminus\bigcup_{Z\in\mathcal{G}'} Z,\qquad
Y^\circ := Y\cap X^\circ \;=\; Y\smallsetminus \bigcup_{\substack{Z\in\mathcal{G}'\\ Z\cap Y\neq\emptyset}}(Z\cap Y).
\]
In particular, $X^\circ\subseteq X$ and $Y^\circ\subseteq Y$ are open immersions. 

\smallskip
\textbf{Step 1: Ambient products.}
Consider the following finite products:
\begin{gather*}
A := \prod_{Z\in\mathcal{G}'} \Bl_Z X \;=\;
\Big(\prod_{\substack{Z\in\mathcal{G}'\\ Z\cap Y\neq\emptyset}} \Bl_Z X\Big)\times
\Big(\prod_{\substack{Z\in\mathcal{G}'\\ Z\cap Y=\emptyset}} \Bl_Z X\Big),\\
B := \Big(\prod_{\substack{Z\in\mathcal{G}'\\ Z\cap Y\neq\emptyset}} \Bl_{Z\cap Y} Y\Big)\times
\Big(\prod_{\substack{Z\in\mathcal{G}'\\ Z\cap Y=\emptyset}} Y\Big), 
\qquad C := \prod_{\substack{Z\in\mathcal{G}'\\ Z\cap Y\neq\emptyset}} \Bl_{Z\cap Y} Y, 
\qquad D := \prod_{W\in \mathcal{G}|_Y}\Bl_W Y,
\end{gather*}
where in $D$ we remove repetitions among the various intersections $Z\cap Y$ appearing in $C$ and index by the distinct members $W\in \mathcal{G}|_Y$.
The products indexed by $\{Z \in \mathcal{G}' | Z \cap Y = \emptyset\}$ appearing as the second factors in $A$ and $B$ are omitted if that index set is empty. 

Each factor carries a canonical embedding of $Y^\circ$, hence $Y^\circ$ embeds diagonally into $A$, $B$, $C$ and $D$.

\smallskip
\textbf{Step 2: Closure description and the three comparison maps.}
By the closure definition of the wonderful compactification, the following hold:
\begin{itemize}
\item $X_{\mathcal{G}'}$ is the closure of $X^\circ$ in $A$; the strict transform $\widetilde Y\subseteq X_{\mathcal{G}'}$ is the closure of $Y^\circ$ in $A$.
\item $Y_{\mathcal{G}|_Y}$ is the closure of $Y^\circ$ in $D$.
\end{itemize}
Thus it suffices to show that the closures of $Y^\circ$ yield canonically isomorphic varieties when taken in $A$, $B$, $C$ and $D$.
We produce morphisms
\[
B\overset{i}{\hookrightarrow} A,\qquad B \;\;\overset{p}{\underset{j}{\rightleftarrows}}\;\; C,\qquad C \;\;\overset{q}{\underset{s}{\rightleftarrows}}\;\; D
\]
compatible with the $Y^\circ$-embeddings and check that they induce canonical isomorphisms on the closures of $Y^\circ$.

\smallskip
\emph{(2.1) From \(B\) to $A$.}
Given $Z\in\mathcal{G}'$ with $Z\cap Y\neq\emptyset$,
by Lemma~\ref{lemma.strict.transforms.as.blowups} 
the closed immersion $Y\hookrightarrow X$ induces a canonical closed immersion $\Bl_{Z\cap Y}Y\hookrightarrow \Bl_Z X$. 
Given $Z\in\mathcal{G}'$ with $Z\cap Y=\emptyset$, then $Y\subseteq X\smallsetminus Z$ and the open immersion $X\smallsetminus Z\hookrightarrow \Bl_Z X$ yields a canonical closed immersion $Y\hookrightarrow \Bl_Z X$.
Taking products over all $Z\in\mathcal{G}'$ gives a closed immersion $i: B\hookrightarrow A$ compatible with the diagonal embeddings of $Y^\circ$ into $A$ and $B$.
Therefore the closure of $Y^\circ$ in $A$ is the image (under $i$) of the closure of $Y^\circ$ in $B$, and the two closures coincide via~$i$.

\smallskip
\emph{(2.2) Between $B$ and $C$.}
Write $B=C\times Y^r$, where $r$ is the cardinality of $\{Z\in\mathcal{G}' : Z\cap Y=\emptyset\}$, and let $p: B \to C$ be the projection.
If $r=0$, then $B=C$ and there is nothing to prove.
Suppose $r\geq 1$. Since $\mathcal{G}|_Y$ is nonempty, we may choose $Z_0\in\mathcal{G}'$ with $Z_0\cap Y\neq\emptyset$.
Let $f: C \to Y$ be the composition
\[
C \xrightarrow{\operatorname{pr}_{Z_0}}
\Bl_{Z_0\cap Y}Y \xrightarrow{\pi_{Z_0}} Y,
\]
where $\operatorname{pr}_{Z_0}$ is the projection and $\pi_{Z_0}$ is the blowdown map, and define
\[
j=(\operatorname{id}_C,\Delta_r\circ f) : C \longrightarrow C\times Y^r=B,
\]
where $\Delta_r: Y \to Y^r$ is the diagonal map. By construction, $p\circ j=\operatorname{id}_C$, and both $p$ and $j$ are compatible with the embeddings of $Y^\circ$. Thus they restrict to morphisms between the closures of $Y^\circ$ in $B$ and $C$.
On the closure in $B$, the composition $j\circ p$ agrees with the identity on the dense subset $Y^\circ$, hence everywhere.
The restricted maps are therefore inverse isomorphisms.
Since the restricted projection $p$ is canonical, so is this isomorphism, independently of the choice of $Z_0$.

\smallskip
\emph{(2.3) Between $C$ and $D$.}
If distinct $Z,Z'\in\mathcal{G}'$ satisfy $Z\cap Y=Z'\cap Y=W$, then the two factors $\Bl_{Z\cap Y}Y$ and $\Bl_{Z'\cap Y}Y$ of $C$ canonically identify with the single factor $\Bl_W Y$ of~$D$.
Define $q:C\to D$ by choosing one representative from each class of repeated factors and projecting onto those representatives, and define $s:D\to C$ by the diagonal embedding into each repeated class. 
These maps are compatible with the $Y^\circ$-embeddings; on the $Y^\circ$-locus they are mutually inverse, hence they induce inverse isomorphisms on the closures.
Since $s$ is canonical, the induced isomorphisms are independent of the representatives chosen to define $q$. 

\smallskip
\textbf{Step 3: Conclusion.}
By (2.1)–(2.3), the closures of $Y^\circ$ yield canonically isomorphic varieties when taken in any of $A$, $B$, $C$ and $D$.
By definition, the closure in $A$ is \(\widetilde Y\subseteq X_{\mathcal{G}'}\), and the closure in $D$ is \(Y_{\mathcal{G}|_Y}\).
Therefore \(\widetilde Y \cong Y_{\mathcal{G}|_Y}\) canonically.
\end{proof}

%%%%%%%%%%%%%%%%%%%%%%%%%%%%%%%%%%%%%%%%%%%%%%%%%%%%%%%%%%%%%%%%

\subsection{Restricted building sets and iterated blow ups} \label{subsection.restricted.building.set.iterated.blowup}

In the next result, Theorem~\ref{theorem.compatibility.restriction.iterated.blowups}, we show that in the setting of Theorem~\ref{strict.transforms.are.wonderful}, the  iterated blowup constructions of the wonderful compactifications $X_{\mathcal{G}'}$ and $Y_{\mathcal{G}|_Y}$ are compatible.  
We use this to establish some properties of the restriction of boundary divisors in Corollary~\ref{corollary.restriction.boundary.divisors} and to characterize the nested subsets of a restricted building set in Proposition~\ref{proposition.nested.ambient.vs.restriction}.

\begin{theorem}[Compatibility of the iterated blowup constructions]
\label{theorem.compatibility.restriction.iterated.blowups}
Let $X$ be a smooth variety and let
$\mathcal{G}=\{Y_1,\ldots,Y_n,Y_{n+1}\}$
be a building set on $X$, ordered compatibly with containments (i.e. $Y_i \subseteq Y_j$ implies $i \leq j$).
Set $Y:=Y_{n+1}$, which is maximal in $\mathcal{G}$ since it is the last element, set
$\mathcal{G}':=\mathcal{G}\smallsetminus\{Y\}$,
and let $\widetilde Y\subseteq X_{\mathcal{G}'}$ denote the strict
transform of $Y$.
The order on $\mathcal{G}|_Y$ induced by the bijection $\Psi_Y:\mathcal{G}_Y\rightarrow\mathcal{G}|_Y$, $\Psi_Y(Y_i): = Y_i\cap Y$, respects containments and therefore satisfies the hypothesis of
Theorem~\ref{thm: li blow-up building set}.
Moreover, after omitting the steps indexed by $Y_i\notin\mathcal{G}_Y$, the iterated blowup construction of $X_{\mathcal{G}'}$
along the dominant transforms of $Y_1,\ldots,Y_n$ induces on $Y$
precisely the iterated blowup construction of $Y_{\mathcal{G}|_Y}$
provided by that theorem in the induced order.
The resulting isomorphism
$\widetilde Y\cong Y_{\mathcal{G}|_Y}$
is the canonical isomorphism of
Theorem~\ref{strict.transforms.are.wonderful}.
\end{theorem}

\begin{proof}
The order on $\mathcal{G}|_Y$ respects containments by
Proposition~\ref{proposition.structure.restricted.building.set}\ref{proposition.structure.restricted.building.set.ii}.
Since $\mathcal{G}|_Y$ is a building set by Theorem~\ref{theorem.restriction.of.building.set.is.building.set}, Lemma~\ref{lemma.LiLi.compatible.orders.on.building.sets}\ref{item.building.sets.initial.segments}
shows that both orders satisfy the hypothesis of
Theorem~\ref{thm: li blow-up building set}.

For $0\leq m\leq n$, let $X^{(m)}$ be the variety obtained by
blowing up successively the dominant transforms of
$Y_1,\ldots,Y_m$, let $Y_k^{(m)}$ be the dominant transform of $Y_k$ in $X^{(m)}$ for $1\leq k\leq n+1$, and set $Y^{(m)}:=Y_{n+1}^{(m)}$, which is the strict transform of $Y$.
By Definition~\ref{bg.def.transform.of.building.set},
$\mathcal{G}^{(m)}=\{Y_1^{(m)},\ldots,Y_n^{(m)},Y^{(m)}\}$ is a building set, so every nonempty intersection of its members is a smooth irreducible variety.
For $k\leq m$, the subvariety $Y_k^{(m)}$ is a divisor, being a dominant transform of the exceptional divisor of the $k$\textsuperscript{th} blowup. 
Since $\mathcal{G}'=\{Y_1,\ldots,Y_n\}$ is a building set by Lemma~\ref{lemma.LiLi.compatible.orders.on.building.sets}\ref{item.building.sets.initial.segments}, ordered compatibly with containments, Proposition~\ref{bg.prop.transform.of.building.set} gives $X^{(n)}=X_{\mathcal{G}'}$, and hence $Y^{(n)}=\widetilde Y$.

For $1\leq m\leq n$, let $\pi_m:X^{(m)}\to X^{(m-1)}$ be the $m$\textsuperscript{th} blowup and let $\rho_m:Y^{(m)}\to Y^{(m-1)}$ be its restriction.
For $Y_k\in\mathcal{G}_Y$ set $y_k:=Y_k\cap Y$.
We prove by induction on $m$ the following claim, in which $y_k^{(0)}:=y_k$ and, for $m\geq 1$, the subvariety $y_k^{(m)}\subseteq Y^{(m)}$ is defined using (a): it is $\rho_m^{-1}\big(y_k^{(m-1)}\big)$ if $Y_m\notin\mathcal{G}_Y$, and the dominant transform of $y_k^{(m-1)}$ under the blowup $\rho_m$ along $y_m^{(m-1)}$ if $Y_m\in\mathcal{G}_Y$.  
Thus $y_k^{(m)}$ is the dominant transform of $y_k$ in the construction induced on $Y$ after $m$ steps; it is irreducible and maps onto $y_k$.

\begin{claim}
For each $0\leq m\leq n$ the following hold:
\begin{enumerate}[label=(\alph*),ref=(\alph*)]
    \item if $m\geq 1$, then $\rho_m$ is an isomorphism when $Y_m\notin\mathcal{G}_Y$, and it is the blowup of $Y^{(m-1)}$ along $y_m^{(m-1)}$ when $Y_m\in\mathcal{G}_Y$;
    \item $Y_k^{(m)}\cap Y^{(m)}=y_k^{(m)}$ for every $Y_k\in\mathcal{G}_Y$;
    \item $Y_k^{(m)}\cap Y^{(m)}=\emptyset$ for every $Y_k\notin\mathcal{G}_Y$ with $k\leq m$.
\end{enumerate}
\end{claim}

The case $m=0$ holds by definition.
Assume the claim for $m-1$ and set $C:=Y_m^{(m-1)}$, the center of $\pi_m$.
Since $Y\not\subseteq Y_m$, the subvariety $Y^{(m)}$ is the strict transform of $Y^{(m-1)}$.

Suppose first that $Y_m\in\mathcal{G}_Y$.
By (b) for $m-1$ we have $C\cap Y^{(m-1)}=y_m^{(m-1)}$, so Lemma~\ref{lemma.strict.transforms.as.blowups} shows that $\rho_m$ is the blowup of $Y^{(m-1)}$ along $y_m^{(m-1)}$, which is (a).
Since $Y_m^{(m)}=\pi_m^{-1}(C)$, we get
$Y_m^{(m)}\cap Y^{(m)}=\rho_m^{-1}\big(C\cap Y^{(m-1)}\big)=\rho_m^{-1}\big(y_m^{(m-1)}\big)=y_m^{(m)}$,
which is (b) for $k=m$.
Now let $k\neq m$ with $Y_k\in\mathcal{G}_Y$.
We show that $y_k^{(m-1)}\not\subseteq y_m^{(m-1)}$ and $Y_k^{(m-1)}\not\subseteq C$.
If $k>m$, either containment would give $Y_k\subseteq Y_m$ after taking images in $X$ (using Proposition~\ref{proposition.structure.restricted.building.set}\ref{proposition.structure.restricted.building.set.ii} for the first one), contradicting the order.
If $k<m$, then $Y_k^{(m-1)}$ is a divisor not containing $Y^{(m-1)}$, so by (b) for $m-1$ the subvariety $y_k^{(m-1)}=Y_k^{(m-1)}\cap Y^{(m-1)}$ is a divisor of $Y^{(m-1)}$.
Since $y_m^{(m-1)}$ and $C$ are proper irreducible subvarieties of $Y^{(m-1)}$ and $X^{(m-1)}$, respectively, either containment would be an equality, and taking images would give $y_k=y_m$ or $Y_k=Y_m$, contradicting the injectivity of $\Psi_Y$ or $k\neq m$.
Hence $Y_k^{(m)}$ and $y_k^{(m)}$ are the strict transforms of $Y_k^{(m-1)}$ and $y_k^{(m-1)}$.
Since $\pi_m$ is an isomorphism over $X^{(m-1)}\smallsetminus C$, both $Y_k^{(m)}\cap Y^{(m)}$ and $y_k^{(m)}$ agree over this open subset with $y_k^{(m-1)}\smallsetminus y_m^{(m-1)}$, which is nonempty.
Moreover, $Y_k^{(m)}\cap Y^{(m)}$ is a nonempty intersection of members of $\mathcal{G}^{(m)}$, hence a smooth irreducible variety, and $y_k^{(m)}$ is a reduced irreducible variety.
Both are therefore the closure of the same nonempty open subset, so they are equal. This proves (b).
Finally, if $Y_k\notin\mathcal{G}_Y$ with $k\leq m$, then $k<m$ and hence $Y_k^{(m-1)}\cap Y^{(m-1)}=\emptyset$ by (c) for $m-1$, and (c) follows since $Y_k^{(m)}\subseteq\pi_m^{-1}\big(Y_k^{(m-1)}\big)$ and $Y^{(m)}\subseteq\pi_m^{-1}\big(Y^{(m-1)}\big)$.

Suppose instead that $Y_m\notin\mathcal{G}_Y$.
We first show that $C\cap Y^{(m-1)}=\emptyset$.
This is clear if $Y_m\cap Y=\emptyset$, so we assume $Y_m\cap Y\neq\emptyset$.
The definition of $\mathcal{G}_Y$ gives $Y_m\cap Y=Y_j$ for some $Y_j\in\mathcal{G}$.
Since $Y_m\not\subseteq Y$ (otherwise $Y_m\subsetneq Y$ and $Y_m\in\mathcal{G}_Y$) and $Y$ is maximal, $Y_j$ is properly contained in both $Y_m$ and $Y$; in particular, $j<m$.
The open subset $U:=X\smallsetminus\bigcup_{l<j}Y_l$ meets $Y_j$ because the order respects containments, and the first $j-1$ blowups are isomorphisms over $U$.
Hence $Y_m^{(j-1)}\cap Y^{(j-1)}$ and $Y_j^{(j-1)}$ agree over the preimage of $U$.
Both are irreducible, the former being a nonempty intersection of members of $\mathcal{G}^{(j-1)}$, so
$Y_m^{(j-1)}\cap Y^{(j-1)}=Y_j^{(j-1)}$.
By Lemma~\ref{lemma: useful lemmas from li's}\ref{lemma: useful lemmas from li's item 2}, blowing up this intersection separates the strict transforms of $Y_m$ and $Y$, and their dominant transforms remain disjoint under the subsequent blowups. Thus $C\cap Y^{(m-1)}=\emptyset$.
Consequently, $\pi_m$ is an isomorphism over a neighborhood of $Y^{(m-1)}$, so $\rho_m$ is an isomorphism, which is (a).
Over $X^{(m-1)}\smallsetminus C$ every dominant transform coincides with the preimage, so
$Y_k^{(m)}\cap Y^{(m)}=\rho_m^{-1}\big(Y_k^{(m-1)}\cap Y^{(m-1)}\big)$ for $k\neq m$, and (b) and (c) for $k\neq m$ follow from the inductive hypothesis.
For $k=m$, the divisor $Y_m^{(m)}=\pi_m^{-1}(C)$ is disjoint from $Y^{(m)}$, which gives (c).
This completes the proof of the claim.

By (a), after omitting the steps with $Y_m\notin\mathcal{G}_Y$, the morphism $\widetilde Y=Y^{(n)}\to Y$ is the iterated blowup of $Y$ along the dominant transforms of the elements of $\mathcal{G}|_Y$ in the induced order.
By Theorem~\ref{thm: li blow-up building set}, the resulting variety is $Y_{\mathcal{G}|_Y}$.
The resulting isomorphism
$\widetilde Y\cong Y_{\mathcal{G}|_Y}$ agrees with that of
Theorem~\ref{strict.transforms.are.wonderful}, since both
restrict to the identity on the common dense open subset
$Y^\circ:=Y\smallsetminus\bigcup_{i=1}^n(Y_i\cap Y)$.
\end{proof}

%%%%%%%%%%%%%%%%%%%%%%%%%%%%%%%%%%%%%%%%%%%%%%%%%%%%%%%%

The assumption in Theorem~\ref{theorem.compatibility.restriction.iterated.blowups} that the order of $\mathcal{G}$ respects containments cannot be weakened to the condition in \cite[Theorem~1.3]{LiLi} that the initial segments of $\mathcal{G}$ are building sets, as the following example shows.

\begin{example}\label{example.order.respecting.containments.needed}
Let $X=\mathbb{A}^4$ with coordinates $a,b,c,d$, and let $Z_1=V(b,c)$, $Z_2=V(d,a+b)$, $W=V(a,b,c,d)$ and $Y=V(a)$.
Then $\mathcal{G}=\{Z_1,Z_2,W,Y\}$ is a building set on $X$, $Y$ is maximal in $\mathcal{G}$, and the initial segments of the order $Z_1,Z_2,W,Y$ are building sets.
However, this order does not respect containments, since $W\subsetneq Z_1$.
Let $X^{(1)}$ be the blowup of $X$ along $Z_1$, let $X^{(2)}$ be the blowup of $X^{(1)}$ along the strict transform $Z_2^{(1)}$ of $Z_2$, and let $Y^{(1)}$ and $Y^{(2)}$ be the strict transforms of $Y$ in $X^{(1)}$ and $X^{(2)}$, respectively.
In the chart of $X^{(1)}$ with coordinates $a,c,v,d$, where $b=cv$, we have $Y^{(1)}=V(a)$ and $Z_2^{(1)}=V(d,a+cv)$.
Hence $Z_2^{(1)}\cap Y^{(1)}$ is defined by the ideal $\langle d,cv\rangle$ in the coordinates $c,v,d$ of $Y^{(1)}$, and it is the union of two lines meeting at a point.
By Lemma~\ref{lemma.strict.transforms.as.blowups}, $Y^{(2)}$ is the blowup of $Y^{(1)}$ along this intersection. The chart of this blowup in which $d$ generates the pulled-back ideal is the hypersurface $V(dt-cv)\subseteq\mathbb{A}^4$, with coordinates $c,v,d,t$. This hypersurface has an ordinary double point at the origin.
Therefore $Y^{(2)}$ is singular, and the construction induced on $Y$ is not an iterated blowup along smooth centers.
This reflects the fact that $\{Y\cap Z_1,Y\cap Z_2\}$ is not a building set on $Y$: the two lines $Y\cap Z_1$ and $Y\cap Z_2$ meet at the origin but do not intersect transversally in $Y$.
The strict transform of $Y$ in $X_{\mathcal{G}\smallsetminus\{Y\}}$ is still isomorphic to $Y_{\mathcal{G}|_Y}$ by Theorem~\ref{strict.transforms.are.wonderful}.
\end{example}

%%%%%%%%%%%%%%%%%%%%%%%%%%%%%%%%%%%%%%%%%%%%%%%%%%%%%%%%

\begin{corollary}[Restriction of boundary divisors]
\label{corollary.restriction.boundary.divisors}
In the setting of Theorem~\ref{theorem.compatibility.restriction.iterated.blowups}, identify $Y_{\mathcal{G}|_Y}$ with $\widetilde Y$ and let $\iota:Y_{\mathcal{G}|_Y}\hookrightarrow X_{\mathcal{G}'}$ denote the inclusion.
For $1\leq i\leq n$, let $E_i$ be the boundary divisor of $X_{\mathcal{G}'}$ corresponding to $Y_i$, namely the dominant transform of $Y_i$ in $X_{\mathcal{G}'}$.
If $Y_i\in\mathcal{G}_Y$, let $e_i$ be the boundary divisor of $Y_{\mathcal{G}|_Y}$ corresponding to $Y_i\cap Y$, namely the dominant transform of $Y_i\cap Y$ under the iterated blowup $Y_{\mathcal{G}|_Y}\to Y$ of Theorem~\ref{theorem.compatibility.restriction.iterated.blowups}.
Then
\[
E_i\cap Y_{\mathcal{G}|_Y}
=
\begin{cases}
e_i, & \textnormal{if }Y_i\in\mathcal{G}_Y,\\
\emptyset, & \textnormal{if }Y_i\notin\mathcal{G}_Y,
\end{cases}
\qquad
\iota^*(E_i)
=
\begin{cases}
e_i, & \textnormal{if }Y_i\in\mathcal{G}_Y,\\
0, & \textnormal{if }Y_i\notin\mathcal{G}_Y,
\end{cases}
\]
where the second equality holds in Chow under the pullback $\iota^*:A^\bullet(X_{\mathcal{G}'}) \rightarrow A^\bullet(Y_{\mathcal{G}|_Y})$.
\end{corollary}

\begin{proof}
In the notation of the proof of Theorem~\ref{theorem.compatibility.restriction.iterated.blowups}, we have $E_i=Y_i^{(n)}$ and $Y_{\mathcal{G}|_Y}=Y^{(n)}$, and $e_i=y_i^{(n)}$ when $Y_i\in\mathcal{G}_Y$.
Hence the first assertion is the case $m=n$ of parts (b) and (c) of the claim in that proof.

Since each $E_i$ is an effective Cartier divisor and $Y_{\mathcal{G}|_Y}$ is not contained in $E_i$ (as $Y\not\subseteq Y_i$), its pullback is the class of $E_i\cap Y_{\mathcal{G}|_Y}$. Passing to Chow classes proves the second assertion.
\end{proof}

%%%%%%%%%%%%%%%%%%%%%%%%%%%%%%%%%%%%%%%%%%%%%%%%%%%%%%%%

As an application of the previous results, we give a characterization of $\mathcal{G}|_Y$-nested sets in terms of $\mathcal{G}$-nested sets when $Y$ is a maximal element of a building set $\mathcal{G}$ on a smooth variety $X$.

\begin{proposition} \label{proposition.nested.ambient.vs.restriction}
    Let $X$ be a smooth variety, $\mathcal{G}$ be a building set on $X$, and $Y$ be a maximal element of $\mathcal{G}$. 
    Let $\mathcal{G}_Y = \{ Z \in \mathcal{G}  \ | \  \emptyset \neq Y \cap Z  \notin \mathcal{G}  \ \textnormal{ or } \ Z  \subsetneq Y\}$ and $\Psi_Y: \mathcal{G}_Y \rightarrow \mathcal{G}|_Y$ be the bijection $\Psi_{Y}(Z) := Y \cap Z$ from Proposition~\ref{proposition.structure.restricted.building.set}.
    Let $\mathcal{T}$ be a nonempty subset of $\mathcal{G}$ with $Y \notin \mathcal{T}$. Then, the following are equivalent:
    \begin{enumerate}[label=(\roman*),ref=(\roman*)]
        \item The set $\mathcal{T}\cup\{Y\}$ is $\mathcal{G}$-nested.
        \item $\mathcal{T} \subseteq \mathcal{G}_Y$ and the set $\Psi_Y(\mathcal{T})$ is $\mathcal{G}|_Y$-nested.
    \end{enumerate}
\end{proposition}
\begin{proof}
Let $\mathcal{G} = \mathcal{G}_{n+1} = \{Y_1,\ldots,Y_{n+1} \}$, $\mathcal{T} = \{Y_{i_1},\ldots,Y_{i_r} \}$, and define $\mathcal{G}_{n} = \{Y_1,\ldots,Y_{n} \}  \subseteq \mathcal{G}_{n+1}$. We may assume that the order of the $Y_i$ is compatible with containments (i.e. $Y_i \subseteq Y_j$ implies $i \leq j$).  
We may also assume that $Y = Y_{n+1}$ by moving $Y$ to the end of the list and relabeling if necessary. 

By Theorem~\ref{strict.transforms.are.wonderful},
the wonderful compactification $Y_{\mathcal{G}|_Y}$ of $\mathcal{G}|_Y$ on $Y$ is the strict transform $\widetilde{Y}$ of $Y$ in the wonderful compactification $X_{\mathcal{G}_n}$ of $\mathcal{G}_n$ on $X$. 
Let $E_1,\ldots,E_n$ be the divisors corresponding to $Y_1,\ldots,Y_n$ in $X_{\mathcal{G}_n}$. 
The building set $\mathcal{G}|_Y$ consists of the subvarieties of $Y$ of the form $W \cap Y$ for $W \in \mathcal{G}_Y$. 
Moreover, by Corollary~\ref{corollary.restriction.boundary.divisors} the boundary divisors of $Y_{\mathcal{G}|_Y}$ arise as $E_i \cap Y_{\mathcal{G}|_Y}$ for each $i$ such that $Y_i \in \mathcal{G}_Y$.

\emph{(i)$\Rightarrow$(ii).}
For each $W \in \mathcal{T}$, 
we have that $\{W \}\cup\{Y\}$ is $\mathcal{G}$-nested since it is a subset of the $\mathcal{G}$-nested set $\mathcal{T}\cup\{Y\}$. 
Since $Y \not\subseteq W$ (because $Y$ is maximal in $\mathcal{G}$) and 
$\{W \}\cup\{Y\}$ is $\mathcal{G}$-nested, then $W \in \mathcal{G}_Y$ by Proposition~\ref{proposition.structure.restricted.building.set}\ref{proposition.structure.restricted.building.set.iii}. 
Then $\mathcal{T} \subseteq \mathcal{G}_Y$. 

Since $\mathcal{T} \cup \{ Y \}$ is $\mathcal{G}$-nested, then the divisors in  the wonderful compactification $X_{\mathcal{G}_{n+1}} = \operatorname{Bl}_{\widetilde{Y}} X_{\mathcal{G}_{n}}$ corresponding to each of the varieties in $\mathcal{T} \cup \{ Y \}$ have a nonempty intersection. 
Mapping those divisors to $X_{\mathcal{G}_{n}}$, we deduce that $E_{i_1},\ldots,E_{i_r}$ and $\widetilde{Y} = Y_{\mathcal{G}|_Y}$ have a nonempty intersection in $X_{\mathcal{G}_{n}}$.
Then, the divisors $E_{i_1} \cap Y_{\mathcal{G}|_Y},\ldots,E_{i_r} \cap Y_{\mathcal{G}|_Y}$ on $Y_{\mathcal{G}|_Y}$ have a nonempty intersection. 
Then, the divisors on the wonderful compactification $Y_{\mathcal{G}|_Y}$ corresponding to the elements in $\Psi_Y(\mathcal{T}) \subseteq \mathcal{G}|_Y$ have a nonempty intersection. 
Hence, $\Psi_Y(\mathcal{T})$ is $\mathcal{G}|_Y$-nested. 

    \smallskip

\emph{(ii)$\Rightarrow$(i).} 
Since $\mathcal{T}\subseteq\mathcal{G}_Y$, the varieties $Y_{i_1}\cap Y,\ldots,Y_{i_r}\cap Y$ belong to $\mathcal{G}|_Y$, and their corresponding boundary divisors in $Y_{\mathcal{G}|_Y}$ are $E_{i_1}\cap Y_{\mathcal{G}|_Y},\ldots,
E_{i_r}\cap Y_{\mathcal{G}|_Y}$. 
Since $\Psi_Y(\mathcal{T})$ is $\mathcal{G}|_Y$-nested then the intersection of the divisors $E_{i_1} \cap Y_{\mathcal{G}|_Y},\ldots,E_{i_r} \cap Y_{\mathcal{G}|_Y}$ in $Y_{\mathcal{G}|_Y}$ is nonempty. 
Then, the intersection of $E_{i_1},\ldots,E_{i_r}$ and $Y_{\mathcal{G}|_Y}$ is nonempty on $X_{\mathcal{G}_n}$, and we call this intersection $Z$.
Notice that $E_{i_1},\ldots,E_{i_r}$ and $Y_{\mathcal{G}|_Y}$ are all elements of the building set $\mathcal{G}^{(n)}_{n+1}$ in $X_{\mathcal{G}_n}$ obtained by iterated dominant transforms of the elements of $\mathcal{G}$ along the iterated blowup of $X$ along the strict transforms of $Y_1,\ldots,Y_n$. 

The elements of $\mathcal{G}_{n+1}^{(n)} \smallsetminus \{Y_{\mathcal{G}|_Y}\}$ are distinct divisors (they are precisely $E_1,\ldots,E_n$), and none of them contains $Y_{n+1}^{(n)}$ (because $Y_{\mathcal{G}|_Y}  \subseteq E_j$ on $X_{\mathcal{G}_n}$ for some $1 \leq j \leq  n$ would imply $Y_{n+1}  \subseteq Y_{j}$ on $X$, which is not possible). 
In particular, each of $E_1,\ldots,E_n$ is minimal in $\mathcal{G}_{n+1}^{(n)}$. 
Also, $Y_{\mathcal{G}|_Y}$ is minimal in $\mathcal{G}_{n+1}^{(n)}$, because if it contained $E_j$ for some $1 \leq j \leq n$, then $Y_{\mathcal{G}|_Y} = E_j$ since $E_j$ is a divisor, and hence mapping down to $X$ we would get $Y_{n+1}=Y_j$, which is a contradiction. 
By \cite[Lemma~2.5]{LiLi}, this implies that $E_{i_1},\ldots,E_{i_r},Y_{\mathcal{G}|_Y}$ are the $\mathcal{G}^{(n)}_{n+1}$-factors of $Z$. 
Then, considering the one-step flag $\{Z\}$, we deduce that $E_{i_1},\ldots,E_{i_r},Y_{\mathcal{G}|_Y}$ is $\mathcal{G}^{(n)}_{n+1}$-nested. 
By repeated application of \cite[Proposition~2.8(iii)]{LiLi}, this means that $Y_{i_1},\ldots,Y_{i_r},Y_{n+1}$ is $\mathcal{G}$-nested.
Hence, $\mathcal{T}\cup\{Y\}$ is $\mathcal{G}$-nested, completing the proof. 
\end{proof}

%%%%%%%%%%%%%%%%%%%%%%%%%%%%%%%%%%%%%%%%%%%%%%%%%%%%%%%%%%%%%%%%

%%%%%%%%%%%%%%%%%%%%%%%%%%%%%%%%%%%%%%%%%%%%%%%%%%%%%%%%%%%%%%%%
%%%%%%%%%%%%%%%%%%%%%%%%%%%%%%%%%%%%%%%%%%%%%%%%%%%%%%%%%%%%%%%%
%%%%                                                        %%%%
%%%%      4. The Chow ring of a wonderful compactification    
%%%%                                                        %%%%
%%%%%%%%%%%%%%%%%%%%%%%%%%%%%%%%%%%%%%%%%%%%%%%%%%%%%%%%%%%%%%%%
%%%%%%%%%%%%%%%%%%%%%%%%%%%%%%%%%%%%%%%%%%%%%%%%%%%%%%%%%%%%%%%%

\section{The Chow ring of a wonderful compactification} \label{section.wonderful}

In this section we prove the main result, Theorem~\ref{theorem.main.intro}, describing the Chow ring of a wonderful compactification under suitable surjectivity assumptions. The proof combines the iterated blowup construction of wonderful compactifications with Keel's formula for the Chow ring of a blowup.  
We first show in Section~\ref{subsection.Chow.surjectivity} that the surjectivity assumptions are preserved throughout the iterated blowup construction. 
We then prove Theorem~\ref{theorem.main.intro} in Section~\ref{subsection.main.proof} by strong induction on the cardinality of the building set $\mathcal G$ on the smooth variety $X$. 
The base case is Keel's formula. For the inductive step, we fix a maximal element $Y\in\mathcal{G}$ and set $\mathcal{G}'=\mathcal{G}\smallsetminus\{Y\}$. The wonderful compactification $X_{\mathcal{G}}$ is the blowup of $X_{\mathcal{G}'}$ along the strict transform $\widetilde{Y}$ of $Y$, and $\widetilde{Y}$ is itself the wonderful compactification associated to the restricted building set $\mathcal{G}|_Y$. The induction hypothesis therefore gives presentations for the Chow rings of both $X_{\mathcal{G}'}$ and $\widetilde{Y}$. The technical heart of the argument is the computation of the kernel of the restriction map
$A^\bullet(X_{\mathcal{G}'})\to A^\bullet(\widetilde{Y})$.
Once this kernel is determined, a final application of Keel's formula yields the desired presentation of $A^\bullet(X_{\mathcal{G}})$ and completes the induction.

\subsection{Chow surjectivity under blow-ups} \label{subsection.Chow.surjectivity}

Let $X$ be a smooth variety, let $\mathcal{G}$ be a building set with induced
arrangement $\mathcal{A}$, and let $F\in \mathcal{G}$ be a minimal element. Let $\pi:\widetilde X:=\operatorname{Bl}_{F}X\longrightarrow X$ 
be the blow-up of $X$ along $F$, and let $E\subseteq \widetilde X$ be the
exceptional divisor. For every $S\in\mathcal{A}$, denote by $\widetilde S$ its
dominant transform in $\widetilde X$. Define
\[
        \widetilde{\mathcal{A}}
        :=
        \{\widetilde S\}_{S\in\mathcal{A}}
        \cup
        \{\widetilde S\cap E\}_{\substack{S\in\mathcal{A} \\ \emptyset\subsetneq S\cap F\subsetneq S}}
\qquad
\textnormal{and}
\qquad
        \widetilde{\mathcal{G}}
        :=
        \{\widetilde S\}_{S\in\mathcal{G}}.
\]
By \cite[Proposition~2.8]{LiLi}, the collection $\widetilde{\mathcal{A}}$ is an arrangement of subvarieties of $\widetilde X$ and the collection
$\widetilde{\mathcal{G}}$ is a building set for $\widetilde{\mathcal{A}}$. 
The following lemma shows that the hypothesis of Chow-surjectivity from the ambient variety to the elements of the arrangement is preserved after blowing up the minimal element
$F\in\mathcal{G}$. 

\begin{lemma}
\label{lemma.chow.surjectivity.after.minimal.blowup}
With the notation introduced above, assume that, for every $S\in\mathcal{A}$,
the restriction homomorphism $A^\bullet(X)\longrightarrow A^\bullet(S)$ is surjective. Then, for every $T\in\widetilde{\mathcal{A}}$, the restriction homomorphism $A^\bullet(\widetilde X)\longrightarrow A^\bullet(T)$ is surjective.
\end{lemma}
\begin{proof}
By definition, every $T\in\widetilde{\mathcal{A}}$ has the form $\widetilde S$ for some $S\in\mathcal{A}$, or $\widetilde S\cap E$ for some $S\in\mathcal{A}$ with $\emptyset\subsetneq S\cap F\subsetneq S$. We consider these two cases.
We use Hartshorne's convention for projective bundles \cite[Section~II.7]{Hartshorne}.  

First suppose that $T=\widetilde S$ for some $S\in\mathcal{A}$, and set $C:=S\cap F$. There are three subcases: $C=\emptyset$, $C=S$ (equivalently, $S\subseteq F$), or $\emptyset\subsetneq C\subsetneq S$.

If $C=\emptyset$, then $\widetilde S\simeq S$, and the surjectivity of $A^\bullet(X)\to A^\bullet(S)$ implies that $A^\bullet(\widetilde X)\to A^\bullet(\widetilde S)$ is already surjective on $\pi^*A^\bullet(X)$.

If $S\subseteq F$, then, by the definition of dominant transform, $\widetilde S=\pi^{-1}(S)\simeq\mathbb P(N_{F/X}^\vee|_S)$. By the projective bundle formula \cite[Theorem~3.3(b)]{Fulton-book}, its Chow ring is generated by pullbacks of classes in $A^\bullet(S)$ and the class $\xi:=c_1(\mathcal O_{\mathbb P(N_{F/X}^\vee|_S)}(1))$. The former are restrictions of classes in $\pi^*A^\bullet(X)$ by hypothesis, and $\xi=-E|_{\widetilde S}$. Hence $A^\bullet(\widetilde X)\to A^\bullet(\widetilde S)$ is surjective.

Suppose now that $\emptyset\subsetneq C\subsetneq S$. Then $C\in\mathcal{A}$, and the dominant transform $\widetilde S$ is the strict transform of $S$. Lemma~\ref{lemma.strict.transforms.as.blowups} therefore identifies $\widetilde S$ with $\operatorname{Bl}_C S$, whose exceptional divisor is $E_S:=\widetilde S\cap E$. Since the surjective restriction $A^\bullet(X)\to A^\bullet(C)$ factors through $A^\bullet(S)$, the restriction $A^\bullet(S)\to A^\bullet(C)$ is also surjective. Keel's blow-up formula shows that $A^\bullet(\widetilde S)$ is generated by pullbacks of classes in $A^\bullet(S)$ and $E_S$. The former are restrictions of classes in $\pi^*A^\bullet(X)$, and $E_S=E|_{\widetilde S}$. Thus $A^\bullet(\widetilde X)\to A^\bullet(\widetilde S)$ is surjective.

Now suppose that $T=\widetilde S\cap E$ for some $S\in\mathcal{A}$ with $\emptyset\subsetneq C:=S\cap F\subsetneq S$. Then $C\in\mathcal{A}$, and $T$ is the exceptional divisor of $\operatorname{Bl}_C S\to S$. Hence $T\simeq\mathbb P(N_{C/S}^\vee)$. By the projective bundle formula, $A^\bullet(T)$ is generated by pullbacks of classes in $A^\bullet(C)$ and $\eta:=c_1(\mathcal O_{\mathbb P(N_{C/S}^\vee)}(1))$. Since $C\in\mathcal{A}$, the former are restrictions of classes in  $\pi^*A^\bullet(X)$, while $\eta=-E|_T$. Hence $A^\bullet(\widetilde X)\to A^\bullet(T)$ is surjective.

These cases exhaust the elements of $\widetilde{\mathcal{A}}$, completing the proof.
\end{proof}

%%%%%%%%%%%%%%%%%%%%%%%%%%%%%%%%%%%%%%%%%%%%%%%%%%%%%

\subsection{Proof of Theorem \ref{theorem.main.intro}} \label{subsection.main.proof}
We recall the statement for the reader's convenience. 

\begin{theorem}[Theorem~\ref{theorem.main.intro}] \label{theorem.main.body}
    Let $\mathcal{G}$ be a building set on a smooth variety $X$ and let $\mathcal{A}$ be its associated arrangement.
    Assume that the pullback $A^{\bullet}(X) \rightarrow A^{\bullet}(Z)$ is surjective for all $Z \in \mathcal{A}$. 
    Then, the Chow ring of the wonderful compactification $X_\mathcal{G}$ is  $A^\bullet(X_{\mathcal{G}}) \cong A^\bullet(X)[E_Y \, | \, Y \in \mathcal{G}]/I$, where $I$ is the ideal generated by:  
    \begin{enumerate}
    \setlength\itemsep{0.7em}
        \item \ \ $\prod_{Y\in\mathcal{T}} E_Y$ \quad for all $\mathcal{T}\subseteq\mathcal{G}$ such that $\mathcal{T}$ is \emph{not} $\mathcal{G}$-nested; 
        \item \ \ $J_Z \cdot E_Z$  \quad for all $Z$ in $\mathcal{A}$;
        \item \ \ $P_{X,Y}\, \raisebox{-1.5mm}{$\Biggl($}
        -\mathlarger{\sum}\nolimits_{\substack{Y'\in\mathcal{G}\\Y'\subseteq Y}} \, E_{Y'}
        \raisebox{-1.5mm}{$\Biggr)$}$
        \quad for all $Y\in\mathcal{G}$; 
        \item \ \ $E_Z \cdot P_{X,Z,Y}\, \raisebox{-1.5mm}{$\Biggl($}
        -\mathlarger{\sum}\nolimits_{\substack{Y'\in\mathcal{G}\\Y'\subseteq Y}} \, E_{Y'}
        \raisebox{-1.5mm}{$\Biggr)$}$
        \quad
        for all $Z\in\mathcal{A}$ and $Y\in\mathcal{G}$
        such that $Y\subsetneq Z$,
    \end{enumerate}
    where $J_Z \subseteq A^{\bullet}(X)$ is the kernel of the pullback map $A^{\bullet}(X) \rightarrow A^{\bullet}(Z)$, $E_Z$ denotes the product of the $E_Y$ as $Y$ ranges over the $\mathcal{G}$-factors of $Z$,  
    and $P_{X,Y}(t),P_{X,Z,Y}(t) \in A^{\bullet}(X)[t]$ denote Chern polynomials. 
\end{theorem}

In the statement we have fixed specific choices of Chern polynomials $P_{X,Y}(t)\in A^\bullet(X)[t]$ for all $Y\in\mathcal{G}$, and $P_{X,Z,Y}(t)\in A^\bullet(X)[t]$ for all $Z\in\mathcal{A}$ and $Y\in\mathcal{G}$ such that $Y\subsetneq Z$. In particular, the statement holds for any such choice.

\noindent\emph{Proof of Theorem~\ref{theorem.main.intro}} 
The statement holds trivially for the empty building set. We therefore assume that the building set $\mathcal{G}$ is nonempty and proceed by strong induction on its cardinality $n$.
 
\smallskip
\textbf{Step 1: Base case.}
    For $n=1$, we have that $\mathcal{A}=\mathcal{G}= \{ Y_1 \}$ and the relations (1) and (4) are vacuous. 
    The wonderful compactification $X_{\mathcal{G}}$ is the blowup of $X$ along $Y_1$, and the claim holds since the presentation in the statement becomes exactly Keel's formula in Lemma~\ref{lemma: keel's blowup formula} for the Chow ring of a blow up,
    \[
    A^{\bullet}(\operatorname{Bl}_{Y_1}X)=\frac{A^{\bullet}(X)[E_1]}{\langle \, J_{Y_1} \cdot E_1, \, P_{X,Y_1}(-E_1) \, \rangle} ,
    \]
    under the assumption that $A^{\bullet}(X) \rightarrow A^{\bullet}(Y_1)$ is surjective, like in our case. 

\smallskip
\textbf{Step 2: Inductive step.}
We assume that the conclusion holds for all wonderful compactifications of building sets with cardinality at most $n$ on any smooth variety.  
Let us fix a smooth variety $X$ and a building set $\mathcal{G}$  
on $X$ with cardinality $n+1$, and denote its associated arrangement by $\mathcal{A}$.
We must show that the presentation in the statement holds for $A^{\bullet}(X_{\mathcal{G}})$. 
Let us set some notation and recall some definitions. 

\begin{itemize}[leftmargin=2em]
    \item Let $Y_1,\ldots,Y_n,Y_{n+1}$ denote the elements of $\mathcal{G}$.
    We assume that $\mathcal{G}$ is ordered respecting containments, that is, if $Y_i\subseteq Y_j$, then $i\leq j$ (for example, any order such that $\operatorname{dim}Y_i$ is weakly increasing works). Denote $Y:=Y_{n+1}$.    
    \item Let $\mathcal{G}|_Y$ denote the building set obtained by restriction of $\mathcal{G}$ to $Y$ and let $\mathcal{A}_{Y}$ denote its arrangement.         
    \item Let $y_i:= Y_i \cap Y$ for each $Y_i \in \mathcal{G}_Y$, where 
    $\mathcal{G}_Y := \{ Z \in \mathcal{G}  \ | \  \emptyset \neq Y \cap Z  \notin \mathcal{G}  \ \textnormal{ or } \ Z  \subsetneq Y\}$. Recall that, $Y_i \mapsto y_i$ provides a bijection between $\mathcal{G}_Y$ and $\mathcal{G}|_Y$ by Proposition~\ref{proposition.structure.restricted.building.set}. 
    \item For each $1 \leq m \leq n+1$, let $A^{\bullet}(X)[E_1,\ldots,E_m]$ denote the polynomial ring in the variables $E_i$, for $1 \leq i \leq m$. We consider these rings as subrings of $A^{\bullet}(X)[E_1,\ldots,E_{n+1}]$ in the natural way.   
    \item For each $1 \leq m \leq n+1$, let $\mathcal{G}_m:= \{Y_1,\ldots,Y_m\}$, which is a building set on $X$ by Lemma~\ref{lemma.LiLi.compatible.orders.on.building.sets}\ref{item.building.sets.initial.segments}. We denote by $\mathcal{A}_m$ the arrangement induced by $\mathcal{G}_m$. Notice that $\mathcal{A}_m \subseteq \mathcal{A}$ and that by our assumption on the order of $\mathcal{G}$, $\operatorname{Fact}(\mathcal{G}_{m},Z)=\operatorname{Fact}(\mathcal{G},Z)$ for all $Z \in \mathcal{A}_m$ (indeed, each $W \in \operatorname{Fact}(\mathcal{G}_{m},Z)$ is minimal over $Z$ in $\mathcal{G}_{m}$, and hence it is also minimal over $Z$ in $\mathcal{G}$, by the chosen order; hence $W \in \operatorname{Fact}(\mathcal{G},Z)$, which shows $\operatorname{Fact}(\mathcal{G}_{m},Z) \subseteq \operatorname{Fact}(\mathcal{G},Z)$; and there cannot be additional $\mathcal{G}$-factors of $Z$ by the additivity of codimensions on a transversal intersection).
    \item Given any building set 
    $\mathcal{H}=\{W_1,\ldots,W_s\}$, we will say that $I\subseteq\{1,\ldots,s\}$ is $\mathcal{H}$-nested if and only if $\{W_i\mid i\in I\}$ is $\mathcal{H}$-nested. \item Define $X^{(0)}:=X$ and $Y_i^{(0)}:=Y_i$ for each $1 \leq i \leq n+1$. Inductively define for each $1 \leq m \leq n+1$, $X^{(m)}$ to be the blow up of $X^{(m-1)}$ along $Y_{m}^{(m-1)}$, and $Y_i^{(m)} \subseteq X^{(m)}$ to be the dominant transform of $Y_i^{(m-1)}$, for each $1\leq i \leq n+1$ (here, the dominant transforms $Y_i^{(m)}$ are strict transforms when $i \neq m$ and the exceptional divisor when $i = m$). 
    By \cite[Proposition~2.8]{LiLi} (reviewed in Definition~\ref{bg.def.transform.of.building.set}), the set $\mathcal{G}^{(m)} :=\left\{Y_1^{(m)},\ldots,Y_{n+1}^{(m)}\right\}$ is a building set on $X^{(m)}$, for each $0 \leq m \leq n+1$. 
\end{itemize}

\smallskip

%%%%%%%%%%%%%%%%%%%%%%%%%%%%%%%%%%%%%%%%%%%%%%%%%%%%%

\begin{definition}\label{def:relations}
For each $1\leq m\leq n+1$, define the following subsets $R_1^{(m)}$, $R_2^{(m)}$, $R_3^{(m)}$, $R_4^{(m)}$ of 
$A^\bullet(X)[E_1, \cdots, E_m]$, which we also consider as subsets of $A^\bullet(X)[E_1, \cdots, E_{n+1}]$: 
\begin{align}
R_1^{(m)} 
&= 
\left\{ \; \prod_{i \in I} E_i \, 
\bigg| \, 
\text{$I \subseteq \{1,\ldots, m\}$ such that $I$ is not $\mathcal{G}_m$-nested
}
\right\}  
\label{eq:r_1}
\\
R^{(m)}_2 
&= 
\bigcup_{Z \in \mathcal{A}_m}
J_Z \cdot E_Z \\
R_3^{(m)}  &= 
\left\{ P_{X,Y_i} \raisebox{-1.5mm}{$\Biggl($} -\sum_{Y_j \subseteq Y_i} E_j \raisebox{-1.5mm}{$\Biggr)$} 
\; \bigg| \;
i \in \{1,\ldots, m\}
\right\}
\label{eq:r_3}
\\
R_4^{(m)} &= 
\left\{ E_Z \cdot P_{X,Z,Y_i} \raisebox{-1.5mm}{$\Biggl($} -\sum_{Y_j \subseteq Y_i } E_j \raisebox{-1.5mm}{$\Biggr)$}
\, \bigg| \,
\text{ for each $Z \in \mathcal{A}_m$ and $Y_i \subsetneq Z$}
\right\} 
\end{align}
\end{definition}

Since the order of each $\mathcal{G}_m$ is compatible with containments, by Lemma~\ref{lemma.LiLi.compatible.orders.on.building.sets}\ref{item.building.sets.initial.segments} 
and Li's Theorem~\ref{thm: li blow-up building set}, 
$X_{\mathcal{G}_m}$ is the iterated blowup of $X$ along the dominant transforms of the elements of $\mathcal{G}_m$ in the given order.  
Then, the variety $X^{(m)}$ is the wonderful compactification of $\mathcal{G}_m$ on $X$. 
Moreover, for each $Z \in \mathcal{A}_m$ we have $Z \in \mathcal{A}_{n+1} = \mathcal{A}$, and then the pullback homomorphism $A^\bullet(X) \rightarrow A^\bullet(Z)$ is surjective by assumption. 
Then, for each $1 \leq m \leq n$, by the induction hypothesis we have 
\[
A^\bullet\left(X^{(m)}\right)
=
\frac{A^\bullet(X)[E_1,\ldots,E_m]}{
\left\langle
R_1^{(m)}, R_2^{(m)}, R_3^{(m)}, R_4^{(m)}
\right\rangle
},
\]
and the desired conclusion is the case $m=n+1$ of this formula. 

\begin{definition}\label{defn:preimage_Rk}
Define the following subsets $r_1$, $r_2$ and $r_4$ of 
$A^\bullet(X)[E_1, \cdots, E_{n}]$:
\begin{align*}
r_1 &=  \left\{
\prod_{i\in I}E_i \ \middle| \ I\subseteq\{1,\dots,n\} \text{ such that } I\cup\{n+1\} \text{ is not $\mathcal{G}_{n+1}$-nested}
\right\}, 
\\
r_2 &= 
J_{Y_{n+1}}
\; \cup \; 
\; 
\bigcup_{\substack{
  Z \in \mathcal{A}_{n+1} \\
  Z \subsetneq Y_{n+1}
}}
\;
\raisebox{-2.0mm}{$\Biggl($}
 \ J_Z  \ \ \cdot \! \! \! \! \! \! \prod_{\substack{
  1 \leq i \leq n \\
  Y_i \in \operatorname{Fact}(\mathcal{G}_{n+1},Z)
}}
  \! \! \! \! \! \!   E_i  \
\raisebox{-2.0mm}{$\Biggr)$}, 
\\ 
r_4 &=  
\left\{
\frac{E_Z}{E_{n+1}} \cdot P_{X,Z,Y_i} \raisebox{-1.5mm}{$\Biggl($} -\sum_{Y_j\subseteq Y_i} E_j \raisebox{-1.5mm}{$\Biggr)$} \; \bigg| \;
Z\in\mathcal{A}_{n+1}, Y_{n+1} \in \operatorname{Fact}(\mathcal{G}_{n+1},Z), Y_i\subsetneq Z 
\right\}.
\end{align*}
\end{definition}

%%%%%%%%%%%%%%%%%%%%%%%%%%%%%%%%%%%%%%%%%%%%%%%%%%%%%

\begin{lemma}\label{lemma:equalities_R_k}
With the notation in Definition~\ref{def:relations} and Definition~\ref{defn:preimage_Rk}, the following equalities hold in 
$A^\bullet(X)[E_1, \cdots, E_{n+1}]$: 
\begin{align}
R_1^{(n+1)} & = R_1^{(n)} \cup E_{n+1}\cdot r_1,
\label{eq:r1_equality}
\\
\left\langle R_2^{(n+1)} 
\right\rangle 
& = 
\left\langle
R_2^{(n)} \cup E_{n+1}\cdot r_2 
\right\rangle, 
\label{eq:r2_equality}
\\
R_3^{(n+1)} & = R_3^{(n)} \cup  \raisebox{-1.5mm}{$\Biggl\{$} P_{X,Y_{n+1}} \raisebox{-1.5mm}{$\Biggl($} -\sum_{Y_j \subseteq Y_{n+1}} E_j \raisebox{-1.5mm}{$\Biggr)$} \raisebox{-1.5mm}{$\Biggr\}$} ,
\label{eq:r3_equality}
\\
R_4^{(n+1)} & = R_4^{(n)} \cup E_{n+1} \cdot r_4,
\label{eq:r4_equality}
\end{align}
Moreover, if we denote by $r'_1$, $r'_2$ and $r'_4$ the images of $r_1$, $r_2$ and $r_4$ in 
\[
A^\bullet\left(X^{(n)}\right)
=
\frac{A^\bullet(X)[E_1,\ldots,E_n]}{
\left\langle
R_1^{(n)} , \, R_2^{(n)} , \, R_3^{(n)} , \, R_4^{(n)}
\right\rangle
}, 
\]
and denote $J_{n+1} := 
\ker\left(
A^\bullet\left( X^{(n)}\right) \rightarrow
A^\bullet\left( Y^{(n)}_{n+1} \right)
\right)$, 
then we have the equality
\begin{align}
J_{n+1} &= \left\langle r'_1 , r'_2 , r'_4 \right\rangle. 
\label{eq:kernel_equality}
\end{align}
\end{lemma}
Consider the blow-up construction of $X_\mathcal{G}=X^{(n+1)}$, where $X^{(m+1)} = \operatorname{Bl}_{Y_{m+1}^{(m)}}X^{(m)}$ for each $0\leq m \leq n$, as in Definition~\ref{bg.def.transform.of.building.set}. At each one of these blow-ups, the element $Y_{m+1}^{(m)}\in\mathcal{G}^{(m)}$ is minimal, since all elements that come before it in the chosen order are divisors, and those that come after cannot be contained in it because the order of $\mathcal{G}$ respects containments. Then, a repeated application of Lemma~\ref{lemma.chow.surjectivity.after.minimal.blowup} shows that the pullback map $A^\bullet\left(X^{(n)}\right)\to A^\bullet\left(Y_{n+1}^{(n)}\right)$ is surjective. Hence Keel's formula (Lemma~\ref{lemma: keel's blowup formula}) can be used to compute the Chow ring of $X_\mathcal{G}=X^{(n+1)} = \operatorname{Bl}_{Y_{n+1}^{(n)}}X^{(n)}$. 
By Lemma~\ref{lemma.Chern.polynomials.blowup} and Lemma~\ref{lemma: useful lemmas from li's}\ref{lemma: useful lemmas from li's item 1}, a Chern polynomial for $Y_{n+1}^{(n)}$ on $X^{(n)}$ is given by 
\[
P_{n+1}(t):= P_{X,Y_{n+1}} \raisebox{-1.5mm}{$\Biggl($} t \ - \sum_{Y_j \subsetneq Y_{n+1}} E_j \raisebox{-1.5mm}{$\Biggr)$}  \in A^\bullet\left(X^{(n)}\right)[t]. 
\]
Moreover, $Y_{n+1}^{(n)}$ is not contained in any of the divisors $Y_i^{(n)}$ with $i\leq n$, since otherwise taking images in $X$ would give $Y_{n+1}\subseteq Y_i$, contrary to the chosen order. Hence the pullback to $X^{(n+1)}$ of the class of $Y_i^{(n)}$ is the class of its strict transform $Y_i^{(n+1)}$, so the identification of each variable $E_i$ with a boundary divisor class is compatible with the blowup. 
The theorem then follows from Lemma~\ref{lemma:equalities_R_k}: 

\begin{align*}
 A^\bullet \left(X^{(n+1)}  \right) &=
\frac{ A^\bullet \left(X^{(n)} \right)[E_{n+1}]}
{ \left\langle E_{n+1}\cdot J_{n+1} , P_{n+1}(-E_{n+1}) \right\rangle }  
\hspace{7.8em}
\text{by Keel's formula}
\\[4mm]
&=
\frac{ A^\bullet \left(X^{(n)} \right)[E_{n+1}]}
{ \left\langle E_{n+1}\cdot \left( r'_1 \cup  r'_2 \cup  r'_4 \right) , P_{n+1}(-E_{n+1}) \right\rangle }
\hspace{4em}
\text{by Equation~\eqref{eq:kernel_equality}}
\\[4mm]
 &=
\frac{ \raisebox{-1.5mm}{$\Biggl($}
 \frac{A^\bullet(X)[E_1, \ldots , E_n]}
  { \left\langle R_1^{(n)}, \, R_2^{(n)}, \,  R_3^{(n)}, \, R_4^{(n)} \right\rangle } \raisebox{-1.5mm}{$\Biggr)$}[E_{n+1}]}
{  \left\langle E_{n+1}\cdot \left( r'_1 \cup  r'_2 \cup  r'_4 \right) , P_{n+1}(-E_{n+1}) \right\rangle } 
\hspace{4em}
\text{by induction hypothesis}
 \\[4mm]
 &=
\frac{A^\bullet(X)[E_1, \cdots , E_{n+1}]}
 { \left\langle R_1^{(n)}, \, R_2^{(n)}, \, R_3^{(n)}, \, R_4^{(n)},   
 E_{n+1}\cdot \left( r_1 \cup  r_2 \cup  r_4 \right),
 \, P_{X,Y_{n+1}}\left(-\sum_{Y_j \subseteq Y_{n+1}} E_j \right)
 \right\rangle 
 }
  \\[4mm]
 &=
\frac{A^\bullet(X)[E_1, \cdots , E_{n+1}]}
 { \left\langle R_1^{(n)} \cup E_{n+1} \cdot r_1, \, R_2^{(n)} \cup E_{n+1} \cdot r_2, \, R_3^{(n)} \cup 
 \left\{\ P_{X,Y_{n+1}}\raisebox{-1.5mm}{$\Biggl($}-\sum_{Y_j \subseteq Y_{n+1}} E_j \raisebox{-1.5mm}{$\Biggr)$} \right\}, 
 \, R_4^{(n)} \cup E_{n+1} \cdot r_4 
 \right\rangle 
 }
  \\[4mm]
&=
  \frac{A^\bullet(X)[E_1, \cdots , E_{n+1}]}
 { \left\langle R_1^{(n+1)}, \, R_2^{(n+1)}, \,  R_3^{(n+1)}, \, R_4^{(n+1)} \right\rangle }
\hspace{5.5em}
\text{by~\eqref{eq:r1_equality}, \eqref{eq:r2_equality},
\eqref{eq:r3_equality}, and \eqref{eq:r4_equality}},
\end{align*} 

\noindent which would complete the inductive step. It remains to prove Lemma~\ref{lemma:equalities_R_k}, which we do by establishing Equations~\eqref{eq:r1_equality}--\eqref{eq:kernel_equality} in turn.

\smallskip

%%%%%%%%%%%%%%%%%%%%%%%%%%%%%%%%%%%%%%%%%%%%%%%%%%%%%

\textbf{Proof of Equation~\eqref{eq:r1_equality}}. Let us first prove that $R_1^{(n+1)}\subseteq R_1^{(n)}\cup E_{n+1}\cdot r_1$. Consider a subset $I\subseteq\{1,\dots,n+1\}$ that is not $\mathcal{G}_{n+1}$-nested, which defines an element $\prod_{i\in I}E_i$ of $R_1^{(n+1)}$. We consider two cases: either $n+1\in I$ or $n+1\not\in I$. 

First suppose that $n+1\in I$ and define $I'=I\smallsetminus \{n+1\}$. Then, since $I'\cup\{n+1\}$ is not $\mathcal{G}_{n+1}$-nested, it follows that $\prod_{i\in I'} E_i$ is in $r_1$. In particular, $\prod_{i\in I} E_i = E_{n+1}\cdot \prod_{i\in I'} E_i$ is in $E_{n+1}\cdot r_1$.

If $n+1\not\in I$, then $I\subseteq\{1,\dots,n\}$ and we claim that $I$ is not $\mathcal{G}_n$-nested, so that $\prod_{i\in I}E_i \in R_1^{(n)}$. Indeed, assume for a contradiction that $I$ is $\mathcal{G}_n$-nested and set $I=\{i_1,\dots,i_s\}$. Then, $Y_{i_1}^{(n)},\ldots,Y_{i_s}^{(n)}$ have a nonempty intersection $Z:=Y_{i_1}^{(n)}\cap\cdots\cap Y_{i_s}^{(n)}$ in $X^{(n)}$. 
Note that all elements of $\mathcal{G}_{n+1}^{(n)}$ but $Y_{n+1}^{(n)}$ are distinct divisors, and none of them contains $Y_{n+1}^{(n)}$ (because $Y_{n+1}^{(n)}  \subseteq Y_{j}^{(n)}$ on $X^{(n)}$ for some $1 \leq j \leq  n$ would imply $Y_{n+1}  \subseteq Y_{j}$ on $X$, which is not possible). In particular, each of $Y_{i_1}^{(n)},\ldots,Y_{i_s}^{(n)}$ is minimal in $\mathcal{G}_{n+1}^{(n)}$. 
By \cite[Lemma~2.5]{LiLi}, this implies that $Y_{i_1}^{(n)},\ldots,Y_{i_s}^{(n)}$ are the $\mathcal{G}^{(n)}_{n+1}$-factors of $Z$. Then, considering the one-step flag $\{Z\}$, we deduce that $I$ is $\mathcal{G}^{(n)}_{n+1}$-nested. By repeated application of \cite[Proposition~2.8(iii)]{LiLi}, this means that $I$ is $\mathcal{G}_{n+1}$-nested, which is a contradiction. Therefore, $I$ is not $\mathcal{G}_n$-nested. It follows that $\prod_{i\in I} E_i$ is in $R_1^{(n)}$.

We have proved that $R_1^{(n+1)}\subseteq R_1^{(n)}\cup E_{n+1}\cdot r_1$. We proceed to show that $R_1^{(n)}\cup E_{n+1}\cdot r_1\subseteq R_1^{(n+1)}$. This containment is significantly simpler. Indeed, $I\subseteq\{1,\dots,n\}$ is not $\mathcal{G}_n$-nested if and only if the boundary divisors associated to the subvarieties $\{Y_i\}_{i\in I}$ have empty intersection in $X^{(n)}$. If this is the case, then their corresponding boundary divisors in $X^{(n+1)}$ have empty intersection as well, so $I$ is not $\mathcal{G}_{n+1}$-nested. Similarly, every element in $E_{n+1}\cdot r_1$ is in $R_1^{(n+1)}$ by definition of $r_1$. This finishes the proof of Equation~\eqref{eq:r1_equality}.
\qed

\smallskip

\textbf{Proof of Equation~\eqref{eq:r2_equality}}. 
We first show that $R_2^{(n+1)}  \subseteq  R_2^{(n)} \cup E_{n+1}\cdot r_2 $. 

Consider the set $J_Z\cdot E_Z \subseteq R_2^{(n+1)}$ corresponding to $Z\in\mathcal{A}_{n+1} = \mathcal{A}$. 
Let $Y_{i_1}, \cdots, Y_{i_s}$ be the distinct $\mathcal{G}_{n+1}$-factors of $Z$. If $Y_{n+1}$ is not a $\mathcal{G}_{n+1}$-factor of $Z$, then $Z\in\mathcal{A}_n\subseteq\mathcal{A}_{n+1}$, and $Y_{i_1}, \cdots, Y_{i_s}$ are also the $\mathcal{G}_n$-factors of $Z$. It follows that $J_Z \cdot E_Z$ is contained in $R_2^{(n)}$. 
Now, suppose that $Y_{n+1}$ is a $\mathcal{G}$-factor of $Z$, say $Y_{n+1}=Y_{i_s}$.
Then $E_{n+1}=E_{i_s}$ is a factor of $E_{Z}$, so 
$J_Z \cdot E_Z =  J_Z\cdot E_{i_1} E_{i_2}\cdots E_{i_{s-1}}E_{n+1}$. 
If $s=1$, then $Z= Y_{n+1}$ and $J_Z \cdot E_Z = J_{Y_{n+1}}\cdot E_{n+1}$ is contained in $E_{n+1}\cdot r_2$. 
Similarly, if $s \geq 2$ we have $Z \subsetneq Y_{n+1}$ and $J_Z\cdot E_{i_1} E_{i_2}\cdots E_{i_{s-1}}$ is in $r_2$, hence 
\[
J_Z\cdot E_Z =  
\left(J_Z\cdot E_{i_1} E_{i_2}\cdots E_{i_{s-1}}\right)\cdot E_{n+1}
\subseteq E_{n+1} \cdot r_2.
\]

Then we have $\left\langle R_2^{(n+1)} 
\right\rangle
 \subseteq  
\left\langle
R_2^{(n)} \cup E_{n+1}\cdot r_2 
\right\rangle$. 
To prove the containment in the other direction it is enough to show that $R_2^{(n)}\subseteq R_2^{(n+1)}$ and $E_{n+1}\cdot r_2\subseteq\left\langle R_2^{(n+1)} \right\rangle$.

To prove that $R_2^{(n)}\subseteq R_2^{(n+1)}$, consider the set $J_Z\cdot E_Z \subseteq R_2^{(n)}$ corresponding to $Z\in\mathcal{A}_n\subseteq\mathcal{A}_{n+1}$. 
Let $Y_{i_1},\dots,Y_{i_s}$ be the distinct $\mathcal{G}_n$-factors of $Z$.  
Then, $Y_{i_1},\dots,Y_{i_s}$ are also the  $\mathcal{G}_{n+1}$-factors of $Z$ and
$J_Z\cdot E_Z \subseteq R_2^{(n+1)}$, as desired. 

To prove $E_{n+1}\cdot r_2\subseteq\left\langle R_2^{(n+1)} \right\rangle$, 
first consider $J_{Y_{n+1}} \subseteq r_2$, and notice that for $Z=Y_{n+1} \in \mathcal{A}_{n+1}$ we have $E_Z=E_{n+1}$ and hence
$J_{Y_{n+1}} \cdot E_{n+1} = J_Z\cdot E_Z \subseteq R_2^{(n+1)}$.  
Now, consider any $Z \in \mathcal{A}_{n+1}$ with $Z \subsetneq Y$ and the associated subset of $r_2$ given by $J_Z \cdot E_{i_1} \ldots E_{i_s} \subseteq r_2$, where $Y_{i_1}, \ldots, Y_{i_s}$ are the $\mathcal{G}_{n+1}$-factors of $Z$ with $i_1, i_2, \ldots, i_s \neq n+1$ (i.e. with $1 \leq i_1, i_2, \ldots, i_s \leq n$).
We need to show that $\left(J_Z\cdot E_{i_1} \cdots E_{i_s}\right) \cdot E_{n+1}$ is contained in $\left\langle R_2^{(n+1)}\right\rangle$, and we do this considering two cases depending on whether $Y_{n+1}$ is a $\mathcal{G}_{n+1}$-factor of $Z$. 

If $Y_{n+1}$ is a $\mathcal{G}_{n+1}$-factor of $Z$, then $Y_{i_1}, Y_{i_2}, \ldots, Y_{i_s}, Y_{n+1}$ are all the $\mathcal{G}_{n+1}$-factors of $Z$ and hence $E_Z=E_{i_1} E_{i_2} \cdots E_{i_s} \cdot E_{n+1}$. One then obtains that 
\[
\left(J_Z\cdot E_{i_1} \cdots E_{i_s}\right) \cdot E_{n+1}=J_Z\cdot E_Z \subseteq R_2^{(n+1)} \subseteq \left\langle R_2^{(n+1)}\right\rangle. 
\]
On the other hand, if $Y_{n+1}$ is not a $\mathcal{G}_{n+1}$-factor of $Z$, then $Y_{i_1}, Y_{i_2}, \ldots, Y_{i_s}$ are all the $\mathcal{G}_{n+1}$-factors of $Z$. This implies that  $E_Z=E_{i_1} E_{i_2} \cdots E_{i_s}$ and $J_Z\cdot E_Z \subseteq R_2^{(n+1)}$. As a consequence, one obtains
$$
\left(J_Z E_{i_1} \cdots E_{i_s}\right) \cdot E_{n+1}=E_{n+1} \cdot J_Z\cdot E_Z \subseteq E_{n+1} \cdot R_2^{(n+1)} \subseteq \left\langle R_2^{(n+1)}\right\rangle, 
$$
and this completes the proof of Equation~\eqref{eq:r2_equality}.  \qed

\smallskip

%%%%%%%%%%%%%%%%%%%%%%%%%%%%%%%%%%%%%%%%%%%%%%%%%%%%%

\textbf{Proof of Equation~\eqref{eq:r3_equality}}.
The conclusion is clear from the definition of the $R_3^{(m)}$ as follows, 
\begin{align*}
R_3^{(n+1)}
=
\left\{
P_{X,Y_i}\raisebox{-1.5mm}{$\Biggl($}-\sum_{Y_j \subseteq Y_i} E_j\raisebox{-1.5mm}{$\Biggr)$}
\;\bigg|\;
i \in \{1,\ldots,n+1\}
\right\} 
=
R_3^{(n)}
\cup
\left\{
P_{X,Y_{n+1}}\raisebox{-1.5mm}{$\Biggl($}-\sum_{Y_j \subseteq Y_{n+1}} E_j\raisebox{-1.5mm}{$\Biggr)$}
\right\}.
\\[-3ex]
&\qedhere
\end{align*}
\smallskip

\textbf{Proof of Equation~\eqref{eq:r4_equality}}. Let us first prove that $R_4^{(n+1)}$ is contained in $R_4^{(n)}\cup E_{n+1}\cdot r_4$. Fix subvarieties $Z\in\mathcal{A}_{n+1}$ and $Y_i\subsetneq Z$, and let $\alpha$ be their corresponding element in $R_4^{(n+1)}$, that is,
\[
    \alpha 
    = E_Z \cdot P_{X,Z,Y_i}\raisebox{-1.5mm}{$\Biggl($}-\sum_{Y_j\subseteq Y_i} E_j\raisebox{-1.5mm}{$\Biggr)$}.
\]
We consider two cases depending on whether $Y_{n+1}$ is a $\mathcal{G}_{n+1}$-factor of $Z$.

If $Y_{n+1}$ is not a $\mathcal{G}_{n+1}$-factor of $Z$, 
then the $\mathcal{G}_{n+1}$-factors of $Z$ are contained in $\mathcal{G}_n$. 
Hence, $Z$ is an intersection of elements of $\mathcal{G}_n$, 
and $Z\in\mathcal{A}_n$, so $\alpha\in R_4^{(n)}$. If $Y_{n+1}$ is a $\mathcal{G}_{n+1}$-factor of $Z$, then $\alpha$ can be written as
\[
    \alpha 
    = E_{n+1}\cdot\left[\frac{E_Z}{E_{n+1}}\cdot P_{X,Z,Y_i}\raisebox{-1.5mm}{$\Biggl($}-\sum_{Y_j\subseteq Y_i} E_j\raisebox{-1.5mm}{$\Biggr)$}\right]\in E_{n+1}\cdot r_4,
\]
where $E_Z$ is divisible by $E_{n+1}$ because $Y_{n+1}$ is a $\mathcal{G}_{n+1}$-factor of $Z$.

The containment $R_4^{(n)}\cup E_{n+1}\cdot r_4\subseteq R_4^{(n+1)}$ is significantly simpler. Indeed, $E_{n+1}\cdot r_4$ is contained in $R_4^{(n+1)}$ by definition of $r_4$. Similarly, $R_4^{(n)}\subseteq R_4^{(n+1)}$ since $\mathcal{G}_n\subseteq\mathcal{G}_{n+1}$ and $\mathcal{A}_n\subseteq\mathcal{A}_{n+1}$.
\qed
\smallskip

\textbf{Setup and strategy for proving Equation~\eqref{eq:kernel_equality}}. The goal is to compute the kernel $J_{n+1}$ of the pullback morphism $A^\bullet\left(X^{(n)}\right)\to A^\bullet\left(Y^{(n)}_{n+1}\right)$. We have that $X^{(n)}=X_{\mathcal{G}'}$, where $\mathcal{G}'=\mathcal{G}\smallsetminus \{Y\}$. Similarly, by Theorem~\ref{strict.transforms.are.wonderful}, we have that $Y^{(n)}_{n+1}$ is the wonderful compactification of the building set $\mathcal{G}|_Y$ on $Y$.  
Therefore, our goal is to compute the kernel of the restriction map $\iota^*:A^{\bullet}(X_{\mathcal{G}'})\to A^{\bullet}(Y_{\mathcal{G}|_Y})$ induced by the inclusion $\iota:Y_{\mathcal{G}|_Y}\hookrightarrow X_{\mathcal{G}'}$ from Theorem~\ref{strict.transforms.are.wonderful}.

The building sets $\mathcal{G}'$ on $X$ and $\mathcal{G}|_{Y}$ on $Y$ have cardinality at most $n$ and we claim that both satisfy the surjectivity assumptions in the statement. 
This is clear for $\mathcal{G}'$ since $\mathcal{G}' \subseteq \mathcal{G}$, and hence the arrangement of $\mathcal{G}'$ is a subset of $\mathcal{A}$. 
Similarly, the arrangement $\mathcal{A}_Y$ on $Y$ corresponding to $\mathcal{G}|_Y$ is a subset of $\mathcal{A}$, and for each $Z \in \mathcal{A}_Y$ the pullback morphism $A^\bullet(X) \rightarrow A^\bullet(Z)$ factors through the pullback morphism $A^\bullet(Y) \rightarrow A^\bullet(Z)$.
Since all such pullbacks $A^\bullet(X) \rightarrow A^\bullet(Z)$ are surjective by assumption, then the pullbacks $A^\bullet(Y) \rightarrow A^\bullet(Z)$ are surjective for all $Z$ in $\mathcal{A}_Y$.
Therefore, by the induction hypothesis, the Chow rings $A^{\bullet}(X_{\mathcal{G}'})$ and $A^\bullet\left(Y_{\mathcal{G}|_Y}\right)$ are described as in the statement as quotients of the polynomial rings $A^\bullet(X)[E_1,\ldots,E_n]$ and 
$A^{\bullet}(Y)[e_i \mid Y_i \in \mathcal{G}_Y]$ by the four types of relations (1), (2), (3) and (4) applied to the building sets $\mathcal{G}'$ and $\mathcal{G}|_Y$, respectively.

Next we lift $\iota^*:A^{\bullet}(X_{\mathcal{G}'})\to A^{\bullet}(Y_{\mathcal{G}|_Y})$ to the level of polynomial rings before the quotient in the presentations is taken. Consider the algebra homomorphism $\phi$ between the polynomial rings 
\begin{equation}
    \phi:A^\bullet(X)[E_1,\dots,E_n]\to A^{\bullet}(Y)[e_i \mid Y_i \in \mathcal{G}_Y]  \label{equation.morphism.polynomial.level}
\end{equation}
uniquely determined by the pullback $A^\bullet(X)\to A^\bullet(Y)$ on the coefficients, and mapping
\begin{align*}
    E_i\mapsto\left\{
    \begin{aligned}
        & e_i \ \  \text{ if $Y_i \in \mathcal{G}_Y$ \ \  (equivalently if $\{i,n+1\}$ is $\mathcal{G}_{n+1}$-nested)},\;\\
        & 0 \ \ \ \text{ otherwise}. 
    \end{aligned}
    \right.
\end{align*}
By Corollary~\ref{corollary.restriction.boundary.divisors} there is a commutative diagram of algebras 
\begin{equation}\label{diagram.proof.main.thm}
\begin{tikzcd}
A^{\bullet}(X)[E_1,\ldots,E_n]
    \arrow[r, "\phi"]
    \arrow[d, "\pi_X"', two heads]
&
A^{\bullet}(Y)[e_i \mid Y_i \in \mathcal{G}_Y]
    \arrow[d, "\pi_Y", two heads]
\\
A^{\bullet}(X_{\mathcal{G}'})
    \arrow[r, "\iota^*"]
&
A^{\bullet}(Y_{\mathcal{G}|_Y})
\end{tikzcd}
\end{equation}
where the homomorphisms $\pi_X$ and $\pi_Y$ are those uniquely determined by the pullback maps in Chow theory, the condition that 
$\pi_X(E_i) \in A^{\bullet}(X_{\mathcal{G}'})$ is the class of the boundary divisor of $X_{\mathcal{G}'}$ associated to $Y_i$, for each $Y_i \in \mathcal{G}'$, 
and that $\pi_Y(e_i) \in A^{\bullet}(Y_{\mathcal{G}|_Y})$ is the class of the boundary divisor of $Y_{\mathcal{G}|_Y}$ associated to $Y_i \cap Y$, for each $Y_i \in \mathcal{G}_Y$.

Since $\pi_X$ is surjective, we can obtain $J_{n+1}$ by first computing the kernel of $\pi_Y\circ\phi$ and then taking its image in $A^\bullet\left(X_{\mathcal{G}'}\right)$, that is, 
\[
    J_{n+1}=\ker\iota^* = \pi_X(\ker(\pi_Y\circ\phi)).
\]
The kernel of $\pi_Y\circ\phi$ will be computed using the following general result:

\begin{lemma} \label{lemma.kernel.strategy}
    Let $f: R \rightarrow S$ be a surjective ring homomorphism and $J = \langle  b_\alpha \, | \, \alpha \in A\rangle $ an ideal of $S$, where $A$ is some indexing set. 
    For each $\alpha \in A$, let $a_\alpha \in R$ such that $f(a_{\alpha})=b_{\alpha}$.
    Let $\pi: S \rightarrow S/J$ be the quotient homomorphism and let $g$ denote the composition $g=\pi \circ f : R \rightarrow S/J$. Then,
    \[
    \ker(g) = \ker(f) + \langle a_\alpha \, | \, \alpha \in A\rangle . 
    \]
\end{lemma}
\begin{proof}
    Given $r \in \ker(g)$, we have that $f(r) \in J$, hence $f(r)=s_1 b_{\alpha_1}+\ldots+s_m b_{\alpha_m}$ for some $s_1,\ldots,s_m \in S$ and some indices $\alpha_1,\ldots,\alpha_m \in A$. 
    Choose any $r_1,\ldots,r_m \in R$ such that $f(r_i)=s_i$ for each $1 \leq i \leq m$. 
    Hence, 
    \[f\left(r-\sum_{i=1}^m r_i a_{\alpha_i}\right)
    =f(r)-\sum_{i=1}^m f(r_i)f(a_{\alpha_i})
    =f(r)-\sum_{i=1}^m s_i b_{\alpha_i} = 0.
    \]
    Then, $r - \sum_{i=1}^m r_i a_{\alpha_i} \in  \ker(f)$, 
    and therefore $r \in \ker(f) + \langle a_\alpha \, | \, \alpha \in A\rangle$. 
    Conversely, the inclusion $\ker(f) + \langle a_\alpha \, | \, \alpha \in A\rangle  \subseteq \ker(g)$ clearly holds, and this completes the proof. 
\end{proof}

%%%%%%%%%%%%%%%%%%%%%%%%%%%%%%%%%%%%%%%%%%%%%%%%%%%%%%%%%%%%%%%%%%%%%%%%%%%%%%%%%%%%%%

Within the proof of Equation~\eqref{eq:kernel_equality} we will evaluate the images under $\phi$ of the evaluated Chern polynomials appearing in relations (3) and (4), and for this the following lemma will be useful. 

\begin{lemma}  \label{lemma.lifting.Chern.polynomials}
Let $\phi: A^{\bullet}(X)[E_1,\ldots,E_n]\rightarrow A^{\bullet}(Y)[e_i \mid Y_i \in \mathcal{G}_Y]$ be the homomorphism in \eqref{equation.morphism.polynomial.level},
mapping $E_i$ to $e_i$ if $Y_i \in \mathcal{G}_Y$ and to $0$ if $Y_i \not\in \mathcal{G}_Y$.  
Denote $y_j := Y_j \cap Y$ for each $Y_j \in\mathcal{G}_Y$. 
Then, for every $Y_i\in\mathcal{G}_Y$, we have
\[
\phi
\raisebox{-1.5mm}{$\Biggl($}
\sum_{\substack{1\leq j\leq n\\Y_j\subseteq Y_i}}E_j
\raisebox{-1.5mm}{$\Biggr)$}
=
\sum_{\substack{Y_j\in\mathcal{G}_Y\\Y_j\subseteq Y_i}}e_j
=
\sum_{\substack{y_j\in\mathcal{G}|_Y\\y_j\subseteq y_i}}e_j.
\]
\end{lemma}

\begin{proof}
We apply $\phi$ to the sum term by term. 
The terms indexed
by $Y_j\notin\mathcal{G}_Y$ restrict to zero, while those indexed by
$Y_j\in\mathcal{G}_Y$ restrict to $e_j$. 
Moreover, 
Proposition~\ref{proposition.structure.restricted.building.set} gives a bijection $Y_j \mapsto y_j$ between $\mathcal{G}_Y$ and $\mathcal{G}|_{Y} $ satisfying that 
$Y_j\subseteq Y_i$ 
if and only if
$y_j\subseteq y_i$, for all $Y_i,Y_j\in\mathcal{G}_Y$. 
From these observations the conclusion follows. 
\end{proof}

%%%%%%%%%%%%%%%%%%%%%%%%%%%%%%%%%%%%%%%%%%%%%%%%%%%%%%%%%%%%%%%%%%%%%%%%%%%%%%%%%%%%%%

Lastly,  before proceeding with the proof of Equation~\eqref{eq:kernel_equality}, we fix Chern polynomials on $Y$ compatible by restriction with those we have fixed on $X$.

%%%%%%%%%%%%%%%%%%%%%%%%%%%%%%%%%%%%%%%%%%%%%%%%%%%%%%%%%%%%%%%%%%%%%%%%%%%%%%%%%%%%%%

\textbf{Chern polynomials on $Y$.}
Consider any $y_i \in \mathcal{G}|_Y$ (in particular, $Y_i \in \mathcal{G}_Y$) and, if applicable, also consider any $Z$ in $\mathcal{A}_Y$ satisfying $y_i \subsetneq Z \subsetneq Y$.
For the remainder of this section let us fix specific Chern polynomials $P_{Y,y_i}(t) \in A^\bullet(Y)[t]$ and $P_{Y,Z,y_i}(t) \in A^\bullet(Y)[t]$, as follows, considering the two  cases: (i) $y_i$ is in $\mathcal{G}$, or  (ii) $y_i$ is not in $\mathcal{G}$.

\begin{itemize}[leftmargin=1em]
    \item If $y_i \in \mathcal{G}$, then $y_i = Y_i \in \mathcal{G}$ by Proposition~\ref{proposition.structure.restricted.building.set}, and we fix the choice 
\[
P_{Y,y_i}(t)=P_{Y,Y_i}(t):=P_{X,Y,Y_i}(t)|_Y \in A^\bullet(Y)[t], 
\]
which works by the definition of $P_{X,Y,Y_i}$, and also fix the choice 
\[
P_{Y,Z,y_i}(t)=P_{Y,Z,Y_i}(t):= P_{X,Z,Y_i}(t)|_Y \in A^\bullet(Y)[t], 
\]
which works by the functoriality of pullbacks. 

\item If $y_i \notin \mathcal{G}$, then $y_i = Y_i \cap Y$ by Proposition~\ref{proposition.structure.restricted.building.set}, and $Y_i$ and $Y$ are the $\mathcal{G}$-factors of $y_i$ by Lemma~\ref{transversality.within.building.set}. 
In particular, $Y_i$ and $Y$ intersect transversally. 
Then, by Lemma~\ref{lemma.chern.polynomial.identities}\ref{lemma.chern.polynomial.identities.ii} we can choose 
\[
P_{Y,y_i}(t) = P_{Y,Y_i \cap Y}(t) := 
P_{X,Y_i}(t)|_Y \in {A^\bullet(Y)[t]}.
\]
We have that $Y$ is minimal over $y_i$ in $\mathcal{G}$, as it is one of its $\mathcal{G}$-factors. 
Since $y_i \subsetneq Z \subsetneq Y$, then $Y$ is also minimal over $Z$ in $\mathcal{G}$ 
and hence 
$Y$ is a $\mathcal{G}$-factor of $Z$. 
Let $W$ be the intersection of the $\mathcal{G}$-factors of $Z$ other than $Y$. 
We claim that $Y_i \subsetneq W$. 
To see this, consider any $\mathcal{G}$-factor $W_j$ of $Z$, not equal to $Y$. 
Since $y_i \subseteq Z \subseteq W_j$, then $W_j$ must contain one of the $\mathcal{G}$-factors of $y_i$, which are $Y$ and $Y_i$. 
Since $W_j$ cannot contain $Y$, because they are distinct $\mathcal{G}$-factors of $Z$, then 
$W_j$ contains $Y_i$. 
Therefore, $Y_i \subseteq W$, and the containment is strict 
since equality would imply
$y_i = Y_i \cap Y = W \cap Y \supseteq Z \supsetneq y_i$, which is a contradiction. 
By Lemma~\ref{lemma.properties.transversal.intersections}\ref{lemma.properties.transversal.intersections.item.disjoint.subsets} we have that $W$ and $Y$ intersect transversally along $Z$. 
Since $Y$ is transverse to both $W$ and $Y_i$, Corollary~\ref{corollary.chern.polynomial.remove.transversal.factor} implies that the restriction of a Chern polynomial 
$P_{W,Y_i}(t)\in A^\bullet(W)[t]$
to $W\cap Y$ is a Chern polynomial for $Y_i\cap Y$ in $W \cap Y$.
Then, we can choose
\[
P_{Y,Z,y_i}(t) = P_{Y,W \cap Y, Y_i \cap Y}(t) := 
P_{X,W,Y_i}(t)|_Y \in {A^\bullet(Y)[t]}.
\]
\end{itemize}

%%%%%%%%%%%%%%%%%%%%%%%%%%%%%%%%%%%%%%%%%%%%%%%%%%%%%%%%%%%%%%%%%%%%%%%%%%%%%
\smallskip

\textbf{Proof of Equation~\eqref{eq:kernel_equality}}.
We want to show $J_{n+1} = \langle r'_1,r'_2,r'_4 \rangle$. 
Using the notation introduced in the setup for this proof above, we claim that it is enough to prove the equality 
\[
    \ker(\pi_Y\circ\phi) = \langle r_1,r_2,r_4, R^{(n)}_3, R^{(n)}_4 \rangle.
\]
Indeed, since $\pi_X$ in Diagram~\eqref{diagram.proof.main.thm} is surjective, the desired conclusion would follow at once:  
\[
    J_{n+1}=\ker\iota^* 
    = \pi_X(\ker(\pi_Y\circ\phi))
    = \langle \pi_X(r_1),\pi_X(r_2),\pi_X(r_4), \pi_X(R^{(n)}_3), \pi_X(R^{(n)}_4) \rangle
    = \langle r'_1,r'_2,r'_4 \rangle.  
\]

Let us first verify the containment $\langle r_1,r_2,r_4, R^{(n)}_3, R^{(n)}_4 \rangle \subseteq \ker(\pi_Y\circ\phi)$.
Notice that $R^{(n)}_3, R^{(n)}_4 \subseteq \ker(\pi_Y\circ\phi)$ by the commutativity of Diagram~\eqref{diagram.proof.main.thm}.  
To complete this containment, it suffices to show that the elements in $\phi(r_1),\phi(r_2),\phi(r_4) \subseteq A^{\bullet}(Y)[e_i \mid Y_i \in \mathcal{G}_Y]$ are either zero or relations of type (1)--(4) for the building set $\mathcal{G}|_Y$, since such relations lie automatically in $\ker\pi_Y$.

Consider an element $\alpha_1:=\prod_{i\in I}E_i$ in $r_1$, where $I\subseteq\{1,\dots,n\}$ and $I\cup\{n+1\}$ is not $\mathcal{G}_{n+1}$-nested. 
By the definition of $\phi$, 
if $Y_i \not\in \mathcal{G}_Y$ then $\phi(E_i)$ is zero, hence
if $\{Y_i \mid i \in I\}$ is not contained in $\mathcal{G}_Y$ then $\phi(\alpha_1)= \prod_{i\in I} \phi(E_i)$ is zero, as desired. 
On the other hand, if $\{Y_i \mid i \in I\} \subseteq \mathcal{G}_Y$, 
by Proposition~\ref{proposition.nested.ambient.vs.restriction} we have that $I$ is not $\mathcal{G}|_Y$-nested since $I \cup \{n+1\}$ is not $\mathcal{G}_{n+1}$-nested.
Since $I$ is not $\mathcal{G}|_{Y}$-nested, $\phi(\alpha_1)=\prod_{i\in I} e_i$ is a relation of type (1) for $\mathcal{G}|_Y$, as desired. 
We have shown $r_1 \subseteq \ker(\pi_Y\circ\phi)$.

Now consider an element $\alpha_2$ in $r_2$. 
If $\alpha_2$ is in $J_{Y_{n+1}}$, it maps to zero in $A^\bullet(Y)$ and hence $\phi(\alpha_2)$ is zero in $A^{\bullet}(Y)[e_i \mid Y_i \in \mathcal{G}_Y]$, as desired. 
If $\alpha_2$ is not in $J_{Y_{n+1}}$, then by the definition of $r_2$ there exists $Z \in \mathcal{A}_{n+1}$, such that  $Z \subsetneq Y_{n+1}$ and $\alpha_2$ is in 
\begin{equation} 
J_Z \cdot \prod_{\substack{
  1 \leq i \leq n \\
  Y_i \in \operatorname{Fact}(\mathcal{G}_{n+1},Z)
}} E_i.
\label{equation.r2.set}
\end{equation}
For each $1 \leq i \leq n$ such that $Y_{i} \in \operatorname{Fact}(\mathcal{G}_{n+1},Z)$, we have that either $Y_i \subsetneq Y$ or $\emptyset \neq Y_{i} \cap Y \not\in \mathcal{G}$, and then $Y_{i} \in \mathcal{G}_Y$. 
Moreover, since the pullback map $A^\bullet(X) \rightarrow A^\bullet(Y)$ is surjective, then  the image of $J_Z \subseteq A^\bullet(X)$ in $A^\bullet(Y)$ is precisely the kernel 
$J_{Y,Z}:=\operatorname{ker} \left( A^\bullet(Y) \rightarrow A^\bullet(Z) \right)$. 
Therefore, the image under $\phi$ of the set in \eqref{equation.r2.set} is given by 
\begin{equation}
J_{Y,Z} \cdot \prod_{\substack{
  1 \leq i \leq n \\
  Y_i \in \operatorname{Fact}(\mathcal{G}_{n+1},Z)
}} e_i.   
\label{equation.r2.set.restricted}
\end{equation}
Since $Z\in\mathcal{A}_{n+1}$ and $Z\subsetneq Y$, we have $Z\in\mathcal{A}_Y$. 
Moreover, by Remark~\ref{remark.restriction.of.g.factors} the product in \eqref{equation.r2.set.restricted} is indexed precisely by the $\mathcal{G}|_Y$-factors of $Z$. 
Then, the set in \eqref{equation.r2.set.restricted} is a relation of type (2) for $\mathcal{G}|_Y$, as desired. 
We have shown $r_2 \subseteq \ker(\pi_Y\circ\phi)$. 

Next, consider an element $\alpha_4 := \frac{E_Z}{E_{n+1}} \cdot P_{X,Z,Y_{i_0}}\left(-\sum_{Y_j\subseteq Y_{i_0}} E_j\right)$ in $r_4$, where 
$Z\in\mathcal{A}_{n+1}$ is such that $Y \in \operatorname{Fact}(\mathcal{G}_{n+1},Z)$ and $Y_{i_0} \subsetneq Z$. 
Suppose first that $Z=Y$. Then $\operatorname{Fact}(\mathcal{G}_{n+1},Z)=\{Y\}$, so $E_Z/E_{n+1}=1$. Moreover, $Y_{i_0}\subsetneq Y$, hence $Y_{i_0}\in\mathcal{G}_Y$ and $y_{i_0}=Y_{i_0}\cap Y=Y_{i_0}$. By our choice of $P_{Y,y_{i_0}}(t)=P_{X,Y,Y_{i_0}}(t)|_Y$ and Lemma~\ref{lemma.lifting.Chern.polynomials}, we obtain
\begin{equation}
    \phi(\alpha_4)=P_{Y,y_{i_0}}\raisebox{-1.5mm}{$\Biggl($}
-\sum_{\substack{y_j \subseteq y_{i_0} \\ y_j \in \mathcal{G}|_Y}} e_j
\raisebox{-1.5mm}{$\Biggr)$}.
\label{equation.r3.a} 
\end{equation}
The expression in \eqref{equation.r3.a} is a relation of type (3) for $\mathcal{G}|_Y$, so $\alpha_4\in\ker(\pi_Y\circ\phi)$.

We now assume that $Z\subsetneq Y$. Since the map $\phi$ from $A^\bullet(X)[E_1,\ldots,E_n]$ to $A^{\bullet}(Y)[e_i \mid Y_i \in \mathcal{G}_Y]$ is an algebra homomorphism, to evaluate the image of $\alpha_4$ we can first evaluate separately the image of the polynomial, of its argument, and of its preceding factor, and put them back together. 
Let $Y_{i_1},\ldots,Y_{i_s},Y$ be the distinct $\mathcal{G}$-factors of $Z$. 
We have that $Z$ is in the arrangement $\mathcal{A}_Y$ of $\mathcal{G}|_Y$ and by Remark~\ref{remark.restriction.of.g.factors} the $\mathcal{G}|_Y$-factors of $Z$ are $Y_{i_1} \cap Y,\ldots, Y_{i_s} \cap Y$. 
For each $1 \leq j \leq s$, we have $\emptyset \neq Y_{i_j} \cap Y \not\in \mathcal{G}$, so $Y_{i_j} \in \mathcal{G}_Y$. 
It follows that the image of $\frac{E_Z}{E_{n+1}}$ in $A^{\bullet}(Y)[e_i \mid Y_i \in \mathcal{G}_Y]$ is precisely $e_Z$, the product of the $e_i$ corresponding to the $\mathcal{G}|_Y$-factors of $Z$. 
Also, since $Y_{i_0}\subsetneq Y$, we have $Y_{i_0}\in\mathcal{G}_Y$ and $y_{i_0}=Y_{i_0}\cap Y=Y_{i_0}$. 
By our choice of $P_{Y,Z,y_{i_0}}(t)=P_{Y,Z,Y_{i_0}}(t)=P_{X,Z,Y_{i_0}}(t)|_Y$ and Lemma~\ref{lemma.lifting.Chern.polynomials}, the image $\phi(\alpha_4)$ of $\alpha_4$ in $A^{\bullet}(Y)[e_i \mid Y_i \in \mathcal{G}_Y]$ is 
\begin{equation}
e_Z \cdot P_{Y,Z,y_{i_0}}
\raisebox{-1.5mm}{$\Biggl($}
-\sum_{\substack{y_j \subseteq y_{i_0} \\ y_j \in \mathcal{G}|_Y}} e_j
\raisebox{-1.5mm}{$\Biggr)$}. 
\label{equation.r3.b}
\end{equation} 
The expression in \eqref{equation.r3.b} is a relation of type (4) for $\mathcal{G}|_Y$, as desired. 
Combining the two cases $Z=Y$ and $Z \subsetneq Y$, we have shown $r_4 \subseteq \ker(\pi_Y\circ\phi)$.

%%%%%%%%%%%%%%%%%%%%%%%%%%%%%%%%%%%%%%%%%%%%%%%%%%%%%%%%%%%%%%%%%%%%%%%%%%%%%

Therefore, we have proved the containment $\langle r_1,r_2,r_4, R^{(n)}_3, R^{(n)}_4 \rangle \subseteq \ker(\pi_Y\circ\phi)$ and we now proceed to prove the other containment. 
Recalling Diagram~\eqref{diagram.proof.main.thm}, $\ker(\pi_Y\circ\phi)$ is precisely the kernel of the composition on the top row of the diagram 
\[
\begin{tikzcd}
A^\bullet(X)[E_1,\ldots,E_n]
\ar[r,twoheadrightarrow,"\phi"] \ar[d,twoheadrightarrow,"\pi_X"'] &
A^{\bullet}(Y)[e_i \mid Y_i \in \mathcal{G}_Y]
\ar[r,twoheadrightarrow,"\pi_Y"] &
A^\bullet\left(Y^{{(n)}}_{n+1}\right) 
=
\frac{A^{\bullet}(Y)[e_i \mid Y_i \in \mathcal{G}_Y]}{\left\langle \text{Relations (1), (2), (3), (4) on $\mathcal{G}|_Y$} \right\rangle}\\
A^\bullet\left(X^{(n)}\right) \ar[rru,twoheadrightarrow,"\iota^*"'] 
\end{tikzcd}.
\]

Applying Lemma~\ref{lemma.kernel.strategy} with $f=\phi$ and $\pi = \pi_Y$, it follows that $\ker(\pi_Y\circ\phi)$ is generated by $K_0, K_1, K_2, K_3, K_4 \subseteq A^\bullet(X)[E_1,\ldots,E_n]$, where 
\[
K_0 := \ker \phi =  \ker\left( A^\bullet(X)[E_1,\ldots,E_n]
\xlongrightarrow{\phi}
A^{\bullet}(Y)[e_i \mid Y_i \in \mathcal{G}_Y] \right)
\]
and $K_1$, $K_2$, $K_3$, $K_4$ are any lifts (i.e. preimages) under $\phi$ of the relations of types (1), (2), (3), (4) on $\mathcal{G}|_Y$
from $A^{\bullet}(Y)[e_i \mid Y_i \in \mathcal{G}_Y]$ to $A^\bullet(X)[E_1,\ldots,E_n]$. 
Therefore, the proof of the containment $\ker(\pi_Y\circ\phi) \subseteq \langle r_1,r_2,r_4, R^{(n)}_3, R^{(n)}_4 \rangle$, and also the proof of the theorem, will be complete if we show that $K_0 \subseteq \langle r_1,r_2,r_4, R^{(n)}_3, R^{(n)}_4 \rangle$ and that 
we can choose the lifts 
$K_1$, $K_2$, $K_3$, $K_4$ 
inside $\langle r_1,r_2,r_4, R^{(n)}_3, R^{(n)}_4 \rangle$. 
We do this one by one below, sometimes showing that a subset of $\langle r_1,r_2,r_4, R^{(n)}_3, R^{(n)}_4 \rangle$ suffices.

\underline{$K_0\subseteq\langle r_1,r_2\rangle$.} We have that $K_0$ is generated as an ideal of $A^\bullet(X)[E_1,\ldots,E_n]$ by $J_{Y_{n+1}}$ and the set 
\begin{equation}
\left\{
E_i \, | \, 1\leq  i \leq n, \text{ and } \{i, n+1\} \text{ is not } \mathcal{G}_{n+1}\text{-nested}  
\right\}. 
\label{equation.variables.mapping.to.zero}
\end{equation}
Since $J_{Y_{n+1}}$ is contained in $r_2$ and the set in \eqref{equation.variables.mapping.to.zero} is contained in $r_1$, it follows that
$K_0 \subseteq \langle r_1,r_2 \rangle$. 
\smallskip 

\underline{We prove that we can choose $K_1\subseteq\langle r_1\rangle$.} Let us consider a relation $\beta_1$ of type (1) for the building set $\mathcal{G}|_Y$. 
Such a relation has the form $\beta_1= e_{i_1} \cdots e_{i_s}$ where $I :=\{i_1,\ldots,i_s \} \subseteq \{1,\ldots,n\}$ is not $\mathcal{G}|_Y$-nested. 
Since the indices $i_1,\ldots,i_s$ correspond to elements of the restricted building set $\mathcal{G}|_Y$, then $\{Y_i \mid i\in I \} \subseteq \mathcal{G}_Y$.
By Proposition~\ref{proposition.nested.ambient.vs.restriction} we deduce that 
$I \cup \{ n+1 \}$ is not $\mathcal{G}$-nested, since $\{Y_i \mid i\in I \} \subseteq \mathcal{G}_Y$ and $I$ is not $\mathcal{G}|_Y$-nested.
Then, 
$E_{i_1} \cdots E_{i_s} \in A^\bullet(X)[E_1,\ldots,E_n]$ is an element in $r_1$ and its image in $A^{\bullet}(Y)[e_i \mid Y_i \in \mathcal{G}_Y]$ is precisely $\beta_1$.
Hence, we can take the lifts
$K_1$ from $\langle r_1 \rangle$, as claimed.
\smallskip 

\underline{We prove that we can choose $K_2\subseteq\langle  r_2 \rangle$.} Let us consider a relation $\beta_2$ of type (2) for the building set $\mathcal{G}|_Y$. 
Then, there exists an element $Z \in \mathcal{A}_Y$ with distinct $\mathcal{G}|_Y$-factors $y_{i_j} := Y_{i_j} \cap Y$ for some $\{Y_{i_1},\ldots,Y_{i_s}\} \subseteq \mathcal{G}_{Y}$, 
such that $\beta_2$ is in the set
\[
J_{Y,Z} \cdot e_{i_1} \cdots e_{i_s} 
\]
where $J_{Y,Z} := \ker(A^\bullet(Y) \rightarrow A^\bullet(Z))$. 
Since the pullback $A^\bullet(X) \rightarrow A^\bullet(Y)$ is surjective and the pullback $A^\bullet(X) \rightarrow A^\bullet(Z)$ factors as $A^\bullet(X) \rightarrow A^\bullet(Y) \rightarrow A^\bullet(Z)$, 
we deduce that all elements of $J_{Y,Z}$ can be lifted to elements of $J_Z$, and conversely $J_Z$ gets mapped into $J_{Y,Z}$. 
We have that $Z \in \mathcal{A}_{n+1}$ and by Remark~\ref{remark.restriction.of.g.factors} the distinct $\mathcal{G}_{n+1}$-factors of $Z$ are $Y_{i_1},\ldots, Y_{i_s}$, and possibly also $Y_{n+1}$,
but in any case the set $J_Z \cdot E_{i_1} \cdots E_{i_s} \in A^\bullet(X)[E_1,\ldots,E_n]$ is contained in $r_2$ and contains a lift of $\beta_2$ to $A^\bullet(X)[E_1,\ldots,E_{n}]$, since $E_{i_1} \cdots E_{i_s}$
maps to $e_{i_1} \cdots e_{i_s}$ and $J_Z$ maps onto $J_{Y,Z}$.  
Hence, we can take the lifts
$K_2$ from $\langle r_2 \rangle$, as claimed.
\smallskip

\underline{We prove that we can choose $K_3\subseteq\langle r_4 , R^{(n)}_3 \rangle$.}  
This is done by constructing a lift $\beta_3^{'} \in\langle r_4, R^{(n)}_3 \rangle \subseteq A^\bullet(X)[E_1,\dots,E_n]$ of any relation $\beta_3 \in A^{\bullet}(Y)[e_i \mid Y_i \in \mathcal{G}_Y]$ of type (3) on $\mathcal{G}|_Y$. 
The element $\beta_3$ has the form
\[
\beta_3 = P_{Y, y_i}
\raisebox{-1.5mm}{$\Biggl($} -\sum_{y_j \subseteq y_i} e_j \raisebox{-1.5mm}{$\Biggr)$}
\]
where $y_i = Y_i \cap Y \in \mathcal{G}|_{Y}$  
and the sum runs over all $y_j \in \mathcal{G}|_{Y}$ satisfying $y_j \subseteq y_i$.  
 We consider two cases: (i) $y_{i} $ is in the building set $\mathcal{G}$, or (ii) $y_{i} $ is not in the building set $\mathcal{G}$. 

We start with case (i). If $y_{i}= Y_i \cap Y$ is in the building set $\mathcal{G}$, since $Y_i \in \mathcal{G}_Y$ then $Y_i \subsetneq Y$ 
by Proposition~\ref{proposition.structure.restricted.building.set}. In particular, in this case $y_i = Y_i \cap Y = Y_i$. 
By our choice of $P_{Y, y_i}(t)$ and 
Lemma~\ref{lemma.lifting.Chern.polynomials}, a lifting of $\beta_3$ to $A^\bullet(X)[E_1,\ldots,E_n]$ is given by 
\[
\beta'_3  
= P_{X,Y, y_i}\raisebox{-1.5mm}{$\Biggl($} -\sum_{Y_j \subseteq Y_i} E_j \raisebox{-1.5mm}{$\Biggr)$}
= \frac{E_{n+1}}{E_{n+1}} \cdot P_{X,Y, Y_i}\raisebox{-1.5mm}{$\Biggl($} -\sum_{Y_j \subseteq Y_i} E_j \raisebox{-1.5mm}{$\Biggr)$}. 
\]
Since $Y_i \subsetneq Y$, then $\beta'_3$ is in $r_4$ and it can be taken as the lift of $\beta_3$ in $\langle r_4, R^{(n)}_3 \rangle$ we are looking for. 
This concludes the case when $y_i$ is in the building set $\mathcal{G}$.

Let us now consider the case (ii) in which $y_i = Y_i \cap Y \not\in\mathcal{G}$. 
By Lemma~\ref{transversality.within.building.set}, in this case $Y$ and $Y_i$ are the $\mathcal{G}$-factors of $y_i$, and in particular $Y$ and $Y_i$ intersect transversally on $X$.
Then, $\beta_3$ can be written as
\[
    \beta_3 
    = P_{Y, y_i}\raisebox{-1.5mm}{$\Biggl($} -\sum_{y_j \subseteq y_i} e_j \raisebox{-1.5mm}{$\Biggr)$}
    = P_{Y, Y_i\cap Y}\raisebox{-1.5mm}{$\Biggl($} -\sum_{y_j \subseteq y_i} e_j \raisebox{-1.5mm}{$\Biggr)$}.
\]
By our choice of $P_{Y, y_i}(t)$ and  Lemma~\ref{lemma.lifting.Chern.polynomials}, a lifting of $\beta_3$ to $A^\bullet(X)[E_1,\ldots,E_n]$ is given by 
\[
\beta'_3 
= P_{X, Y_i}
\raisebox{-1.5mm}{$\Biggl($} -\sum_{Y_j \subseteq Y_i} E_j \raisebox{-1.5mm}{$\Biggr)$}. 
\]
We observe that $\beta'_3$ is in $R^{(n)}_3$, so it can be taken as the lift of $\beta_3$ in $\langle r_4, R^{(n)}_3 \rangle$ we are looking for. 
This concludes the case when $y_i$ is not in the building set $\mathcal{G}$.

\smallskip

%%%%%%%%%%%%%%%%%%%%%%%%%%%%%%%%%%%%%%%%%%%%%%%%%%%%%

\underline{We prove that we can choose $K_4\subseteq\langle r_4 , R^{(n)}_4 \rangle $.}  
This is done by constructing a lift $\beta_4^{'} \in\langle r_4, R^{(n)}_4 \rangle \subseteq A^\bullet(X)[E_1,\dots,E_n]$ of any relation $\beta_4 \in A^{\bullet}(Y)[e_i \mid Y_i \in \mathcal{G}_Y]$ of type (4) on $\mathcal{G}|_Y$. 
The element $\beta_4$ has the form
\[
\beta_4 = e_Z \cdot P_{Y,Z, y_i}
\raisebox{-1.5mm}{$\Biggl($} -\sum_{y_j \subseteq y_i} e_j \raisebox{-1.5mm}{$\Biggr)$}
\]
where $Z \in \mathcal{A}_Y$, $y_i = Y_i \cap Y \in \mathcal{G}|_{Y}$ satisfies $y_i \subsetneq Z$, the sum runs over all $y_j \in \mathcal{G}|_{Y}$ satisfying $y_j \subseteq y_i$, and 
$e_Z$ is the product of the $e_j$ such that $y_j$ is a $\mathcal{G}|_{Y}$-factor of $Z$.  
 We consider two cases: (i) $y_{i} $ is in the building set $\mathcal{G}$, or (ii) $y_{i} $ is not in the building set $\mathcal{G}$.

We start with case (i). If $y_{i}= Y_i \cap Y$ is in the building set $\mathcal{G}$, since $Y_i \in \mathcal{G}_Y$ then $Y_i \subsetneq Y$ 
by Proposition~\ref{proposition.structure.restricted.building.set}. In particular, in this case $y_i = Y_i \cap Y = Y_i$.
We consider two subcases depending on whether $Y$ is a $\mathcal{G}$-factor of $Z$. 
If $Y$ is a $\mathcal{G}$-factor of $Z$, let $Y_{i_1},\ldots,Y_{i_s},Y$ be the distinct $\mathcal{G}$-factors of $Z$. 
Then, by Remark~\ref{remark.restriction.of.g.factors} the $\mathcal{G}|_Y$-factors of $Z$ are 
$y_{i_1},\ldots,y_{i_s}$ and we have 
\[
e_Z = e_{i_1}\cdots e_{i_s} 
\qquad 
\text{ and }
\qquad 
E_Z = E_{i_1}\cdots E_{i_s} \cdot E_{n+1}.  
\]
By our choice of $P_{Y,Z, y_i}(t)$ and Lemma~\ref{lemma.lifting.Chern.polynomials}, a lifting of $\beta_4$ to $A^\bullet(X)[E_1,\ldots,E_n]$ is given by 
\[
\beta'_4  
=E_{i_1}\cdots E_{i_s} \cdot P_{X,Z, y_i}\raisebox{-1.5mm}{$\Biggl($} -\sum_{Y_j \subseteq Y_i} E_j \raisebox{-1.5mm}{$\Biggr)$}
=\frac{E_Z}{E_{n+1}} \cdot P_{X,Z, Y_i}\raisebox{-1.5mm}{$\Biggl($} -\sum_{Y_j \subseteq Y_i} E_j \raisebox{-1.5mm}{$\Biggr)$}. 
\]
Since $Y$ is a $\mathcal{G}$-factor of $Z$ and $Y_i \subsetneq Z$, then $\beta'_4$ is in $r_4$ and it can be taken as the lift of $\beta_4$ in $\langle r_4, R^{(n)}_4 \rangle$ we are looking for. 

If $Y$ is not a $\mathcal{G}$-factor of $Z$, let $Y_{i_1},\ldots,Y_{i_s}$ be the distinct $\mathcal{G}$-factors of $Z$. 
Then, by Remark~\ref{remark.restriction.of.g.factors} the $\mathcal{G}|_Y$-factors of $Z$ are 
$y_{i_1},\ldots,y_{i_s}$ and we have 
\[
e_Z = e_{i_1}\cdots e_{i_s} 
\qquad 
\text{ and }
\qquad 
E_Z = E_{i_1}\cdots E_{i_s}.  
\]
Then, by our choice of $P_{Y,Z, y_i}(t)$ and by Lemma~\ref{lemma.lifting.Chern.polynomials} a lifting of $\beta_4$ to $A^\bullet(X)[E_1,\ldots,E_n]$ is given by 
\[
\beta'_4 = E_{i_1}\cdots E_{i_s} \cdot P_{X,Z, y_i}
\raisebox{-1.5mm}{$\Biggl($}-\sum_{Y_j \subseteq Y_i} E_j \raisebox{-1.5mm}{$\Biggr)$}
= E_Z \cdot P_{X,Z, Y_i}
\raisebox{-1.5mm}{$\Biggl($} -\sum_{Y_j \subseteq Y_i} E_j \raisebox{-1.5mm}{$\Biggr)$}. 
\]
We observe that $\beta'_4$ is in $R^{(n)}_4$, so it can be taken as the lift of $\beta_4$ in $\langle r_4, R^{(n)}_4 \rangle$ we are looking for. 
This concludes the case when $y_i$ is in the building set $\mathcal{G}$.

Let us now consider the case (ii) in which $y_i = Y_i \cap Y \not\in\mathcal{G}$. 
By Lemma~\ref{transversality.within.building.set}, in this case $Y$ and $Y_i$ are the $\mathcal{G}$-factors of $y_i$, and in particular $Y$ and $Y_i$ intersect transversally on $X$. 
We have that $Y$ is minimal over $y_i$ in $\mathcal{G}$, as it is one of its $\mathcal{G}$-factors. 
Since $y_i \subseteq Z \subsetneq Y$, then $Y$ is also minimal over $Z$ in $\mathcal{G}$ 
and hence 
$Y$ is a $\mathcal{G}$-factor of $Z$. 
Let $Y_{i_1},\ldots,Y_{i_s},Y$ be the distinct $\mathcal{G}$-factors of $Z$. 
Notice that $s \geq 1$ since $Z \subsetneq Y$ (where $Z\in\mathcal{A}_Y$ guarantees the strict containment). 
Then, by Remark~\ref{remark.restriction.of.g.factors} the $\mathcal{G}|_Y$-factors of $Z$ are 
$y_{i_1},\ldots,y_{i_s}$ and we have 
\[
e_Z = e_{i_1}\cdots e_{i_s} 
\qquad 
\text{ and }
\qquad 
E_Z = E_{i_1}\cdots E_{i_s} \cdot E_{n+1}.  
\]
Define $W=Y_{i_1}\cap\cdots\cap Y_{i_s}$, so that $Z = W \cap Y$ and this intersection is transverse on $X$. 
By \cite[Lemma~2.4(i)]{LiLi} $Y_{i_1},\ldots,Y_{i_s}$ are the $\mathcal{G}$-factors of $W$, 
and in particular $E_W = E_{i_1} \cdots E_{i_s}$. 
Then, $\beta_4$ can be written as
\[
    \beta_4 
    = e_Z \cdot P_{Y,Z, y_i}\raisebox{-1.5mm}{$\Biggl($} -\sum_{y_j \subseteq y_i} e_j \raisebox{-1.5mm}{$\Biggr)$}
    = e_{i_1}\cdots e_{i_s}  \cdot P_{Y,W\cap Y, Y_i\cap Y}\raisebox{-1.5mm}{$\Biggl($} -\sum_{y_j \subseteq y_i} e_j \raisebox{-1.5mm}{$\Biggr)$}.
\]

Then, by our choice of $P_{Y,Z, y_i}(t)$ and Lemma~\ref{lemma.lifting.Chern.polynomials} a lifting of $\beta_4$ to $A^\bullet(X)[E_1,\ldots,E_n]$ is given by 
\[
\beta'_4  
= E_{i_1}\cdots E_{i_s}  \cdot P_{X,W, Y_i}
\raisebox{-1.5mm}{$\Biggl($} -\sum_{Y_j \subseteq Y_i} E_j \raisebox{-1.5mm}{$\Biggr)$}
= E_{W} \cdot P_{X,W, Y_i}
\raisebox{-1.5mm}{$\Biggl($} -\sum_{Y_j \subseteq Y_i} E_j \raisebox{-1.5mm}{$\Biggr)$}. 
\]
We observe that $\beta'_4$ is in $R^{(n)}_4$, so it can be taken as the lift of $\beta_4$ in $\langle r_4, R^{(n)}_4 \rangle$ we are looking for. 
This concludes the case when $y_i$ is not in the building set $\mathcal{G}$.

This completes the proof of Equation~\eqref{eq:kernel_equality},
and hence also of Lemma~\ref{lemma:equalities_R_k} and Theorem~\ref{theorem.main.intro}.  \qed

%%%%%%%%%%%%%%%%%%%%%%%%%%%%%%%%%%%%%%%%%%%%%%%%%%%%%%%%%%%%%%%%
%%%%%%%%%%%%%%%%%%%%%%%%%%%%%%%%%%%%%%%%%%%%%%%%%%%%%%%%%%%%%%%%
%%%%                                                        %%%%
%%%%  5. Simplified presentations and sharpness of the Theorem  
%%%%                                                        %%%%
%%%%%%%%%%%%%%%%%%%%%%%%%%%%%%%%%%%%%%%%%%%%%%%%%%%%%%%%%%%%%%%%
%%%%%%%%%%%%%%%%%%%%%%%%%%%%%%%%%%%%%%%%%%%%%%%%%%%%%%%%%%%%%%%%

\section{Simplified presentations and sharpness of Theorem~\ref{theorem.main.intro}} \label{section.examples}

Many prominent examples of wonderful compactifications in the literature come from building sets endowed with additional properties that simplify their Chow rings. These include (weighted) Fulton--MacPherson spaces \cite{Fulton-MacPherson,routis2014weighted,petersen2017chow}, the moduli spaces of pointed stable rooted trees of projective spaces of Chen--Gibney--Krashen \cite{Chen-Gibney-Krashen}, and the moduli spaces of rational curves with a cyclic action of Clader--Damiolini--Li--Ramadas \cite{clader2023wonderful}. Another large class is given by Hu's compactifications of \emph{regular simple arrangements of subvarieties} \cite{hu2003compactification}. In this section we identify properties responsible for some of these simplifications in the context of Theorem~\ref{theorem.main.intro}, and we provide a correction to the presentation in \cite[Theorem~1.8]{hu2003compactification} of the Chow ring of Hu's compactifications, showing that additional relations are needed in general, and identifying a natural hypothesis under which the presentation holds as stated there.

For an arbitrary building set Theorem~\ref{theorem.main.intro} is sharp in the following sense. Relations (2) and (4) cannot be indexed by $\mathcal{G}$ instead of $\mathcal{A}$, as shown in Examples~\ref{example.relations.2} and \ref{example.relations.4}, respectively. Relations (1) may be restricted to those arising from containment-minimal subsets that are not $\mathcal{G}$-nested, 
but a proper subfamily of these latter relations does not suffice in general, as Example~\ref{example.relations.1} illustrates.

\subsection{Conditions for simplified presentations} In this subsection we present a few  natural conditions on a building set and its arrangement that help simplify the presentation given in Theorem~\ref{theorem.main.intro}. 

\begin{definition}\label{definition.complete.regular.retractable}
    Let $\mathcal{G}$ be a building set on a smooth variety $X$ with induced arrangement $\mathcal{A}$. 
    \begin{itemize}
        \item We say that $\mathcal{G}$ is \emph{complete} if $\mathcal{G}=\mathcal{A}$.

        \item Let $Z,V,W\in\mathcal{A}$ with $Z=V\cap W$ and $Z\subsetneq V$, $Z\subsetneq W$. We say the decomposition $Z=V\cap W$ is \emph{retractable} if, after possibly renaming $V$ and $W$, there exists a morphism $r_V\colon X\to V$ with 
        \[
            r_V\circ i_{X,V}=\operatorname{id}_V\quad\text{and}\quad r_V(W)\subseteq Z,
        \]
        where $i_{X,V}:V\to X$ is the inclusion. We say that $\mathcal{A}$ is \emph{retractable} if every $Z\in\mathcal{A}\smallsetminus\mathcal{G}$ admits at least one retractable decomposition $Z=V\cap W$ with $V,W\in\mathcal{A}$.

        \item Given $Y\in\mathcal{G}$ and $Z\in\mathcal{A}$ with $Y\subsetneq Z$, a \emph{transversal witness} for the pair $(Y,Z)$ is an element $W\in\mathcal{G}$ such that $Y=W\cap Z$ is a transverse intersection. We say that $\mathcal{G}$ is \emph{transversally regular} if every such pair admits a transversal witness.
    \end{itemize}
\end{definition}

\begin{remark} 
    While a complete building set is vacuously retractable, no other implication holds among the conditions in Definition~\ref{definition.complete.regular.retractable}.
    \begin{enumerate}[label=(\roman*),ref=(\roman*)]
        \item Consider $\mathcal{G}=\{p,L\}$ with $p\in L$ a point on a line in $\mathbb{P}^2$. This is
        complete and (vacuously) retractable, but not transversally regular.
        \item In $\mathbb{P}^1\times\mathbb{P}^1$, take $\mathcal{G}=\{\mathbb{P}^1\times\{0\},
        \{0\}\times\mathbb{P}^1\}$. This is retractable, via the projections, but not complete.
        \item In $\mathbb{P}^2$, take the arrangement $\mathcal{A}=\{L_1,L_2,L_1~\cap~L_2\}$ given by two lines and their intersection point,
        with building set $\mathcal{G}=\{L_1,L_2\}$. This is transversally regular but neither complete nor retractable, since $\mathbb{P}^2$ admits no nonconstant map to $\mathbb{P}^1$.      
    \end{enumerate}
\end{remark}

\begin{proposition}\label{prop.complete.relations.1}
    Let $\mathcal{G}$ be a complete building set. Then, relations (1) in Theorem~\ref{theorem.main.intro} are equivalent to
    \begin{itemize}
        \item[(1')]  \ $E_Y\cdot E_{Y'}$ \quad for all $Y,Y'\in \mathcal{G}$ such that $\{Y,Y'\}$ is \emph{not} $\mathcal{G}$-nested.
    \end{itemize}
\end{proposition}
\begin{proof}
    Every relation of type (1') is of type (1), so it is enough to prove that if $\mathcal{G}$ is complete and $\mathcal{T}\subseteq\mathcal{G}$ is not $\mathcal{G}$-nested, then there exist $Y,Y'\in\mathcal{T}$ such that $\{Y,Y'\}$ is not $\mathcal{G}$-nested.

    For a complete building set the nested sets are precisely the chains. Thus, in this setting a subset $\mathcal{T}\subseteq\mathcal{G}$ is $\mathcal{G}$-nested if and only if every pair of its elements is $\mathcal{G}$-nested. The result follows from this.
\end{proof}

\begin{proposition}\label{prop.reduction.relations.2}
    Let $\mathcal{G}$ be a building set with retractable arrangement $\mathcal{A}$. Then, relations (2) in Theorem~\ref{theorem.main.intro} are equivalent to
    \begin{itemize}
        \item[(2')] \ \ $J_Y \cdot E_Y$  \quad for all $Y$ in $\mathcal{G}$.
    \end{itemize}
\end{proposition}

For the proof we use the following lemma. 

\begin{lemma} \label{lemma.kernels.generated.by.factors}
    Let $\mathcal{G}$ be a building set on a smooth variety $X$ whose induced arrangement
    $\mathcal{A}$ is retractable in the sense of Definition~\ref{definition.complete.regular.retractable}.
    Then, for every $Z\in\mathcal{A}$, 
    \[
            J_{X,Z} = \sum_{Y\in\operatorname{Fact}(\mathcal{G},Z)} J_{X,Y}.
    \]
\end{lemma}

\begin{proof}
    If $Y\in\mathcal{A}$ contains $Z$, then 
    the pullbacks under the inclusions satisfy 
    $i_{X,Z}^{*} = i_{Y,Z}^{*}\circ i_{X,Y}^{*}$ so $J_{X,Y}\subseteq J_{X,Z}$. It follows that the right-hand side is always contained in the left-hand side. For the reverse inclusion we argue by descending induction on the poset structure of $\mathcal{A}$ under containment. If $Z\in\mathcal{G}$, then $Z$ is its own unique $\mathcal{G}$-factor and the equality is immediate. 
    
    If $Z\in\mathcal{A}\smallsetminus\mathcal{G}$ let $Z = V\cap W$ be a retractable decomposition, with retraction $r_V:X\to V$ satisfying $r_V(W)\subseteq Z$. The claim will follow once we prove that $J_{X,Z} = J_{X,V}+J_{X,W}$. Indeed, since $Z\subsetneq V$ and $Z\subsetneq W$, the induction hypothesis gives
    \[
        J_{X,Z}
        = J_{X,V}+J_{X,W}
        = \sum_{Y\in\operatorname{Fact}(\mathcal{G},V)} J_{X,Y}
        + \sum_{Y\in\operatorname{Fact}(\mathcal{G},W)} J_{X,Y}
        \subseteq \sum_{Y\in\operatorname{Fact}(\mathcal{G},Z)} J_{X,Y},
    \]
    because every $G\in\operatorname{Fact}(\mathcal{G},V) \cup \operatorname{Fact}(\mathcal{G},W)$ contains some $G'\in\operatorname{Fact}(\mathcal{G},Z)$, so $J_{X,G}\subseteq J_{X,G'}$.
    
    To establish the equality $J_{X,Z} = J_{X,V}+J_{X,W}$, it suffices to prove the inclusion ``$\subseteq$", since the reverse inclusion is immediate. Consider $\alpha\in J_{X,Z}$ and define $a:=i_{X,V}^*(\alpha)$. Since $i_{X,Z} = i_{X,V}\circ i_{V,Z}$, functoriality gives $i_{V,Z}^{*}(a) = i_{X,Z}^{*}(\alpha) = 0$.

    Put $\theta := r_V^{*}(a)\in A^\bullet(X)$. The condition $r_V(W)\subseteq Z$ implies that $r_V\circ i_{X,W}$ factors as $i_{V,Z}\circ\rho_W$ for some morphism $\rho_W\colon W\to Z$ by \cite[Tag~0356]{stacks-project}. Then
    \[
        i_{X,W}^{*}(\theta) = \left(r_V\circ i_{X,W}\right)^{*}(a) = (\rho_W^{*}\circ i_{V,Z}^{*})(a) = 0,
    \]
    so $\theta\in J_{X,W}$. On the other hand, $i_{X,V}^{*}(\theta) = i_{X,V}^{*}(\alpha)$ because $r_V\circ i_{X,V}=\operatorname{id}_V$. Thus, $\alpha-\theta \in J_{X,V}$. Therefore $\alpha = (\alpha-\theta) + \theta \in J_{X,V}+J_{X,W}$ and the result follows.
\end{proof}

\begin{proof}[Proof of Proposition~\ref{prop.reduction.relations.2}]
    Clearly relations (2') are a subset of relations (2). Suppose $\mathcal{A}$ is a retractable arrangement and consider an element $Z\in\mathcal{A}$. Then, by Lemma~\ref{lemma.kernels.generated.by.factors}, 
    \[
        E_Z\cdot J_{X,Z} = \raisebox{-1.5mm}{$\Biggl($} \prod_{Y\in\operatorname{Fact}(\mathcal{G},Z)} E_Y \raisebox{-1.5mm}{$\Biggr)$} \cdot  \sum_{Y\in\operatorname{Fact}(\mathcal{G},Z)} J_{X,Y}.
    \]
    The right-hand side is clearly contained in the ideal $\langle E_Y\cdot J_{X,Y}\,\vert\, Y\in\mathcal{G}\rangle$ and the result follows.    
\end{proof}

The building sets of the weighted Fulton--MacPherson spaces satisfy a slightly weaker property than being transversally regular, which motivates the following definition.

\begin{definition}
Given $Y\in\mathcal{G}$ and $Z\in\mathcal{A}$ with $Y\subsetneq Z$, a \emph{weak transversal witness} for the pair $(Y,Z)$ is a smooth subvariety $W\subseteq X$ such that $Y=W\cap Z$ as a transverse intersection, and either $W\in\mathcal{G}$ or $W\cap Y_i\in\mathcal{G}$ for some $\mathcal{G}$-factor $Y_i$ of $Z$ 
(so every transversal witness is a weak transversal witness).  
We say that $\mathcal{G}$ \emph{admits weak transversal witnesses} if every such pair admits a weak transversal witness.
\end{definition}

\begin{lemma}\label{lemma.local.relations.4}
In the setting of Theorem~\ref{theorem.main.intro}, 
let $Y\in\mathcal{G}$ and $Z\in\mathcal{A}$ with $Y\subsetneq Z$, and suppose that the pair $(Y,Z)$ admits a weak transversal witness $W$. For each $T\in\mathcal G$, set $\sigma_T := \mathlarger{\sum}\nolimits_{\substack{Y'\in\mathcal{G}\\ Y'\subseteq T}} E_{Y'}$.  
\begin{enumerate}[label=(\roman*),ref=(\roman*)]
    \item\label{item.lemma.local.relations.4.i} If $W\in\mathcal{G}$ (i.e., if $W$ is a transversal witness), then the relation $E_Z\cdot P_{X,Z,Y}(-\sigma_Y)$ of type (4) associated to $Y\subsetneq Z$ lies in the ideal generated by the relations (1), (2) and (3). 
    \item\label{item.lemma.local.relations.4.ii} Suppose that $W\cap Y_i\in\mathcal{G}$ for some $\mathcal{G}$-factor $Y_i$ of $Z$. Then, modulo relations (1) and (2), the relation $E_Z\cdot P_{X,Z,Y}(-\sigma_Y)$ is implied by the type (4) relation $E_{Y_i}\cdot P_{X,Y_i,W\cap Y_i}(-\sigma_{W\cap Y_i})$ associated to $W\cap Y_i\subsetneq Y_i$.
\end{enumerate} 
\end{lemma}

\begin{proof}
Set $f := E_Z \cdot P_{X,Z,Y}(-\sigma_Y)$. We first make an observation that will be used in the proofs of both parts. We claim that $E_Z\cdot E_{Y'}$ lies in the ideal generated by relations (1) for all $Y'\in\mathcal{G}$ such that $Y'\subseteq W$ and $Y'\not\subseteq Y$. Indeed, consider such a $Y'$, and suppose $\mathcal{T}:=\{Y'\}\cup\operatorname{Fact}(\mathcal{G},Z)$ is $\mathcal{G}$-nested and induced by a flag $Z_1\subseteq\cdots\subseteq Z_r$ of elements of $\mathcal{A}$. Then $Z_1=\bigcap_{T\in\mathcal{T}}T=Y'\cap Z \subseteq W\cap Z=Y$, so some $\mathcal{G}$-factor $T$ of $Z_1$ is contained in $Y$, as $Y\in\mathcal{G}$. But $T\in\mathcal{T}$ and $Y'\not\subseteq Y$, while every $F\in\operatorname{Fact}(\mathcal{G},Z)$ contains $Z\supsetneq Y$. Hence $\mathcal{T}$ is not $\mathcal{G}$-nested, proving the claim.
We also observe that, modulo relations (2), $E_Z\cdot P_{X,Z,Y}(-\sigma_Y)$ is independent of the choice of $P_{X,Z,Y}(t)$. Indeed, the difference of two choices has constant term in $J_Z$ and all other coefficients in $J_Y$. Any discrepancy in the constant term is annihilated by $E_Z$, while $J_Y\subseteq J_{Y'}$ for every $Y'\subseteq Y$, so any discrepancy in the positive-degree terms vanishes after evaluation at $-\sigma_Y$.

To show \ref{item.lemma.local.relations.4.i}, suppose that $W\in\mathcal{G}$. Since $Y\subseteq W$, we may write
\[
    f = E_Z\cdot P_{X,Z,Y}(-\sigma_W) + E_Z\cdot\underbrace{\raisebox{-1.5mm}{$\Biggl($}\mathlarger{\sum}\nolimits_{\substack{Y'\in\mathcal{G} \\ Y'\subseteq W \\ Y'\not\subseteq Y}} E_{Y'}\raisebox{-1.5mm}{$\Biggr)$}}_{\sigma_W-\sigma_Y}\cdot Q
\]
for some residual polynomial $Q$. The latter summand lies in the ideal generated by relations (1), 
since the claim applies to the monomials $E_Z\cdot E_{Y'}$. 
Moreover, by Lemma~\ref{lemma.chern.polynomial.identities}\ref{lemma.chern.polynomial.identities.ii}  applied to the transversal intersection $Y=W\cap Z$, the polynomial $P_{X,W}(t)$ is an admissible choice of $P_{X,Z,Y}(t)$. 
By the independence just established and the initial claim, $f\equiv E_Z\cdot P_{X,W}(-\sigma_Y)
 \equiv E_Z\cdot P_{X,W}(-\sigma_W)$ modulo relations (1) and (2). The conclusion follows since $P_{X,W}(-\sigma_W)$ is a relation of type (3).

To show \ref{item.lemma.local.relations.4.ii}, suppose that $W\cap Y_i\in\mathcal{G}$ for some $\mathcal{G}$-factor $Y_i$ of $Z$. Since $Z\subseteq Y_i$, we have $Y\subseteq W\cap Y_i\subsetneq Y_i$. The polynomial $P_{X,Y_i,W\cap Y_i}(t)$ is an admissible choice of $P_{X,Z,Y}(t)$. Indeed, this is immediate if $Z=Y_i$. If $Z\subsetneq Y_i$, then at each point of $Y$, transversality of $Z$ and $W$ in $X$ expresses every tangent vector to $Y_i$ as a sum of tangent vectors to $Z$ and $W$. Since $Z\subseteq Y_i$, the latter vector is also tangent to $Y_i$, hence to $W\cap Y_i$ by Lemma~\ref{tangent.intersection}. Thus $Z$ and $W\cap Y_i$ intersect transversally inside $Y_i$, and the assertion follows from Lemma~\ref{lemma.chern.polynomial.identities}\ref{lemma.chern.polynomial.identities.ii}. 
Every $Y'\in\mathcal{G}$ such that $Y'\subseteq W\cap Y_i$ and $Y'\not\subseteq Y$ is contained in $W$. Thus, by the independence of choices and the initial claim, $f\equiv E_Z\cdot P_{X,Y_i,W\cap Y_i}(-\sigma_Y) \equiv E_Z\cdot P_{X,Y_i,W\cap Y_i}(-\sigma_{W\cap Y_i})$ modulo relations (1) and (2). 
Since $Y_i$ is a $\mathcal{G}$-factor of $Z$, the right-hand side is a multiple of $E_{Y_i}\cdot P_{X,Y_i,W\cap Y_i}(-\sigma_{W\cap Y_i})$, proving the result.
\end{proof}

The following two propositions provide criteria for simplifying relations (4) in Theorem~\ref{theorem.main.intro}. Both are immediate consequences of the preceding lemma. We apply these propositions below to recover some Chow ring presentations appearing in the literature.

\begin{proposition}\label{prop.transversally.regular.relations.4}
    In the setting of Theorem~\ref{theorem.main.intro}, if the building set $\mathcal{G}$ is transversally regular, then relations (4) are implied by relations (1), (2) and (3).  
\end{proposition}

\begin{proposition}\label{prop.weak.transversally.regular.relations.4}
    In the setting of Theorem~\ref{theorem.main.intro}, if the building set $\mathcal{G}$ admits weak transversal witnesses, then, in the presence of relations (1)--(3), relations (4) can be replaced by the following analogue (4'), involving only the building set: 
    \begin{itemize}
        \item[(4')] \ \ $E_Z \cdot P_{X,Z,Y}\raisebox{-1.5mm}{$\Biggl($}
    -\mathlarger{\sum}\nolimits_{\substack{Y'\in\mathcal{G}\\ Y'\subseteq Y}} E_{Y'}\raisebox{-1.5mm}{$\Biggr)$}$ for all $Z,Y\in\mathcal{G}$ such that $Y\subsetneq Z$.
    \end{itemize}
\end{proposition}

\smallskip

\begin{example} \label{example.table.of.classes}
    The following table records which of the conditions in
    Definition~\ref{definition.complete.regular.retractable} hold for some well-known classes of
    wonderful compactifications.
    Here Y means that the property holds throughout the family, whereas N means that it does not hold in general. 

    \begin{center}
    \begin{tabular}{|l|c|c|c|}
        \hline 
        \multicolumn{1}{|c|}{Family} & Complete & Retractable & Transversally regular \\  
        \hline
        Fulton--MacPherson spaces $X^{[n]}$ & N & Y & Y \\
        Weighted Fulton--MacPherson spaces $X_{\mathbf a}^{[n]}$ & N & Y & N \\
        Ulyanov's compactifications $X\langle n\rangle$ & Y & Y & Y \\
        Chen--Gibney--Krashen spaces $T_{d,n}$ & N & N & Y \\
        Moduli spaces of points in flags $T_{\mathbf w}^{\mathbf a}$ & N & N & N \\ 
        Polystellahedral compactifications $\mathbb P_H^{[\mathbf a]}$ & Y & Y & Y \\ 
        Clader--Damiolini--Li--Ramadas spaces $\overline{\mathcal{L}}^r_n$ & Y & Y & Y \\
        \hline
    \end{tabular}
    \end{center}

    \begin{itemize}[leftmargin=2em]
    \item For a smooth variety $X$, write $\Delta_I:=\{(x_1,\dots,x_n)\in X^n\mid x_i=x_j \text{ for all }i,j\in I\}$ for each $I\subseteq\{1,\dots,n\}$ with $|I|\geq2$. The Fulton--MacPherson space $X^{[n]}$ is the wonderful compactification of the building set consisting of all such diagonals; its arrangement consists of their intersections, called polydiagonals. For a weight vector $\mathbf a=(a_1,\dots,a_n)$ with $0<a_i\leq1$, the weighted Fulton--MacPherson space $X_{\mathbf a}^{[n]}$ is the wonderful compactification of the building set consisting of those $\Delta_I$ with $\sum_{i\in I}a_i>1$; its arrangement is a subcollection of the polydiagonals.  
    These building sets are not complete in general. The building sets of Fulton--MacPherson spaces are transversally regular, whereas those of weighted Fulton--MacPherson spaces in general admit only weak transversal witnesses. The building set of the weighted Fulton--MacPherson space has a retractable arrangement for every weight vector, with the retractions obtained by copying coordinates across the blocks.
    \item Ulyanov's compactifications $X\langle n\rangle$ \cite{ulyanov} are wonderful compactifications whose building sets consist of all polydiagonals in the self-product $X^n$ of a smooth variety $X$. In particular, the building set is complete, retractable and, just as with Fulton--MacPherson spaces, transversally regular.
    \item The space $T_{d,n}$ of \cite{Chen-Gibney-Krashen} compactifies the space of configurations of $n$ distinct labeled points in $\mathbb{A}^d$ up to translation and scaling. 
    The space $T_{d,n}$ is a wonderful compactification of a building set on $\mathbb{P}^{d(n-1)-1}$ constructed as follows. Modulo translation, every point configuration is equivalent to one of the form $p_1,\dots,p_{n-1},0\in\mathbb{A}^d$. If we write $p_i=(x_{i1},\dots,x_{id})$ for all $i=1,\dots,n-1$, then the equivalence class of this point configuration corresponds to 
    \[
        [x_{11} : \cdots : x_{1d}: x_{21} : \cdots : x_{2d} : \cdots : x_{(n-1)1} : \cdots : x_{(n-1)d}]\in\mathbb{P}^{d(n-1)-1}.
    \]
    The space $T_{d,n}$ is the wonderful compactification of the building set consisting of the subvarieties $\delta_I$ with $I\subsetneq\{1,\dots,n\}$ and $|I|\geq 2$, where 
    \begin{equation*}
    \delta_{I} :=\begin{cases} 
    \bigcap_{i,j \in I}V(x_{i1}-x_{j1},\ldots , x_{id}-x_{jd}) & \ \text{if } n\notin I;
    \\
    \bigcap_{i\in I\smallsetminus \{n\} }V( x_{i1} ,\ldots ,x_{id})     & \  \text{if } n\in I.  
    \end{cases}
    \end{equation*}
    This building set is not complete in general. 
    Suitable diagonals $\delta_J$ in this building set provide transversal witnesses, 
    but there is no retraction from projective space onto a positive-dimensional proper linear subspace.

    \item The spaces $T_{\mathbf w}^{\mathbf a}$ of \cite{polymatroids-moduli}, where $\mathbf a=(a_1,\dots,a_n)$ with $a_1\geq\cdots\geq a_n>0$, compactify configurations of $n$ distinct labeled points $(p_1,\dots,p_n)$ in the flag $\mathbb C^{a_n}\subseteq\cdots\subseteq\mathbb C^{a_1}$ with $p_i\in\mathbb C^{a_i}$, up to simultaneous scaling and translations preserving the flag. They depend on a weight vector $\mathbf w=(w_1,\dots,w_n)$ with $0<w_i\leq1$ and $\sum_{i=1}^n w_i>1$. After translating the last point to the origin, the resulting configurations up to scaling lie in $\mathbb P^{a_1+\cdots+a_{n-1}-1}$. The space $T_{\mathbf w}^{\mathbf a}$ is the wonderful compactification of the building set on this projective space consisting of the linear collision loci $\delta_I^{\mathbf a}$ with $I\subsetneq\{1,\dots,n\}$ and $\sum_{i\in I}w_i>1$, where $\delta_I^{\mathbf a}$ is defined by the coincidence of the points indexed by $I$. For $\mathbf a=(d,\dots,d)$ and $\mathbf w=(1,\dots,1)$, these spaces recover $T_{d,n}$. For Losev--Manin weights $\mathbf w=(\varepsilon_1,\dots,\varepsilon_{n-1},1)$ with $\sum_{i=1}^{n-1}\varepsilon_i\leq1$, they are isomorphic to polypermutohedral varieties. These building sets are neither complete nor retractable in general, since the family includes $T_{d,n}$. They need not be transversally regular either: for $\mathbf a=(1,1,1,1)$ and $\mathbf w=(3/5,3/5,1/5,1/5)$, the building set consists of a line and two distinct points on it, so neither point-line pair admits a transversal witness.
    
    \item The spaces $\mathbb P_H^{[\mathbf a]}$ of \cite{polymatroids-moduli}, where $\mathbf a=(a_1,\dots,a_n)$ with $a_1\geq\cdots\geq a_n>0$, compactify configurations of $n$ labeled, not necessarily distinct points $(p_1,\dots,p_n)$ in the flag $\mathbb C^{a_n}\subseteq\cdots\subseteq\mathbb C^{a_1}$ with $p_i \in \mathbb C^{a_i}$, and are isomorphic to polystellahedral varieties. They are wonderful compactifications of the building set on $\prod_{i=1}^n\mathbb P^{a_i}$ consisting of all nonempty intersections of the divisors 
    \[
    D_i:=\mathbb P^{a_1}\times\cdots\times H_i\times\cdots\times\mathbb P^{a_n},
    \]
    where $H_i\subseteq\mathbb P^{a_i}$ are the compatible hyperplanes at infinity. This building set is complete and hence vacuously retractable. It is also transversally regular, since any inclusion between its members is obtained by imposing additional hyperplane conditions in distinct factors.

   \item The moduli space $\overline{\mathcal{L}}^r_n$ parametrizing stable $(r,n)$-curves over $\mathbb{C}$ was introduced in \cite[Definition~2.2]{clader2023permutohedral}.
   By \cite[Theorem~3.4]{clader2023wonderful}, this space is the wonderful compactification of the building set on $\left(\mathbb{P}^1\right)^n$ consisting of all nonempty intersections of the hypersurfaces $\{\widetilde{H}_i^j\}$ defined as follows. Fix a primitive $r$-th root of unity $\zeta$. For each $i\in\{1,\dots,n\}$ and $j\in\{0,1,\dots,r-1\}$, define 
    \[
        \widetilde{H}_i^j=\{(p_1,\dots,p_n)\in(\mathbb{P}^1)^n \,\vert\, p_i=\zeta^j\}.
    \] 
    This building set is complete by definition, and thus also vacuously retractable. It is also transversally regular, as can be seen by applying Lemma~\ref{lemma.codimension.transversality.two.subvarieties} appropriately.
    \end{itemize}
\end{example}

\begin{example}[Weighted Fulton--MacPherson Spaces] \label{example.WFM}
    The Chow rings of weighted Fulton--Mac\-Pher\-son spaces were described in \cite[Theorem~8]{routis2014weighted} and \cite[Theorem~3.2]{petersen2017chow} using four families of relations that in our language are (1'), (2'), (3) and (4'). 
    That presentation is recovered from our Theorem~\ref{theorem.main.intro}. 
    Indeed, the surjectivity assumptions hold because every diagonal or polydiagonal is a retract of $X^n$; 
    (1) is equivalent to (1') by \cite[Theorem~2~(2)]{routis2014weighted}; and
    (2') and (4') can be used in place of (2) and (4) by Propositions~\ref{prop.reduction.relations.2} and \ref{prop.weak.transversally.regular.relations.4}, since here the building set has a retractable arrangement and admits weak transversal witnesses. 
    Moreover, when all weights are equal to $1$, this presentation can be simplified to the presentation of Fulton--MacPherson \cite[Theorem~6]{Fulton-MacPherson}, as shown in \cite[Theorem~3.4]{petersen2017chow}. 
\end{example}

\subsection{On complete and transversally regular building sets}
\label{sec:based.hu.example}
In \cite{hu2003compactification}, Y. Hu studied and computed the Chow rings of compactifications of \emph{regular simple arrangements of subvarieties} on a smooth variety $X$. In our terminology, Hu's compactification is the wonderful compactification $X_\mathcal{G}$ of a complete building set $\mathcal{G}$ satisfying the following condition:

\begin{definition}[\cite{hu2003compactification}]
    We say that a building set $\mathcal{G}$ is \emph{regular} if every element of $\mathcal{G}$ is the intersection of containment-maximal elements of $\mathcal{G}$.
\end{definition}

Regularity is related to transversal regularity:

\begin{lemma}
    If a building set $\mathcal{G}$ is transversally regular, then every $Y\in\mathcal{G}$ is the \emph{transversal} intersection of containment-maximal elements of $\mathcal{G}$. In particular, a transversally regular building set is regular.
\end{lemma}

\begin{proof}
    We use induction on $\operatorname{codim}Y$. If $Y$ is containment-maximal the result is trivial. We may therefore assume that $Y$ is not containment-maximal and choose a containment-maximal $M\in\mathcal{G}$ with $Y\subsetneq M$. By transversal regularity there is some $W\in\mathcal{G}$ with $Y=W\cap M$ transverse, so $\operatorname{codim}W<\operatorname{codim}Y$ and by induction $W=M_2\cap\cdots\cap M_k$ transversally, with each $M_i$ containment-maximal. Since $Y=W\cap M$ and $W=M_2\cap\cdots\cap M_k$ are transversal intersections, we have
    \[
    \codim_X(Y)=\codim_X(M)+\codim_X(W)=\codim_X(M)+\sum_{i=2}^k\codim_X(M_i).
    \]
    The subvarieties $M,M_2,\dots,M_k$ belong to the arrangement induced by $\mathcal{G}$, and their intersection is $Y$. Therefore, Lemma~\ref{lemma.codimension.characterization.transversal.intersections} implies that they intersect transversally on $X$.
\end{proof}

In our language, the presentation in \cite[Theorem~1.8]{hu2003compactification} is stated for a complete and regular building set and uses relations (1'), (2) and (3). 
The surjectivity assumptions in \cite[Theorem~1.8]{hu2003compactification} also ensure the applicability of our Theorem~\ref{theorem.main.intro}. 
Hence, the presentation in \cite[Theorem~1.8]{hu2003compactification} is recovered, for example, when the complete building set is transversally regular, whereas regularity alone does not suffice in general. 
Indeed, for a complete building set the relations (1)--(3) of Theorem~\ref{theorem.main.intro} are precisely the three families of relations in \cite[Theorem~1.8]{hu2003compactification}: relations (1) reduce to the pairwise relations (1'), by Proposition~\ref{prop.complete.relations.1}.
Relations (2) and (3), being indexed by $\mathcal{A}=\mathcal{G}$, coincide with Hu's second and third families.  
If the complete building set is additionally transversally regular, relations (4) become redundant by Proposition~\ref{prop.transversally.regular.relations.4}, and the presentation in \cite[Theorem~1.8]{hu2003compactification} holds  under this additional assumption. 

However, regularity alone does not make relations (4) redundant. The example below exhibits a complete and regular building set for which relations (4) do not follow from (1)--(3), so the presentation in \cite[Theorem~1.8]{hu2003compactification} is missing the relations of this type.
 
\begin{example} \label{example.regular.but.relation.four.needed}
    Let $\mathcal{G}=\{L_1,L_2,p\}$ be the complete and regular building set in $\mathbb{P}^3$ consisting of two lines intersecting at a point $p$. The intersection $p=L_1\cap L_2$ is \emph{not} transversal, since 
    \[
      \codim(p)=3 \neq 4 = \codim(L_1)+\codim(L_2).
    \]
    Write $A^\bullet(\mathbb{P}^3) = \mathbb{Z}[h]/\langle h^4\rangle $ with $h$ the hyperplane class. Theorem~\ref{theorem.main.intro} gives $A^\bullet(X_\mathcal{G})$ as the
    quotient of $A^\bullet(\mathbb{P}^3)[E_{L_1},E_{L_2},E_p]$ by the relations:
    \begin{enumerate}
        \item[(1)] \ \ $E_{L_1}\cdot E_{L_2}$;
        \item[(2)] \ \ $h\cdot E_p$ \ and \ $h^2\cdot E_{L_i}$;
        \item[(3)] \ \ $h^3-E_p^3$ \ and \ $\big(h-E_p-E_{L_i}\big)^2$;
        \item[(4)] \ \ $E_{L_i}\big(h-E_p\big)$,
    \end{enumerate}
    for $i=1,2$. The relations (4) in this presentation are irredundant.
    Indeed, mapping $h, E_{L_1}, E_p, E_{L_2}$ to $x, y, 0,0$, respectively, in $\mathbb{F}_2[x, y] /\langle x^2, y^2\rangle$ annihilates all relations of types (1)-(3), while $E_{L_1}\left(h-E_p\right)$ maps to $x y \neq 0$; by symmetry, the same argument applies to the relation $E_{L_2}\left(h-E_p\right)$. 
\end{example}

\subsection{Sharpness of Theorem~\ref{theorem.main.intro}}
In this subsection we give examples showing that relations (2) and (4) cannot in general be indexed only by $\mathcal{G}$ rather than by $\mathcal{A}$; see Examples~\ref{example.relations.2} and \ref{example.relations.4}, respectively. 
We also show that, although relations (1) may be restricted to containment-minimal subsets $\mathcal{T}\subseteq\mathcal{G}$ that are not $\mathcal{G}$-nested, no proper subfamily of these relations suffices in general; see Example~\ref{example.relations.1}. 
In particular, relations (1') cannot replace relations (1) in full generality.

\begin{example} \label{example.relations.2}     
    Let $X=\mathbb{P}^4$ and consider the building set $\mathcal{G}=\{V,W\}$ on $X$, where $V$ and $W$ are two two-dimensional linear subspaces intersecting at a point. Write $A^\bullet(\mathbb{P}^4) = \mathbb{Z}[h]/\langle h^5\rangle$, with $h$ being the hyperplane class, and consider the Chern polynomials of $V$ and $W$: $P_V(t) = P_W(t) = (t+h)^2 = t^2+2ht+h^2$. We have that $J_V = J_W = \langle h^3\rangle $ and $J_{V\cap W} = \langle h\rangle$. Using our Theorem~\ref{theorem.main.intro} or via a direct computation one sees that the Chow ring of $X_{\mathcal{G}}$ is 
    \[
        A^\bullet(X_{\mathcal{G}}) 
        = \frac{\mathbb{Z}[h,E_V,E_W]}{\left\langle h^5, h^3 E_W, h E_V E_W, h^3 E_V, E_V^2-2E_Vh+h^2, E_W^2-2E_Wh+h^2\right\rangle }.
    \]
     The homomorphism $\phi: \mathbb{Z}[h,E_V,E_W] \rightarrow \mathbb{Z}[h]$ sending $E_V$, $ E_W,h$ to $h$ satisfies $\phi(h E_V E_W)=h^3$,~which is not in the ideal $\langle h^4 \rangle$ generated by the images of the other relations. Hence, the type (2) relation $h E_V E_W \in J_{V\cap W}\cdot E_VE_W$, corresponding to the arrangement element $V \cap W$, is independent of the other relations.
\end{example}

A similar situation occurs for the relations of type (4). 

\begin{example} \label{example.relations.4} 
    Let $X=\mathbb{P}^5$ and consider the building set $\mathcal{G}=\{Z,V,W\}$ on $X$, where $V$ and $W$ are two three-dimensional linear subspaces intersecting along a line, and $Z$ is a point on this line. 
    The associated arrangement is $\mathcal{A}=\mathcal{G}\cup\{V \cap W\}$.
    Write $A^\bullet(X) = \mathbb{Z}[h]/\langle h^6\rangle$ with $h$ being the hyperplane class. 
    From Theorem~\ref{theorem.main.intro}, or via direct computation, it follows that the Chow ring of $X_{\mathcal{G}}$ is
    \[
    A^\bullet\left(X_\mathcal{G}\right) = \frac{A^\bullet(X)[E_Z,E_V,E_W]}{\langle \, \text{relations} \, \rangle},
    \]
    where the relations are: 
    \begin{itemize}
        \item[(1)]: There are no relations of this kind in this example.
        \item[(2)]: $h^4E_V,\quad h^4E_W,\quad hE_Z,\quad h^2E_VE_W$.
        \item[(3)]: $(h-E_Z)^5,\quad (h-E_V-E_Z)^2,\quad (h-E_W-E_Z)^2$.
        \item[(4)]: $E_V(h-E_Z)^3,\quad E_W(h-E_Z)^3,\quad E_VE_W(h-E_Z)$.
    \end{itemize}
    In particular, notice that the last relation $E_VE_W\cdot P_{X,V\cap W,Z}(-E_Z) = E_VE_W(h-E_Z) = 0$ of type (4) is accounting for the inclusion of $Z$ in the element $V\cap W$ of the arrangement. This relation is independent of all the others. Indeed, mapping $h,E_V,E_W,E_Z$ to $x,y,z,0$, respectively, in
$\mathbb{F}_2[x,y,z]/\langle x^2,y^2,z^2\rangle$ annihilates every other relation, while $E_VE_W(h-E_Z)$ maps to $xyz\neq 0$.
\end{example}

For the relations of type (1) we have the following.

\begin{example} \label{example.relations.1}
    Let $X=\mathbb{P}^4$ and consider the building set $\mathcal{G}=\{V_1,V_2,V_3\}$, where $V_1,V_2,V_3$ are three two-dimensional linear subspaces such that each $V_i\cap V_j$ is a point, these three points are distinct, and $V_1\cap V_2\cap V_3=\emptyset$. The associated arrangement is $\mathcal{A}=\mathcal{G}\cup \{V_i\cap V_j \}_{i<j}$. Write $A^\bullet(X)=\mathbb{Z}[h]/\langle h^5\rangle$ with $h$ the hyperplane class, and set $E_i:=E_{V_i}$. We have $J_{V_i} = \langle h^3\rangle$, $J_{V_i\cap V_j} = \langle h\rangle$, and we can take $P_{X,V_i}(t)=(t+h)^2$. Then the relations of Theorem~\ref{theorem.main.intro} are:
    \begin{itemize}
        \item[(1)]: $E_1E_2E_3$.
        \item[(2)]: $h^3E_i$ for $i=1,2,3$, \quad $hE_iE_j$ for $i<j$.
        \item[(3)]: $(h-E_i)^2$ for $i=1,2,3$.
        \item[(4)]: There are no relations of this kind in this example.
    \end{itemize}
    Every pair in $\mathcal{G}$ is $\mathcal{G}$-nested, so $\{V_1,V_2,V_3\}$ is the unique non-nested subset and (1) contributes a single relation, which is independent of the others. Indeed, mapping $h,E_1,E_2,E_3$ to $0,x_1,x_2,x_3$, respectively, in $\mathbb{Z}[x_1,x_2,x_3]/\langle x_1^2,x_2^2,x_3^2\rangle$ annihilates all other relations, whereas $E_1E_2E_3$ maps to $x_1x_2x_3\neq0$.
\end{example}

% ====== Bibliography ======

\bibliographystyle{alpha}
\bibliography{biblioChowWonderful.bib}

\end{document}